\documentclass[12pt, a4paper,reqno]{amsart}
\usepackage{amscd,amssymb,amsmath, array,latexsym,amsbsy}
\allowdisplaybreaks[2]
\begin{document}
%%%%%%%%%%%%%%%%%%%%%%%%%%%%%%%%%%
\def\supp{\operatorname{supp}}
\def\tr{\operatorname{tr}}
\def\rk{\operatorname{rk}}
\def\lt{{\operatorname{{lt}}}}
\def\id{\operatorname{{id}}}
\def\interior{\operatorname{Int}} % interior
\def\Ind{\operatorname{Ind}} % induced representation
\def\Prim{\operatorname{Prim}}
\def\Aut{\operatorname{Aut}}
\def\diag{\operatorname{diag}}
%%%%%%%%%%%%%%%%%%%%%%%%%%%%
\def\H{\mathcal{H}} %Hilbert space
\def\K{\mathcal{K}} % Compacts
\def\N{\mathcal{N}} % neighbourhood base at zero in dual
\def\C{\mathbb{C}}
\def\T{\mathbb{T}}
\def\Z{\mathbb{Z}}
\def\R{\mathbb{R}}
\def\P{\mathbb{P}}
\def\NN{\mathbb{N}}
%%%%%%%%%%%%%%%%%%%%%%%%%%%%%%%%%%
\newtheorem{thm}{Theorem}[section]
\newtheorem{cor}[thm]{Corollary}
\newtheorem{prop}[thm]{Proposition}
\newtheorem{lemma}[thm]{Lemma}
\theoremstyle{definition}
\newtheorem{defn}[thm]{Definition}
\newtheorem{remark}[thm]{Remark}
\newtheorem{question}[thm]{Question}
\newtheorem*{remark1}{Remark on hypotheses}
\newtheorem{example}[thm]{Example}
\numberwithin{equation}{section}
%%%%%%%%%%%%%%%%%%%%%%%%%%%%%%%%%%%
\title
[\boldmath Strength of convergence and multiplicities in the spectrum]
{Strength of convergence and multiplicities in the spectrum of a $C^*$-dynamical system}

\author[Archbold]{Robert Archbold}
\address{Department of Mathematical Sciences
\\University of Aberdeen
\\Aberdeen AB24 3UE
\\Scotland
\\United Kingdom
}
\email{r.archbold@maths.abdn.ac.uk}

\author[an Huef]{Astrid an Huef}
\address{School of Mathematics and Statistics
\\University of New South Wales
\\Sydney, NSW 2052
\\Australia}
\email{astrid@unsw.edu.au}

\date{December 7, 2006, with revisions April 3, 2007}

\begin{abstract}
We consider separable $C^*$-dynamical systems $(A,G,\alpha)$ for
which the induced action of the group $G$ on the   primitive ideal space $\Prim A$    of the $C^*$-algebra $A$ is free. We study how
the representation theory of the associated crossed-product $C^*$-algebra $A\rtimes_\alpha G$ depends  on the
representation theory of $A$ and the properties of the action of $G$
on   $\Prim A$ and the spectrum $\hat A$.    Our main tools involve computations of  upper and lower bounds on multiplicity numbers associated to irreducible representations of $A\rtimes_\alpha G$. We apply our techniques to give necessary and sufficient
conditions, in terms of $A$ and   the action of $G$,     for
$A\rtimes_{\alpha}G$ to be (i) a continuous-trace $C^*$-algebra,
(ii) a Fell $C^*$-algebra and  (iii) a bounded-trace $C^*$-algebra.
When $G$ is amenable, we also give necessary and
sufficient conditions for the crossed-product $C^*$-algebra $A\rtimes_{\alpha}G$ to be (iv) a liminal
$C^*$-algebra and (v) a Type I $C^*$-algebra.  The results in (i), (iii)--(v) extend
some earlier special cases in which $A$ was assumed to have the
corresponding property.
\end{abstract}

\thanks{This research was supported by grants from the Australian Research Council, the University of Aberdeen and the University of New South Wales.}

\maketitle
\section{Introduction}
Throughout, $(A,G,\alpha)$ is a separable $C^*$-dynamical system, so
that $A$ is a separable $C^*$-algebra, $G$ is a second countable
locally compact group and $\alpha:G\to\Aut A$ is a strongly continuous homomorphism into the group of automorphisms of $A$.   There are induced actions of $G$ on the
spectrum $\hat A$ and the primitive ideal space $\Prim A$ of $A$ given by
$s\cdot\sigma=\sigma\circ\alpha_{s^{-1}}$ for $\sigma\in\hat A$ and $s\cdot P=\alpha_s(P)$ for $P\in\Prim A$, so that $s\cdot\ker\sigma =\ker(\sigma\circ\alpha_s^{-1})$.  Thus $(G,\hat A)$ and $(G,\Prim A)$ are
jointly continuous transformation groups (see, for example, \cite[Lemma 7.1]{raewill}).    
The representation theory
of the crossed product $C^*$-algebra $A\rtimes_\alpha G$ associated to $(A,G,\alpha)$
depends  on the representation theory of $A$
and the   properties of the action of $G$,     and this dependence has been widely studied  (see, for example, \cite{tak,green2,GR,RR,GL1, GL2,OR, E,AaH}
and, for $A=C_0(X)$, \cite{Goot,green1,W2,W,aH-Fell, aH}).

In this paper we assume that the induced action of $G$ on   $\Prim A$    is free, and study  the effect on lower multiplicity numbers of  inducing sequences  $\pi_n\to\pi$ in $\hat A$ to sequences $\Ind\pi_n\to\Ind\pi$ in  $(A\rtimes_\alpha G)^\wedge$. In particular, we obtain lower and upper bounds on the  lower multiplicity $M_L(\Ind\pi,(\Ind\pi_n))$ of $\Ind\pi$ relative to $(\Ind\pi_n)$.

The lower bounds for  $M_L(\Ind\pi,(\Ind\pi_n))$ are obtained in two situations: (1)
if  $(\pi_n)$ converges $k$-times  to $\pi$ in ${\hat A}/G$ for some positive
integer $k$ (Theorem~\ref{thm-a}) and
(2) if $(\pi_n)$ converges
not only to $\pi$ but also to some other points in the closure of
the orbit of $\pi$ (Theorem~\ref{thm-1a}).
On the other hand, if the spectrum $\hat A$ is Hausdorff, we obtain an upper bound for
$M_L(\Ind\pi,(\Ind\pi_n))$ in terms of the upper multiplicity $M_U(\pi)$ of $\pi$  in $\hat A$ and
measure-theoretic properties of the action of $G$ on  $\hat A$ (Theorem~\ref{3.5'}).
The proof of Theorem~\ref{3.5'} builds  on the proof  for $A=C_0(X)$
from \cite[Theorem 3.5]{AaH}, but  the new technical details to deal with non-commutative $A$ are substantial.

Using both Theorems~\ref{thm-a} and  \ref{3.5'} we show  in Theorem~\ref{thm-circle2} how the representation theory of $A\rtimes_\alpha G$ is related to that of $C_0({\hat A})\rtimes G$ when $\hat A$ is Hausdorff; this result
provides an extension of \cite[Theorem 1.1]{AaH} to the non-commutative case.
It also follows from Theorem~\ref{thm-a} that the upper multiplicity increases under inducing, that  is $M_U(\pi)\leq M_U(\Ind\pi)$ for all $\pi\in\hat A$ (see Corollary~\ref{cor-A}); this gives a new  proof  of Deicke's theorem \cite[Theorem
5.3.2]{D-thesis} without the use of coactions and with weaker hypotheses.

We apply our results to determine necessary and sufficient
conditions, in terms of $A$ and   the induced free   actions of $G$ on
$\Prim A$ and $\hat A$,    for $A\rtimes_{\alpha}G$ to be
(i) a continuous-trace
$C^*$-algebra, (ii) a Fell $C^*$-algebra, (iii) a bounded-trace
$C^*$-algebra. When $G$ is amenable we also give
necessary and sufficient conditions for $A\rtimes_{\alpha}G$ to be
(iv) a liminal $C^*$-algebra, (v) a Type I $C^*$-algebra. Special cases of  the results (i) and  (iii)--(v) have been
previously obtained under the underlying assumption that $A$ already
has the property that is being considered for $A\rtimes_\alpha G$.
For example, in Theorem~\ref{thm-ct} we generalise a theorem obtained  by the combined efforts of
Raeburn and Rosenberg \cite{RR}, Olesen and Raeburn \cite{OR} and Deicke \cite{D} by showing that $A\rtimes_\alpha G$ has
continuous trace if and only if $A$ has continuous trace and $G$
acts properly on $\hat A$. Apart from methods of proof, our new contribution is that $A$ must have continuous trace if $A\rtimes_\alpha G$ has.
 Of course, the freeness of the action of
$G$ on   $\Prim A$
   is a crucial underlying condition:
by Takai duality  the
dual action $\hat\alpha$ of $G=\T$ on   the non-Type I   
irrational rotation algebra $A=C(\T)\rtimes_\alpha \Z$ gives  a crossed product $A\rtimes_{\hat\alpha} G
= C(\T)\otimes \mathcal{K}$ which has continuous trace.
Finally, Example~\ref{ex-fell} illustrates Theorem~\ref{thm-1a} and \S4 is devoted to an example illustrating the use of Theorems~\ref{thm-a} and \ref{3.5'}.

 \subsection*{Notation and preliminaries}
Let $(A,G,\alpha)$ be a separable dynamical system and $\pi\colon A\to B(\H_\pi)$ be  a representation. (All our representations are non-degenerate.)  The group $G$ comes equipped with a left Haar measure $\mu$ and modular function $\Delta$, and also a right Haar measure $\nu$ given by $\nu(E)=\mu(E^{-1})$.

Define $\tilde\pi\colon A\to B(L^2(G,\H_\pi,\nu))$ and $\lambda\colon G\to U(L^2(G,\H_\pi,\nu))$ by
\begin{equation*}
(\tilde\pi(a)\xi)(s)=\pi(\alpha_{s^{-1}}(a))(\xi(s))\quad\text{and}\quad (\lambda_t\xi)(s)=\Delta(t)^{1/2}\xi(t^{-1}s)
\end{equation*}
for $a\in A$, $s,t\in G$ and $\xi\in L^2(G,\H_\pi,\nu)$.  Then $(\tilde\pi,\lambda)$ is covariant for $(A,G,\alpha)$ and  we write $\Ind\pi$ for $\tilde\pi\rtimes\lambda:A\rtimes_\alpha G\to B(L^2(G,\H_\pi,\nu))$. For $f\in C_c(G,A)$ we have
\begin{align*}
(\Ind\pi(f)\xi)(s)&=\left(\int_G\tilde\pi(f(t))\lambda_t\,d\mu(t)\xi\right)(s)\\
&=\int_G(\tilde\pi(f(t)))(\lambda_t\xi)(s)\,d\mu(t)\\
&=\int_G\pi(\alpha_{s^{-1}}(f(t)))\xi(t^{-1}s)\Delta(t)^{1/2}\, d\mu(t)\\
&=\int_G\pi(\alpha_{s^{-1}}(f(t^{-1}))\xi(ts)\Delta(t)^{-1/2}\, d\nu(t).
\end{align*}

  Suppose that the induced action of $G$ on $\Prim A$ is free. If $\pi\in\hat A$ then $\pi$ is a homogeneous representation and $\Ind\pi$ is irreducible  by separability and freeness of the induced action on $\Prim A$ \cite[Proposition~1.7]{EW}. We thank the referee for pointing out to us that freeness on $\hat A$ is not known to be sufficient for \cite[Proposition~1.7]{EW}.

   We will make frequent use of the following result which follows from \cite[Theorem~24]{green2} (see  \cite[\S 4.3]{D-thesis} and  \cite[Chapter~1]{E} for a discussion of this).

%%%%%%%%%%%%%%%%%%%%%%%%%%%%%%%%%%%%%%%%%%%%%%%%%%%%%%%%%%%%%%%%%%%%%%%%%%%%%%%%%%%%%%%%%%%%%%%%%%%%%%%%%%%%%%%%%%%%%%%%%%%%%%%%%%%%%%%%%

\begin{thm} \textup{(Green)}\label{thm-green}
Let $(A,G,\alpha)$ be a separable $C^*$-dynamical system such that the induced action of $G$ on $\hat A$ is free. If,  for each $\pi\in \hat A$,
\begin{enumerate}
\item  the orbit $G\cdot\pi$ is locally closed in $\hat A$; and
\item the map $s\mapsto s\cdot\pi$ is a homeomorphism of $G$ onto $G\cdot\pi$,
\end{enumerate}
then   $A$ and $A\rtimes_\alpha G$ are Type \textrm{I}    and the map $\Ind: \hat A\to (A\rtimes_\alpha G)^\wedge$ induces a homeomorphism of $\hat A/G$ onto $(A\rtimes_\alpha G)^\wedge$.
\end{thm}

\begin{remark}\label{rem-hyp}
 
If the action of $G$ on $\Prim A$ is free then  the action of $G$ on $\hat A$ is free.
On the other hand, 
%if the action of  $G$ on $\hat A$ is free and conditions (1) and  (2) of Theorem~\ref{thm-green} hold, then $\hat A\simeq \Prim(A)$ %(\cite[Remark~4.3.2 (ii)]{D-thesis}). To see this, let $\pi\in \hat A$. Since $G\simeq G\cdot\pi$, the singleton $\{\pi\}$ is closed (and hence %locally closed) in $G\cdot\pi$. Since $G\cdot\pi$ is  locally closed in $\hat A$ it follows that $\{\pi\}$ is  also locally closed in $\hat A$.  %Thus $\hat A$ is $\textrm{T}_0$ and hence $\hat A\simeq \Prim A$ via the kernel map (and $A$ is Type \textrm{I} by separability).
 if $G$ acts freely on $\hat A$ and condition (2) of Theorem~\ref{thm-green} holds, then $G$ acts freely on $\Prim A$ as well.  To see the latter, suppose that  $\ker\pi =s\cdot\ker\pi$ for some $s\neq e$ and $\pi\in\hat A$. Then $s\cdot\pi$ and $\pi$ are distinct by freeness of the action on $\hat A$. Since $s\cdot \ker\pi=\ker(s\cdot\pi)$ every open neighbourhood of $\pi$ in $\hat A$ must contain $s\cdot\pi$ as well.  But  condition (2) implies there are open neighbourhoods $U$ and $V$ of $\pi$ and $s\cdot\pi$, respectively, such that $(U\cap G\cdot\pi)\cap (V\cap G\cdot\pi)=\emptyset$, a contradiction.
 
 In particular, if the action of  $G$ on $\hat A$ is free and conditions (1) and  (2) of Theorem~\ref{thm-green} hold, then $\hat A\simeq\Prim(A)$ via the kernel map and $A$ is Type \textrm{I} by separability (see, for example, \cite[Remark~4.3.2 (ii)]{D-thesis}).
\end{remark}

We write $S(A)$ and
$P(A)$ for the state space and the set of pure states of a
$C^*$-algebra $A$, respectively.
We refer the reader to \cite{A} for the definitions of the upper and lower multiplicities $M_U(\pi)$ and  $M_L(\pi)$ of an irreducible representation $\pi$ of $A$, and to \cite{AS} for the definitions of upper and lower multiplicities $M_U(\pi,(\pi_n))$ and $M_L(\pi,(\pi_n))$ of $\pi$ relative to a net $(\pi_n)$ in $\hat A$. We set $\P=\mathbb N\setminus\{0\}$.

%%%%%%%%%%%%%%%%%%%%%%%%%%%%%%%%%%%%%%%%%%%%%%%%%%%%%%%%%%%%%%%%%%%%%%%%%%%%%%%%%%%%%%%%%%%%%%%%%%%%%%%%%%%%%%%%%%%%%%%%%%%%%%

\section{Strength of convergence in the spectrum and multiplicities of irreducible representations of $A\rtimes_\alpha G$}

%%%%%%%%%%%%%%%%%%%%%%%%%%%%%%%%%%%%%%%%%%%%%%%%%%%%%%%%%%%%%%%%%%%%%%%%%%%%%%%%%%%%%%%%%%%%%%%%%%%%%%%%%%%%%%%%%%%%%%%%%%%%%%%
Our first theorem is a generalisation and improvement of
\cite[Theorem~2.3]{AD}; it shows that two ingredients contribute to
 the relative lower multiplicity
$M_L(\Ind\pi,(\Ind\pi_n))$ of an irreducible representation
$\Ind\pi$ of $A\rtimes_\alpha G$ relative to the net $(\Ind\pi_n)$.
The first ingredient is the possible $k$-times convergence of the sequence
$(\pi_n)$ to $\pi$ and the second  is the consideration of the $k$ relative lower multiplicities associated to
$\pi$ arising from the $k$-times convergence.

Before stating the theorem, we recall the definition of $k$-times
convergence from \cite[Definition~2.2]{AD}. Let $(G,X)$ be a second countable, locally compact
transformation group where $G$ (but not necessarily $X$) is
Hausdorff, and let $k\in\P$.
 A sequence $(x_n)_{n\geq 1}$ in $X$ is
\emph{$k$-times convergent in $X/G$ to $z\in X$} if there exist $k$
sequences $(t_n^{(1)})_n,(t_n^{(2)})_n,\cdots ,(t_n^{(k)})_n\subset
G$, such that
\begin{enumerate}
\item $t_n^{(i)}\cdot x_n\to z$ as $n\to\infty$ for $1\leq i\leq k$, and
\item if $1\leq i<j\leq k$ then $t_n^{(j)}(t_n^{(i)})^{-1}\to\infty$ as
$n\to\infty$.
\end{enumerate}
%%%%%%%%%%%%%%%%%%%%%%%%%%%%%%%%%%%%%%%%%%%%%%%%%%%%%%%%%%%%%%%%%%%%%%%%%%%%%%%%%%%%%%%%%%%%%%%%%%%%%%%%%%%%%%%%%%%%%%%%%%%%%%%
\begin{thm}\label{thm-a}
Let $(A, G,\alpha)$ be a separable $C^*$-dynamical system such that
the induced action of  $G$ on   $\Prim A$    is free.  Let $\pi\in\hat A$,
let $k\in\P$ and suppose that there is a sequence $(\pi_n)_n$ in
$\hat A\setminus\{\pi\}$ which is $k$-times convergent in $\hat A/G$
to $\pi$ with the following additional data, where $m_1,\dots,
m_k\in\P$: there are sequences
$(t_n^{(1)})_n,(t_n^{(2)})_n,\dots,(t_n^{(k)})_n$ in $G$ such that
\begin{enumerate}
\item $t_n^{(j)}\cdot \pi_n\to \pi$ as $n\to\infty$ for $1\leq j\leq k$; and
\item if $1\leq i<j\leq k$ then $t_n^{(j)}(t_n^{(i)})^{-1}\to\infty$ as
$n\to\infty$; and
\item
%$t_n^{(j)}\cdot\pi_n\in \hat A\setminus\{\pi\}$ eventually and
$M_L(\pi, (t_n^{(j)}\cdot\pi_n))\geq m_j$ for $1\leq j\leq k$.
\end{enumerate}
Then $M_L(\Ind\pi, (\Ind\pi_n))\geq m_1+\cdots +m_k$.
\end{thm}
%%%%%%%%%%%%%%%%%%%%%%%%%%%%%%%%%%%%%%%%%%%%%%%%%%%%%%%%%%%%%%%%%%%%%%%%%%%%%%%%%%%%%%%%%%%%%%%%%%%%%%%%%%%%%%%%%%%%%%%%%%%%%%%
To recover \cite[Theorem~2.3]{AD}  from Theorem~\ref{thm-a}   take $m_1=\cdots=m_k=1$ so that
\[
M_U(\Ind\pi)\geq M_U(\Ind\pi,(\Ind\pi_n))\geq M_L(\Ind\pi,(\Ind\pi_n))\geq k.
\]
To prove Theorem~\ref{thm-a} (and Theorem~\ref{thm-1a} below) we require the following sequence version of \cite[Lemma~5.2]{ASS}.

%%%%%%%%%%%%%%%%%%%%%%%%%%%%%%%%%%%%%%%%%%%%%%%%%%%%%%%%%%%%%%%%%%%%%%%%%%%%%%%%%%%%%%%%%%%%%%%%%%%%%%%%%%%%%%%%%%%%%%%%%%%%%%%
\begin{lemma}\label{lem-j} Let $A$ be a $C^*$-algebra, $\pi\in\hat A$ and $\phi$ a pure state associated with $\pi$, and $(\pi_n)_n$ a sequence
in $\hat A$.
\begin{enumerate}
\item  Suppose that $A$ is separable and that $M_L(\pi,(\pi_n))\geq m$ for some $m\in\P$. Then there exists a subsequence $(\pi_{n_j})_j$ and, for each $j\geq 1$, an orthonormal set $\{\xi_j^{(1)},\cdots, \xi_j^{(m)}\}$ in $\H_{\pi_{n_j}}$ such that
\[
\lim_{j\to\infty}\langle\pi_{n_j}(\cdot)\xi_j^{(i)}\,,\, \xi_j^{(i)}\rangle=\phi\quad\quad(1\leq i\leq m).
\]
\item Suppose that for every subsequence $(\pi_{n_j})_j$  of $(\pi_n)_n$ there is a further subsequence $(\pi_{{n_j}_k})_k$ of $(\pi_{n_j})_j$ such that, for each $k\geq 1$, there is an orthonormal set $\{\xi_k^{(1)},\cdots, \xi_k^{(m)}\}$ in the Hilbert space of $\pi_{{n_j}_k}$ such that
\begin{equation}\label{lem-j-condition}
\lim_{k\to\infty}\langle\pi_{{n_j}_k}(\cdot)\xi_k^{(i)}\,,\, \xi_k^{(i)}\rangle=\phi\quad\quad(1\leq i\leq m).
\end{equation}
Then $M_L(\pi,(\pi_n))\geq m$.
\end{enumerate}
\end{lemma}
%%%%%%%%%%%%%%%%%%%%%%%%%%%%%%%%%%%%%%%%%%%%%%%%%%%%%%%%%%%%%%%%%%%%%%%%%%%%%%%%%%%%%%%%%%%%%%%%%%%%%%%%%%%%%%%%%%%%%%%%%%%%%%%
\begin{proof}
Let $\N$ be the standard base of neighbourhoods of $0$ in $A^*$.

(1) Since $A$ is separable there is  a decreasing basic sequence $(N_j)_j$ of w*-neighbourhoods of $0$ in $A^*$ such that $N_j\in\N$ for all $j$.  Since $M_L(\pi,(\pi_n))\geq m$, for each $N\in\N$ we have
\begin{equation}\label{eq-lem-j}
\liminf_n d(\pi_n,\phi,N)\geq m.
\end{equation}
Applying \eqref{eq-lem-j} with $N=N_1,N_2,\dots$ in turn, we can construct an increasing sequence $n_1<n_2<\cdots$ such that for each $j$ there is an orthonormal set $\{\xi_j^{(1)},\cdots, \xi_j^{(m)}\}$ in $\H_{\pi_{n_j}}$ such that
\[\langle\pi_{n_j}(\cdot)\xi_j^{(i)}\,,\,\xi_j^{(i)}\rangle\in\phi+N_j\quad\quad (1\leq i\leq m).
\]
Since $(N_j)_j$ is decreasing,
\[
\lim_{j\to\infty}\langle\pi_{n_j}(\cdot)\xi_j^{(i)}\,,\,\xi_j^{(i)}\rangle=\phi.
\]

(2) Suppose that $M_L(\pi,(\pi_n))=R<m$.  Then there exists $N\in\N$ such that
$
\liminf_n d(\pi_n,\phi, N)=R.
$
So there exists a subsequence $(\pi_{n_j})_j$ of $(\pi_n)_n$ such that
\begin{equation}\label{eq-lem-j2}
d(\pi_{n_j},\phi,N)=R\quad\text{for all $j\geq 1$}.
\end{equation}
By hypothesis there exists a further subsequence $(\pi_{{n_j}_k})_k$, and, for each $k\geq 1$, an orthonormal set $\{\xi_k^{(1)},\cdots, \xi_k^{(m)}\}$ in the Hilbert space of $\pi_{{n_j}_k}$ such that
\[
\lim_{k\to\infty}\langle\pi_{{n_j}_k}(\cdot)\xi_k^{(i)}\,,\, \xi_k^{(i)}\rangle=\phi\quad\quad(1\leq i\leq m).
\]
Since $m$ is finite, there exists $k_0$ such that for $k\geq k_0$ and $1\leq i\leq m$,
\[
\langle\pi_{{n_j}_k}(\cdot)\xi_k^{(i)}\,,\, \xi_k^{(i)}\rangle\in\phi+N.
\]
So for $k\geq k_0$, $d(\pi_{{n_j}_k},\phi,N)\geq m$, contradicting \eqref{eq-lem-j2}. Thus $M_L(\pi,(\pi_n))\geq m$.
\end{proof}
Passing to a subnet increases the lower multiplicity, that is, $M_L(\pi,(\pi_{{n_j}_k}))\geq M_L(\pi,(\pi_{n_j}))$. So if $A$ is separable, Lemma~\ref{lem-j} says that $M_L(\pi,(\pi_n))\geq m$ if and only if for every subsequence $(\pi_{n_j})_j$  of $(\pi_n)_n$ there is a further subsequence $(\pi_{{n_j}_k})_k$ of $(\pi_{n_j})_j$  such that the $m$ vector condition \eqref{lem-j-condition} holds.
%%%%%%%%%%%%%%%%%%%%%%%%%%%%%%%%%%%%%%%%%%%%%%%%%%%%%%%%%%%%%%%%%%%%%%%%%%%%%%%%%%%%%%%%%%%%%%%%%%%%%%%%%%%%%%%%%%%%%%%%%%%%%%%

\begin{proof}[Proof of Theorem~\ref{thm-a}]
Since $(A,G,\alpha)$ is separable and the induced action on   $\Prim A$    is free, if $\sigma\in\hat A$ then $\Ind\sigma\in (A\rtimes_\alpha G)^\wedge$ by \cite[Proposition~1.7]{EW}.

The proof of the theorem builds on the ideas used to prove \cite[Theorem~2.3]{AD}. Let $\xi\in\H_\pi$ with $\|\xi\|=1$ and let $\phi=\langle\pi(\cdot)\xi\,,\,\xi\rangle$.  Let $W$ be a compact symmetric neighbourhood of $e$ in $G$ and set
\[
\eta=\nu(W)^{-1/2}\chi_W(\cdot)\xi\quad\text{and}\quad\psi=\langle\Ind\pi(\cdot)\eta\,,\,\eta\rangle.
\]
Then $\eta$ is a unit vector in $L^2(G,\H_\pi,\nu)$ and $\psi$ is a pure state of $A\rtimes_\alpha G$ associated with the irreducible representation $\Ind\pi$.  Let $\N$ be the standard w*-neighbourhod base of $0$ in $(A\rtimes_\alpha G)^*$.

Suppose $M_L(\Ind\pi,(\Ind\pi_n))=R<m_1+\cdots +m_k$. There exists $N\in\N$ such that $\liminf_n(\Ind\pi_n,\psi,N)=R$.   Note that the data in (1), (2) and (3) is preserved by passing to subsequences. By replacing $(\pi_n)_n$ by a subsequence and replacing the sequences $(t_n^{(j)})_n\ (1\leq j\leq k)$ by the corresponding subsequences as well, we may assume that
\begin{equation}\label{eq-contradiction}
d(\Ind\pi_n,\psi,N)=R\quad(n\geq 1)
\end{equation}
and  items (1), (2) and (3) still hold.

By (3) (with $j=1$) and Lemma~\ref{lem-j} there exist a strictly increasing sequence $n_1<n_2<\cdots<n_r<\cdots$ and for each $r\geq 1$ an orthonormal set
\[
\{  \xi_{n_r}^{(1,1)},\dots,\xi_{n_r}^{(1,m_1)} \}\subset \H_{t_{n_r}^{(1)}\cdot\pi_{n_r}}(=\H_{\pi_{n_r}})
\] of $m_1$ vectors  such that
\[
\lim_{r\to\infty}\langle t_{n_r}^{(1)}\cdot\pi_{n_r}(\cdot)\xi_{n_r}^{(1,i_1)}\,,\,\xi_{n_r}^{(1,i_1)}\rangle=\phi\quad\quad(1\leq i_1\leq m_1).
\]
Noting again that the data in (1), (2) and (3) is preserved by passing to subsequences, we may apply Lemma~\ref{lem-j} to (3) $k-1$ times in turn so that (after reindexing) we have, for each $j\in\{1,\dots,k\}$, an orthonormal set
\[
\{\xi_n^{(j,1)},\dots,\xi_n^{(j,m_j)}\}\subset \H_{t_n^{(j)}\cdot\pi_n}(=\H_{\pi_n})
\]
of $m_j$ vectors such that
\begin{equation}\label{eq-lem-j-1}
\lim_{n\to\infty}\langle t_n^{(j)}\cdot\pi_n(\cdot)\xi_n^{(j,i_j)}\,,\,\xi_n^{(j,i_j)} \rangle=\phi\quad\quad(1\leq i_j\leq m_j).
\end{equation}
For $n\geq 1$, $1\leq j\leq k$ and $1\leq i_j\leq m_j$, we define $\eta_n^{(j,i_j)}\in L^2(G,\H_{\pi_n},\nu)$ by
\[
\eta_n^{(j,i_j)}(s)=\nu(W)^{-1/2}\chi_{W}(s (t_n^{(j)})^{-1})\xi_n^{(j,i_j)}.
\]
 Note that for fixed $n$ and $j$ and for $p\neq q$ in $\{1,\dots,m_j\}$, we have $\langle \eta_n^{(j,p)}\,,\,\eta_n^{(j,q)} \rangle=0$ because $\langle \xi_n^{(j,p)}\,,\,\xi_n^{(j,q)} \rangle=0$.  By (2) there exists $n_0$ such that for all $n>n_0$ and all $j\neq j'$ in  $\{1,\dots,k\}$ we have  $t_n^{(j)}t_n^{(j')^{-1}}\notin W^2$, and hence $\langle \eta_n^{j,p)}\,,\,\eta_n^{j',q)} \rangle=0$  for $1\leq p\leq m_j$ and $1\leq q\leq m_{j'}$.  So for each $n\geq n_0$,
\begin{equation*}
\{\eta_n^{(j,i_j)}\colon 1\leq j\leq k, 1\leq i_j\leq m_j\}
\end{equation*}
is an orthonormal set of $m_1+\cdots+m_k$ vectors in $L^2(G,\H_{\pi_n},\nu)$.

Let $f\in C_c(G,A)$. Then, for $1\leq j\leq k$ and $1\leq p\leq m_j$,
\begin{align*}
\langle& \Ind\pi_n(f)\eta_n^{(j,p)}\,,\,\eta_n^{(j,p)} \rangle
\notag\\
&=\int_G \langle \Ind\pi_n(f)\eta_n^{(j,p)}(v)\,,\,\eta_n^{(j,p)}(v) \rangle\, d\,\nu(v)
\notag\\
&=\int_G \Big\langle\int_G\pi_n(\alpha_{v^{-1}}(f(u^{-1})))(\eta_n^{(j,p)}(uv))\Delta(u)^{-1/2}\, d\, \nu(u)\,,\, \eta_n^{(j,p)}(v) \Big\rangle\, d\,\nu(v)
\notag\\
&=\nu(W)^{-1}\int_G\int_G\langle \chi_W(uv(t_n^{(j)})^{-1})\Delta(u)^{-1/2}\pi_n(\alpha_{v^{-1}}(f(u^{-1})))\xi_n^{(j,p)}\,,\,
\notag\\ &\hskip 8cm
\chi_W(v(t_n^{(j)})^{-1})\xi_n^{(j,p)}\rangle\, d\, \nu(u)d\,\nu(v)
\notag\\
&=\nu(W)^{-1}\int_{t\in W}\int_{s\in W}\Delta(st^{-1})^{-1/2}\langle \pi_n(\alpha_{(t_n^{(j)})^{-1}t^{-1}}(f(ts^{-1})))\xi_n^{(j,p)}\,,\,
\xi_n^{(j,p)}\rangle\, d\, \nu(s)d\,\nu(t)
\notag\\
&=\nu(W)^{-1}\int_{t\in W}\int_{s\in W}\Delta(st^{-1})^{-1/2}\langle t_n^{(j)}\cdot\pi_n(\alpha_{t^{-1}}(f(ts^{-1})))\xi_n^{(j,p)}\,,\,
\xi_n^{(j,p)}\rangle\, d\, \nu(s)d\,\nu(t),
\end{align*}
by two changes of  variables: first set $s=uv(t_n^{(j)})^{-1}$ and then $t=v(t_n^{(j)})^{-1}$.
On the other hand,
\begin{align*}
\psi(f)
&=\langle\Ind\pi(f)\eta\,,\,\eta\rangle\\
&=\int_G\langle \Ind\pi(f)\eta(t)\,,\,\eta(t) \rangle\, d\nu(t)\\
&=\int_G\Big\langle\int_G\pi(\alpha_{t^{-1}}(f(u^{-1})))(\eta(ut))\Delta(u)^{-1/2}\,d\nu(u)\,,\,\eta(t)\Big\rangle\, d\nu(t)\\
&=\nu(W)^{-1}\int_G\int_G\langle\chi_W(ut)\Delta(u)^{-1/2}\pi(\alpha_{t^{-1}}(f(u^{-1})))\xi\,,\,\chi_W(t)\xi\rangle\,d\nu(u)\, d\nu(t)\\
&=\nu(W)^{-1}\int_{t\in W}\int_{s\in W}\Delta(st^{-1})^{-1/2}\langle \pi(\alpha_{t^{-1}}(f(ts^{-1})))\xi\,,\,\xi \rangle\,d\nu(s)\, d\nu(t).
\end{align*}
So by \eqref{eq-lem-j-1} and the Bounded Convergence Theorem, we conclude that
\[\langle \Ind\pi_n(f)\eta_n^{(j,p)}\,,\,\eta_n^{(j,p)} \rangle\to \langle\Ind\pi(f)\eta\,,\,\eta\rangle=\psi(f).
\]  Since $C_c(G,A)$ is norm-dense in $A\rtimes_\alpha G$ it follows that
\[\langle \Ind\pi_n(\cdot)\eta_n^{(j,p)}\,,\,\eta_n^{(j,p)} \rangle\to \langle\Ind\pi(\cdot)\eta\,,\,\eta\rangle=\psi.
\]
Hence there exists $n\geq n_0$ such that, for all $1\leq j\leq k$ and all $1\leq p\leq m_j$,
\[
\langle \Ind\pi_n(\cdot)\eta_n^{(j,p)}\,,\, \eta_n^{(j,p)}\rangle\in\psi+N.
\]
Thus $d(\Ind\pi_n,\psi,N)\geq m_1+\cdots+m_k$, contradicting \eqref{eq-contradiction}.
\end{proof}

%%%%%%%%%%%%%%%%%%%%%%%%%%%%%%%%%%%%%%%%%%%%%%%%%%%%%%%%%%%%%%%%%%%%%%%%%%%%%%%%%%%%%%%%%%%%%%%%%%%%%%%%%%%%%%%%%%%%%%%%%%%%%%%%%%

If $X$ is any topological space with a jointly continuous action of
$G$ on $X$, we say the transformation group $(G,X)$ is \emph{Cartan}
if every point $x\in X$ has a \emph{wandering neighbourhood}, that
is a neighbourhood $U$ of $x$ such that $\{s\in G:s\cdot U\cap U\neq
\emptyset\}$ is relatively compact in $G$ (cf.
\cite[Definition~1.1.2]{palais} where $X$ is assumed to be
completely regular).
%%%%%%%%%%%%%%%%%%%%%%%%%%%%%%%%%%%%%%%%%%%%%%%%%%%%%%%%%%%%%%%%%%%%%%%%%%%%%%%%%%%%%%%%%%%%%%%%%%%%%%%%%%%%%%%

Some of the ideas of the next lemma have already appeared in the literature, but without the terminology of ``$2$-times convergence''.  For example, Green proved in \cite{green1}  that $C_0(X)\rtimes G$ has continuous trace if and only if the action of $G$ on $X$ is proper (for free actions); the key ingredient in one direction is that if the action is not proper but $X/G$ is Hausdorff  then there is a sequence converging $2$-times.  This same idea is also exploited in \cite{W,aH} for transformation groups and in \cite{MW} for groupoids, for example.

%%%%%%%%%%%%%%%%%%%%%%%%%%%%%%%%%%%%%%%%%%%%%%%%%%%%%%%%%%%%%%%%%%%%%%%%%%%%%%%%%%%%%%%%%%%%%%%%%%%%%%%%%%%%%%%
\begin{lemma}\label{lem-notproper}
 Suppose that $(G,X)$ is a second-countable, locally-compact transformation
 group for which $G$ (but not necessarily $X$) is Hausdorff.
 Then $(G,X)$ is not Cartan if and only if there exist $z\in X$ and  a sequence
 $(x_n)_n$ in $X$ which converges $2$-times in $X/G$ to $z$.
\end{lemma}
%%%%%%%%%%%%%%%%%%%%%%%%%%%%%%%%%%%%%%%%%%%%%%%%%%%%%%%%%%%%%%%%%%%%%%%%%%%%%%%%%%%%%%%%%%%%%%%%%%%%%%%%%%%%%%%

\begin{proof} Suppose that $(G,X)$ is not Cartan.  Then there exists a point $x$ in $X$
which has no wandering neighbourhood.  Let $(V_n)_n$ be a decreasing sequence of open
neighbourhoods of $z$ in $X$ and let $(K_n)_n$ be an increasing sequence of compact
neighbourhoods of $e$ in $G$ such that $G=\cup_{n=1}^\infty\interior K_n$.

For each $n\geq 1$, since $V_n$ is not a wandering neighbourhood of
$z$, the set $\{s\in G:s\cdot V_n \cap V_n\neq \emptyset\}$ is not
relatively compact, hence not a subset of $K_n$.  So there exists
$x_n\in V_n$ and $s_n\notin K_n$ such that $s_n\cdot x_n\in V_n$.
Both $x_n\to z$ and $s_n\cdot x_n\to z$ since $(V_n)_n$ is a
sequence of decreasing neighbourhoods of $z$.  Moreover,
$s_n\to\infty$ since any compact subset of $G$ is contained in
$\interior K_n$ for some $n$.
%since $s_n\notin K_n$ and $(K_n)_n$ is an increasing
%sequence of compact neighbourhoods of $e$ in $G$.
Thus $(x_n)_n$
converges $2$-times in $X/G$ to $z$.

Conversely, assume there exists a sequence $(x_n)_n$ converging $2$-times in $X/G$ to some $z\in X$.  We may assume (see \cite[Definition~2.2]{AD}) that $x_n\to z$.  So there exists a sequence $(t_n)_n$ in $G$ such that $t_n\cdot x_n\to z$ and $t_n\to\infty$. Let $V$ be any  neighbourhood of $z$ in $X$.  There exists $n_0$ such that $x_n\in V$  and $t_n\cdot x_n\in V$ whenever $n\geq n_0$.  Thus $t_n\in\{s\in G:s\cdot V\cap V\neq \emptyset\}$ for $n\geq n_0$, and since $t_n\to\infty$, $V$ cannot be wandering. Since  $z$ has no wandering neighbourhood, $(G,X)$ is not Cartan.
\end{proof}

%%%%%%%%%%%%%%%%%%%%%%%%%%%%%%%%%%%%%%%%%%%%%%%%%%%%%%%%%%%%%%%%%%%%%%%%%%%%%%%%%%%%%%%%%%%%%%%%%%%%%%%%%%
Suppose that $G$ acts jointly continuously on a $\textrm{T}_1$ locally
compact Hausdorff space $X$.  Then we say that $G$ acts \emph{integrably}
on $X$ if, for every compact subset $N$ of $X$,
$$ \sup_{x\in N} \nu(\{s\in G : s\cdot x\in N\})<\infty$$
(see \cite[Proposition~3.1]{aH} and \cite[Definition~1.10]{rieffel}). Note that the relevant subsets of $G$ are
indeed $\nu$-measurable (in fact, closed) because of the joint
continuity and the $\textrm{T}_1$ property (cf. \cite[p. 400]{AD}). Note also
that if $G$ is non-compact and acts integrably on $X$ then $X$ is
necessarily non-compact too.

\begin{cor}\label{cor-A}
Let $(A, G,\alpha)$ be a separable $C^*$-dynamical system such that the induced action of $G$ on   $\Prim A$    is free. Then $M_U(\Ind\pi)\geq M_U(\pi)$ for all $\pi\in\hat A$.  Moreover,
\begin{enumerate}
\item if $A\rtimes_\alpha G$ is a Fell algebra, then $A$ is
a Fell algebra and $(G,\hat A)$ is Cartan ; and
\item if $A\rtimes_\alpha G$ has bounded trace, then $A$ has bounded trace and the action
of $G$ on $\hat A$ is integrable.
\end{enumerate}
\end{cor}

\begin{proof}
Let $\pi\in \hat A$.  Since $(A, G,\alpha)$ is separable and the induced action   on $\Prim A$    is free,  $\Ind\pi$ is an irreducible representation of $A\rtimes_\alpha G$ by \cite[Proposition~1.7]{EW}, so $M_U(\Ind\pi)$ makes sense.  By \cite[Lemma~1.2]{AK} there exists a sequence $(\pi_n)$ in $\hat A$ such that $\pi_n\to \pi$ and $M_U(\pi)=M_L(\pi,(\pi_n))$.  Apply Theorem~\ref{thm-a} with $k=1$ and $m_1=M_U(\pi)$ to get
\[M_U(\Ind\pi)\geq M_L(\Ind\pi,(\Ind\pi_n))\geq M_U(\pi).
\]

(1) Suppose that $A\rtimes_\alpha G$ is a Fell algebra. For
$\pi\in\hat A$, we have $M_U(\pi)\leq M_U(\Ind\pi)=1$ and so $A$ is
a Fell algebra by \cite[Theorem~4.6]{A}.  If $(G,\hat A)$ is not
Cartan then by Lemma~\ref{lem-notproper} there exist $\pi\in\hat A$
and a sequence $(\pi_n)_n$ in $\hat A$ converging $2$-times in $\hat
A/G$ to $\pi$. By Theorem~\ref{thm-a}, applied with $k=2$ and
$m_1=1=m_2$, we obtain $M_L(\Ind\pi,(\Ind\pi_n))\geq 2$, which
contradicts the fact that $M_U(\Ind\pi)=1$. Thus $(G,\hat A)$ is
Cartan.

(2) Suppose that $A\rtimes_\alpha G$ has bounded trace. Then for
$\pi\in\hat A$, we have $M_U(\pi)\leq M_U(\Ind\pi)<\infty$ and so
$A$ has bounded trace by \cite[Theorem~2.6]{ASS}. As noted above,
$\Ind\pi$ is irreducible for each $\pi\in\hat A$ and so by
\cite[Theorem 3.4 and p.409]{AD} $G$ acts integrably on $\hat A$.
\end{proof}
%%%%%%%%%%%%%%%%%%%%%%%%%%%%%%%%%%%%%%%%%%%%%%%%%%%%%%%%%%%%%%%%%%%%%%%%%%%%%%%%%%%%%%%%%%%%%%%%%%%%%%%%%%%%%%%%%%%%%%%%%%%%%%%%%%

Deicke proves in \cite[Theorem~5.3.2]{D-thesis} that
$M_U(\Ind\pi)\geq M_U(\pi)$ for all $\pi\in\hat A$ and
$M_L(\Ind\pi)\geq M_L(\pi)$ for all $\pi\in\hat A$ such that
$\{\Ind\pi\}$ is not open in $(A\rtimes_\alpha G)^\wedge$, under the
assumption that $(A,G,\alpha)$ is strongly regular (in the sense
that the induced action of $G$ on $\hat A$ is free, the orbits are
locally closed and are canonically homeomorphic to $G$).    Since strong regularity implies $\Prim A\simeq\hat A$, Corollary~\ref{cor-A} captures Deicke's $M_U$ result.   

Deicke's techniques are those of non-abelian duality, and he apparently needs
the strong regularity to ensure that the dual coaction $\hat\alpha$
of $G$ on $A\rtimes_\alpha G$ is pointwise unitary.  When
$(A,G,\alpha)$ is strongly regular all irreducible representations
of $A\rtimes_\alpha G$ are induced, so we can   also    recover Deicke's $M_L$
result by using Theorem~\ref{thm-a}, \cite[Lemma~A.2]{AaH} and the
fact that $\Ind$ is an open map.

\medskip

We now observe  that finite multiplicity numbers in the spectrum of $A$ impose stiff restrictions on the sort of actions on $A$ that can occur. Nevertheless, as we illustrate in Section~\ref{sec-4}, very interesting examples occur.

%%%%%%%%%%%%%%%%%%%%%%%%%%%%%%%%%%%%%%%%%%%%%%%%%%%%%%%%%%%%%%%%%%%%%%%%%%%%%%%%%%%%%%%%%%%%%%%%%%%%%%%%%%%%%%%%%%
\begin{lemma}\label{lem-restriction-on-actions}
Suppose that $(A,G,\alpha)$ is a $C^*$-dynamical system.  Let $\pi\in\hat A$ and suppose that $M_U(\pi)<\infty$. Then either $M_L(\pi)=1$ or the stability subgroup $S_\pi:=\{s\in G:s\cdot\pi=\pi\}$ is open in $G$.
\end{lemma}

\begin{proof}
Suppose that $S_\pi$ is not open.
Then there exists a net $(s_\beta)_\beta$ in $G\setminus S_\pi$ such that $s_\beta\to e$. Thus $s_\beta\cdot\pi\neq\pi$ and $s_\beta\cdot\pi\to\pi$.  In particular, $\{\pi\}$ is not open in $\hat A$.

Now also suppose that  $M_U(\pi)=m<\infty$. Then $M_U(s_\beta\cdot\pi)=m$ for all $\beta$. By \cite[Theorem~1.5]{AKLSS}
\[
m=M_U(\pi)\geq m M_U(\pi,(s_\beta\cdot\pi)).
\]
Thus $1\leq M_L(\pi)\leq M_U(\pi,(s_\beta\cdot\pi))\leq 1$, and hence $M_L(\pi)=1$.
\end{proof}

%%%%%%%%%%%%%%%%%%%%%%%%%%%%%%%%%%%%%%%%%%%%%%%%%%%%%%%%%%%%%%%%%%%%%%%%%%%%%%%%%%%%%%%%%%%%%%%%%%%%%%%%%%%%%%%%%%%

The next Corollary is immediate from Lemma~\ref{lem-restriction-on-actions} since $\{e\}$ is open in $G$ if and only if $G$ is discrete.

%%%%%%%%%%%%%%%%%%%%%%%%%%%%%%%%%%%%%%%%%%%%%%%%%%%%%%%%%%%%%%%%%%%%%%%%%%%%%%%%%%%%%%%%%%%%%%%%%%%%%%%%%%%%%%%%%%%

\begin{cor}
Suppose that $(A,G,\alpha)$ is a  $C^*$-dynamical system and the induced action of $G$ on $\hat A$ is free.  Let $\pi\in \hat A$ and suppose $M_U(\pi)<\infty$.  Then either $M_L(\pi)=1$ or $G$ is discrete.
\end{cor}

%%%%%%%%%%%%%%%%%%%%%%%%%%%%%%%%%%%%%%%%%%%%%%%%%%%%%%%%%%%%%%%%%%%%%%%%%%%%%%%%%%%%%%%%%%%%%%%%%%%%%%%%%%%%%%%%%%%%%%%%%%%%%%%%%%
In Theorem~\ref{3.1'} we find an upper bound on the lower multiplicity of $\Ind\pi$ relative to a sequence $(\Ind\pi_n)$. To do this we follow \cite{aH} and \cite{AD}
in using a vector-valued version of Mercer's Theorem from
\cite[Chapitre V]{duflo}.  Let $\H$ be a separable Hilbert space. We
may take as the fundamental family $\Lambda$ of \cite{duflo} the set
$C_c(G,\H)$ and then $L^2(\Lambda)=L^2(G,\H,\nu)$. A bounded
operator $T\in B(L^2(G,\H))$ is then said to be defined by a kernel
$K\in B(L^2(G\times G,\H))$ if, for all $\xi,\eta\in C_c(G,\H)$,
\begin{enumerate}
\item the function $(s,t)\mapsto \langle K(s,t)(\xi(t))\,,\,\eta(s)\rangle$ is $\nu\times\nu$ integrable; and
\item $\langle T\xi\,,\,\eta\rangle=\int_G\int_G \langle K(s,t)(\xi(t))\,,\,\eta(s)\rangle\, d\nu(s)\, d\nu(t)$.
\end{enumerate}
The kernel $K$ is said to be continuous if, for all $\xi,\eta\in C_c(G,\H)$, the function $(s,t)\mapsto \langle K(s,t)(\xi(t))\,,\,\eta(s)\rangle$ is continuous.
By combining Th\'eor\`eme~3.3.1, Remarque~3.2.1 and Proposition~3.1.1 of \cite[Chapitre V]{duflo} we obtain:
\begin{thm}\label{thm-duflo}\textup{(Duflo)} Let $T$ be a positive operator on $L^2(G,\H)$ defined by a continuous kernel $K\in L^2(G\times G,B(\H))$.  Then $K(s,s)$ is positive for all $s\in G$. If $s\mapsto \tr(K(s,s))$ is integrable
then $T$ is trace-class with
\[
\tr(T)=\int_G\tr(K(s,s))\, d\nu(s).
\]
\end{thm}

Let $\Psi:C_b(\hat A)\to ZM(A)$
be the Dauns--Hofmann isomorphism, so that
\[\sigma(\Psi(\tau)a)=\tau(\sigma)\sigma(a)\] for $\tau\in C_b(\hat
A)$, $\sigma\in \hat A$, and $a\in A$.  If $\lt:G\to\Aut C_0(\hat
A)$ is the action by left translation induced  from the action of
$G$ on $\hat A$, so that
$\lt_s(\tau)(\sigma)=\tau(s^{-1}\cdot\sigma)$, then
$\alpha_s(\Psi(\tau)a)=\Psi(\lt_s(\tau))\alpha_s(a)$.

For $\sigma\in\hat A$, we let $\phi_\sigma:G\to\hat A$ be the
function $s\mapsto s\cdot\sigma$.     Recall that if $\hat A$ is $\textrm{T}_0$  then $\hat A\simeq\Prim A$ via the kernel map.   

%%%%%%%%%%%%%%%%%%%%%%%%%%%%%%%%%%%%%%%%%%%%%%%%%%%%%%%%%%%%%%%%%%%%%%%%%%%%%%%%%%%%%%%%%%%%%%%%%%%%%%%%%%%%%%%%%%%%%%%%%%%%%%%%%%
\begin{thm}\label{3.1'}
Let $(A,G,\alpha)$ be a separable $C^*$-dynamical system such that   $\hat A$ is Hausdorff and the induced action of $G$ on $\hat A\simeq\Prim A$  is free.    Let $\pi\in\hat A$ such that $G\cdot\pi$ is locally closed in $\hat A$ and  $M_U(\pi)=u<\infty$, and  $(\pi_n)_n$ a sequence in $\hat A$. Let $M\in\R$ with $M\geq 1$.  Suppose that for every open neighbourhood $V$ of $\pi$ in $\hat A$ there exists an open neighbourhood $V_1$ of $\pi$ such that   $V_1\subset V$ and
\[
\nu(\{s\in G:s\cdot\pi_n\in V_1\})\leq M\nu(\{s\in G:s\cdot\pi\in V_1\})
\]
frequently.  Then $M_L(\Ind\pi,(\Ind\pi_n))\leq \lfloor M^2u\rfloor$.
\end{thm}
%%%%%%%%%%%%%%%%%%%%%%%%%%%%%%%%%%%%%%%%%%%%%%%%%%%%%%%%%%%%%%%%%%%%%%%%%%%%%%%%%%%%%%%%%%%%%%%%%%%%%%%%%%%%%%%%%%%%%%%%%%%%%%%%%%%%%%%%%%%%

\begin{proof}   By separability and freeness of the action on $\hat A$, if $\sigma\in\hat A$ then $\Ind\sigma\in (A\rtimes_\alpha G)^\wedge$ by \cite[Proposition~1.7]{EW}.
The proof of the Theorem builds on the proof of \cite[Theorem~3.1]{AaH}; in fact we use the construction of a function $F$ from there verbatim.  Fix $\epsilon>0$ such that
\[
\frac{M^2u(1+\epsilon)^2}{1-\epsilon}<\lfloor M^2u\rfloor +1.
\]
 Using the function $F$, we will build an operator  $E\in C_c(G,A)$
and use the generalised lower semi-continuity result
\[
\liminf_n \tr(\Ind\pi_n(E^**E))\geq
M_L(\Ind\pi,(\Ind\pi_n))\tr(\Ind\pi (E^**E))
\]
of \cite[Theorem~4.3]{AS} to bound $M_L(\Ind\pi,(\Ind\pi_n))$.
In the scalar case $A=\C$ of \cite[Theorem~3.1]{AaH} we were able to do this with an operator $D$ in place of $E$ where $\Ind\pi(D^**D)$ is a rank-one projection so that $\tr(\Ind\pi(D^**D))=1$; in the non-scalar case we have to work a little harder.

Since $M_U(\pi)=u<\infty$, by \cite[Theorem~2.5]{ASS} there exist $a\in A^+$ and an open  neighbourhood $V$ of $\pi$ in $\hat A$ such that $\|a\|=1$, $\pi(a)$ is a rank-one projection, and $\sigma(a)$ is a finite-rank operator with rank at most $u$ for all $\sigma\in V$. Since $\hat A$ is Hausdorff, by shrinking $V$ and using the functional calculus, we may assume that $\sigma(a)$ is a projection of rank at most $u$ for all $\sigma\in V$ (once you know how to do this it is standard: see, for example, the second paragraph of the proof of \cite[Theorem~4.6]{A}).
 Since $G\cdot\pi$ is locally closed, we may assume, by \cite[Lemma~2.1]{AaH} and further shrinkage of $V$ if necessary, that $\phi_\pi^{-1}(V)$ is relatively compact,  and then  $\nu(\phi_\pi^{-1}(V))<\infty$ (recall that $\phi_\pi: s\mapsto s\cdot\pi:G\to G\cdot\pi$).

Since $\alpha$ is strongly continuous, there exists an open neighbourhood $N\subset \phi_\pi^{-1}(V)$ of $e$ in $G$ such that
\begin{equation}\label{action-estimate}
\|s\cdot\pi(a)-\pi(a)\|=\|\pi(\alpha_{s^{-1}}(a))-\pi(a)\|<\epsilon\quad\text{for all $s\in N$}.
\end{equation}
Since $G\cdot \pi$ is locally closed in $\hat A$ it follows from
\cite[Theorem~1]{Gli}, applied to the locally compact Hausdorff
transformation group $(G, G\cdot \pi)$, that $\phi_\pi$ is a
homeomorphism of $G$ onto $G\cdot \pi$. So $N\cdot\pi=V_1\cap G\cdot\pi$ for some open neighbourhood $V_1\subset V$ of $\pi$.  By our hypothesis there exists an open neighbourhood $V_2\subset V_1$ of $\pi$ such that
\[
\nu(\phi_{\pi_n}^{-1}(V_2))\leq M\nu(\phi_\pi^{-1}(V_2))<\infty
\]
frequently.

We now give the construction of the function $F$ from \cite[Theorem~3.1]{AaH}.  Let $\delta>0$ such that
\begin{equation*}
\delta<\frac{\epsilon\nu(\phi_\pi^{-1}(V_2))}{1+\epsilon}<\nu(\phi_\pi^{-1}(V_2)).
\end{equation*}
By the regularity of the measure $\nu$ there exists a compact subset $W$ of the
open set $\phi_\pi^{-1}(V_2)$ such that
\begin{equation*}
0<\nu(\phi_\pi^{-1}(V_2))-\delta<\nu(W).
\end{equation*}
Since $W$ is compact, there is a compact neighbourhood $W_1$ of
$W$ contained in $\phi_\pi^{-1}(V_2)$ and a continuous function
$g:G\to[0,1]$ such that $g$ is identically one on $W$ and is
identically zero off the interior of $W_1$. Then
\begin{equation*}
\nu(\phi_\pi^{-1}(V_2))-\delta<\nu(W)\leq \int_G g(t)^2\, dt=\|g\|_2^2,
\end{equation*}
and hence
\begin{equation}\label{eq-estimate}
\frac{\nu(\phi_\pi^{-1}(V_2))}{\|g\|_2^2}<
1+\frac{\delta}{\|g\|_2^2}<1+\frac{\delta}{\nu(\phi_\pi^{-1}(V_2))-\delta}<1+\epsilon.
\end{equation}
There is a continuous
function $g_1:W_1\cdot \pi\to[0,1]$ such that $g_1(t\cdot \pi)=g(t)$ for
$t\in W_1$. Since $W_1\cdot \pi$ is a compact subset of the locally
compact Hausdorff space $\hat A$, it follows from Tietze's Extension
Theorem (applied to the one-point compactification of $\hat A$ if
necessary) that $g_1$  can be extended to a continuous function
$g_2:\hat A\to [0,1]$. Because $W_1\cdot \pi$ is a compact subset of the
open set $V_2$, there exists a compact neighbourhood $P$ of $W_1\cdot
\pi$ contained in $V_2$ and a continuous function $h:\hat A\to [0,1]$ such
that $h$ is identically one on $W_1\cdot \pi$ and is identically zero
off the interior of $P$.
 Note that $h$ has compact support contained in $P$. We set
\[
f(\sigma)=h(\sigma)g_2(\sigma)\hskip1cm(\sigma\in\hat A).
\]
Then $f\in C_c(\hat A)$ with $0\leq f\leq 1$ and  $\supp f\subset\supp
h\subset P\subset V_2$. Note that
\begin{equation}\label{eq-reciprocal}
 \|f_\pi\|^2_2
 %=\int_G f(t\cdot \pi)^2\, dt
 =\int_G h(t\cdot \pi)^2g_2(t\cdot \pi)^2\, dt\geq\int_{W_1} g(t)^2\, dt =\|g\|^2_2
\end{equation}
since $h$ is identically one on $W_1\cdot z$ and $g$ has support
inside $W_1$. We now set
\[
F(\sigma)=\frac{f(\sigma)}{\|f_\pi\|_2}\hskip1cm(\sigma\in\hat A).
\]
Now $F\in C_c(\hat A)$, $\|F_\pi\|_2=1$ and $F_\sigma(s)=F(s\cdot \sigma)\neq 0$
implies $s\in\phi_\sigma^{-1}(V_2)$ by our choice of $h$. Since
$\phi_\pi^{-1}(V_2)$ is relatively compact, $\supp F_\pi$ is compact.

There exists a subsequence $(\pi_{n_i})_i$ of $(\pi_n)_n$ such that
\begin{equation*}
\nu(\phi_{\pi_{n_i}}^{-1}(V_2))\leq M\nu(\phi_\pi^{-1}(V_2))
\end{equation*}
for all $i\geq 1$.
Using \eqref{eq-reciprocal},
\begin{equation}\label{eq-moved}
\int_G F(s\cdot \pi_{n_i})^2\, d\nu(s)\leq\frac{\nu(\phi_{\pi_{n_i}}^{-1}(V_2))}{\|f_\pi\|^2_2}
\leq\frac{M\nu(\phi_\pi^{-1}(V_2))}{\|g\|^2_2}.
\end{equation}

Choose $b\in C_c(G\times \hat A)$ such that $0\leq b\leq 1$ and $b$ is identically
one on the set $(\supp F_\pi)(\supp F_\pi)^{-1}\times\supp F$.
Set
\[
B(t,\sigma)=F(\sigma)F(t^{-1}\cdot \sigma)b(t^{-1},\sigma)\Delta(t)^{-1/2}\hskip1cm(t\in G, \sigma\in\hat A).
\]
Then $B\in C_c(G\times\hat A)$. Let $\Psi:C_b(\hat A)\to ZM(A)$ be the Dauns--Hofmann isomorphism and define $C:G\to A$ by
\[
C(t)=\Psi(B(t,\cdot))a\hskip1cm(t\in G).
\]
Then $C$ has compact support because $B$ does.  To see that $C$ is continuous, fix $\epsilon_2>0$  and let $t_j\to t$ in $G$.  Let $K$ and  $L$ be compact subsets of $G$ and $\hat A$, respectively, such that $\supp B\subset K\times L$.  Since $B$ is continuous, for each $\sigma\in\hat A$ we can choose an open neighbourhood $W_\sigma$ of $t$ in $G$ and an open neighbourhood $U_\sigma$ of $\sigma$ in $\hat A$ such that
$|B(s,\rho)-B(t,\rho)|<\epsilon_2$ whenever $(s,\rho)\in W_\sigma\times U_\sigma$.
Choose a finite subcover $\{U_{\sigma_1}\,\dots, U_{\sigma_l}\}$ of $L$ and set $W_0=\cap_{i=1}^l W_{\sigma_i}$ and $U_0=\cup_{i=1}^l U_{\sigma_i}$.  Since $\|a\|=1$,
\begin{align*}
\|C(t_j)-C(t)\|&=\|\big(\Psi(B(t_j,\cdot))-\Psi(B(t,\cdot))  \big)a\|\\
&\leq \|B(t_j,\cdot)-B(t,\cdot)\|_\infty\|a\|\\
&=\sup_{\rho\in U}\|B(t_j,\rho)-B(t,\rho)\|\leq\epsilon_2
\end{align*}
whenever $t_j\in W_0$. Thus $C\in C_c(G,A)$. Next, set
\[E=\frac{1}{2}(C+C^*)\]
so that $E$ is self-adjoint.  Note that
\begin{align*}
C^*(t)&=\Delta(t^{-1})\alpha_t(C(t^{-1}))^*=\Delta(t)^{-1}\alpha_t(a)\Psi\big( \lt_t(B(t^{-1},\cdot)) \big)\\
&=\Delta(t)^{-1}\alpha_t(a)\Psi\big(\lt_t \big( F\lt_{t^{-1}}(F)b(t,\cdot)\Delta(t)^{1/2}  \big)  \big)\\
&=\Delta(t)^{-1/2}\Psi\big(F\lt_t(F)\lt_t(b(t,\cdot))\big)\alpha_t(a)
\end{align*}
since $a$ is self-adjoint, $B$ is real-valued and $\Psi$ takes values in $ZM(A)$.  Thus
\[
E(t)=\frac{1}{2}\Delta(t)^{-1/2}\Psi(F\lt_t(F))\big( \Psi(b(t^{-1},\cdot))a+\Psi(\lt_t(b(t,\cdot))\alpha_t(a) \big)
\]
and
\begin{align*}
E*E(t)&=\int_G E(r)\alpha_r(E(r^{-1}t))\, d\mu(r)\\
&=\int_G E(r^{-1})\alpha_{r^{-1}}(E(rt))\, d\nu(r)\\
&=\frac{1}{4}\Delta(t)^{-1/2}\int_G\Psi(F\lt_{r^{-1}}(F))\big(\Psi(b(r,\cdot))a+\Psi(\lt_{r^{-1}}(b(r^{-1},\cdot))\alpha_{r^{-1}}(a)   \big)\\\
&\quad\quad \cdot \Psi(\lt_{r^{-1}}(F)\lt_t(F))\big( \Psi(\lt_{r^{-1}}(b(t^{-1}r^{-1},\cdot))\alpha_{r^{-1}}(a)+\Psi(\lt_t(b(rt,\cdot))\alpha_t(a) \big)\, d\nu(r).
\end{align*}

Fix $\sigma\in\hat A$ such that $\nu(\phi_\sigma^{-1}(V_2))<\infty$. For $\xi\in L^2(G,\H_\sigma,\nu)$,
\begin{align*}
(\Ind\sigma&(E*E))\xi)(s)
=
\int_G s\cdot\sigma(E*E(w^{-1}))\xi(ws)\Delta(w)^{-1/2}d\nu(w)\\
&=\frac{1}{4}\int_G\int_G F(s\cdot\sigma)F(rs\cdot\sigma)^2F(ws\cdot\sigma)\big( b(r,s\cdot\sigma) s\cdot\sigma(a)+b(r^{-1},rs\cdot\sigma)rs\cdot\sigma(a)\big)\\
&\hskip2cm
\cdot\big(  b(wr^{-1},rs\cdot\sigma) rs\cdot\sigma(a)+b(rw^{-1},ws\cdot\sigma)ws\cdot\sigma(a) \big)\, d\nu(r)\xi(ws)\, d\nu(w)\\
&=\frac{1}{4}\int_G\int_G F(s\cdot\sigma)F(rs\cdot\sigma)^2F(t\cdot\sigma)\big( b(r,s\cdot\sigma) s\cdot\sigma(a)+b(r^{-1},rs\cdot\sigma)rs\cdot\sigma(a)\big)\\
&\hskip2.1cm
\cdot\big(  b(ts^{-1}r^{-1},rs\cdot\sigma) rs\cdot\sigma(a)+b(rst^{-1},t\cdot\sigma)t\cdot\sigma(a) \big)\, d\nu(r)\xi(t)\, d\nu(t)\\
\intertext{after the change of variable $t=ws$.  Changing variables again, this time setting $v=rs$, this becomes}
&=\frac{1}{4}\int_G\int_G F(s\cdot\sigma)F(v\cdot\sigma)^2F(t\cdot\sigma)\big( b(vs^{-1},s\cdot\sigma) s\cdot\sigma(a)+b(sv^{-1},v\cdot\sigma)v\cdot\sigma(a)\big)\\
&\hskip2.1cm
\cdot\big(  b(tv^{-1},v\cdot\sigma) v\cdot\sigma(a)+b(vt^{-1},t\cdot\sigma)t\cdot\sigma(a) \big)\, d\nu(v)\xi(t)\, d\nu(t)\\
&=\int_GK_\sigma(s,t)\xi(t)\, d\nu(t),
\end{align*}
where
\begin{align*}
K_\sigma(s,t)&=\frac{1}{4}\int_G F(s\cdot\sigma)F(v\cdot\sigma)^2F(t\cdot\sigma)\big( b(vs^{-1},s\cdot\sigma) s\cdot\sigma(a)+b(sv^{-1},v\cdot\sigma)v\cdot\sigma(a)\big)\\
&\hskip4cm
\cdot\big(  b(tv^{-1},v\cdot\sigma) v\cdot\sigma(a)+b(vt^{-1},t\cdot\sigma)t\cdot\sigma(a) \big)\, d\nu(v).
\end{align*}

We claim that  $\Ind\sigma(E*E)$ is the operator in $B(L^2(G,H_\sigma))$ defined by the continuous kernel $K_\sigma\in C_c(G\times G, B(\H_\sigma))$ in the sense of Theorem~\ref{thm-duflo}. We start by verifying the continuity: that for non-zero fixed $\xi,\eta\in C_c(G,\H_\sigma)$, the function  $(s,t)\mapsto\langle K_\sigma(s,t)(\xi(t))\,,\, \eta(s)\rangle$ is continuous.
Define $T\colon G\times G\times G\to B(\H_\sigma)$ by
\begin{align*}
T(s,t,v)=\frac{1}{4} &F(s\cdot\sigma)F(v\cdot\sigma)^2F(t^{-1}\cdot\sigma)\big( b(vs^{-1},s\cdot\sigma) s\cdot\sigma(a)+ b(sv^{-1},v\cdot\sigma)v\cdot\sigma(a)\big)\\
&\cdot\big(  b(tv^{-1},v\cdot\sigma) v\cdot\sigma(a)+b(vt^{-1},t\cdot\sigma)t\cdot\sigma(a) \big).
\end{align*}
Then $T$ is norm continuous and $\|T(s,t,v)\|\leq \|F\|_\infty^4$.
%Note that $K(s_j,t_j)=\int_G T(s_j,t_j,v)\, d\nu(v)$.
Let $s_j\to s$ and $t_j\to t$. We apply the Dominated Convergence Theorem to the continuous functions \[v\mapsto \langle T(s_j,t_j,v)(\xi(t_j))\,,\,\eta(s_j)\rangle\]
with dominating function $v\mapsto \|F\|_\infty^4\|\xi\|_\infty\|\eta\|_\infty\chi_{\phi_\sigma^{-1}(V_2)}(v)$:
\begin{align}
\lim_{j\to\infty}\Big\langle K_\sigma(s_j,t_j)(\xi(t_j))\,,\,\eta(s_j)\Big\rangle
&=\lim_{j\to\infty}\langle \int_G T(s_j,t_j,v)\, d\nu(v)(\xi(t_j))\,,\,\eta(s_j)\rangle\notag\\
&=\lim_{j\to\infty}\int_G\langle T(s_j,t_j,v)(\xi(t_j))\,,\,\eta(s_j)\rangle\, d\nu(v)\notag\\
&=\int_G\lim_{j\to\infty}\langle T(s_j,t_j,v)(\xi(t_j))\,,\,\eta(s_j)\rangle\, d\nu(v)\label{eq-dct}.
\end{align}
Fix $v\in G$ and $\epsilon_3>0$. By continuity of $T,\xi$ and $\eta$, there exist open neighbourhoods $U_1$ and $U_2$ of $s$ and $t$, respectively, such that
$(s_j,t_j)\in U_1 \times U_2$ implies $\|\eta(s_j)-\eta(s)\|\leq \epsilon_3/(3\|F\|_\infty^4\|\xi\|_\infty)$,
$\|\xi(t_j)-\xi(t)\|\leq \epsilon_3/(3\|F\|_\infty^4\|\eta\|_\infty)$, and $\|T(s_j,t_j,v)-T(s,t,v)\|<\epsilon_3/(3\|\xi\|_\infty \|\eta\|_\infty)$.
Now
\begin{align*}
|\langle &T(s_j,t_j,v)(\xi(t_j))\,,\,\eta(s_j)\rangle-\langle T(s,t,v)(\xi(t))\,,\,\eta(s)\rangle|\\
&\leq |\langle T(s_j,t_j,v)(\xi(t_j))\,,\,\eta(s_j)-\eta(s)\rangle|+|\langle T(s_j,t_j,v)(\xi(t_j)-\xi(t))\,,\,\eta(s)\rangle|\\
&\hskip6cm+|\langle (T(s_j,t_j,v)-T(s,t,v)(\xi(t))\,,\,\eta(s)\rangle|<\epsilon_3
\end{align*}
whenever $(s_j,t_j)\in U_1\times U_2$. Now
\[\eqref{eq-dct}=\int_G\lim_{j\to\infty}\langle T(s_j,t_j,v)(\xi(t_j))\,,\,\eta(s_j)\rangle\, d\nu(v)= \langle K_\sigma(s,t)(\xi(t))\,,\,\eta(s)\rangle
\]
as required, and we have shown $K$ is continuous.
Next,
\begin{align*}
\int_{G\times G} |\langle K(s,t)(\xi(t))\,,&\, \eta(s)\rangle| \, d(\nu\times\nu)(s,t)\\
&\leq \|F\|^4_\infty\|\|\xi\|_\infty\|\eta\|_\infty\nu(\supp\xi)\nu(\supp\eta)\nu(\phi_\sigma^{-1}(V_2))<\infty,
\end{align*}
 so $(s,t)\mapsto \langle K(s,t)(\xi(t))\,,\, \eta(s)\rangle$ is $(\nu\times\nu)$-integrable.  Further,
\begin{align*}
\langle \Ind\sigma(E*E)\xi\,,\,\eta\rangle
&=\int_G\langle (\Ind\sigma(E*E)\xi)(s)\,,\,\eta(s)\rangle\, d\nu(s)\\
&=\int_G\Big\langle \int_G K_\sigma(s,t)(\xi(t))\, d\nu(t)\,,\, \eta(s)  \Big\rangle  \, d\nu(s)\\
%&=\int_G\int_G\langle K_\sigma(s,t)(\xi(t))\, d\nu(t)\,,\, \eta(s)\rangle \, d\nu(t)\, d\nu(s)\\
&=\int_G\int_G\langle K_\sigma(s,t)(\xi(t))\, d\nu(t)\,,\, \eta(s)\rangle \, d\nu(s)\, d\nu(t)
\end{align*}
by Fubini's Theorem.  Thus $\Ind\sigma(E*E)$ is the operator defined by the continuous kernel $K_\sigma$ in the sense of Theorem~\ref{thm-duflo}.
By Theorem~\ref{thm-duflo}, or by inspecting the defining integral, $K_\sigma(s,s)\geq 0$ for all $s\in G$.  Before we can apply  Theorem~\ref{thm-duflo} to compute the trace of $\Ind\sigma(E*E)$ we need to verify that $s\mapsto\tr(K(s,s))$ is integrable.

Let $\{h_k\}$ be an orthonormal basis for $\H_\sigma$.  To see that $s\mapsto \tr(K(s,s))$ is measurable,   Fix $s_0\in G$  and let $U$ be a relatively compact neighbourhood of $s_0$ in $G$.  For $k\geq 1$ choose $\xi_k\in C_c(G,\H_\sigma)$ such that $\xi_k(s)=h_k$ for $s\in U$. For $s\in U$,
\[
\langle K(s,s,)h_k\,,\, h_k\rangle=\langle K(s,s,)\xi_k(s)\,,\, \xi_k(s)\rangle\to\langle K(s_0,s_0,)\xi_k(s_0)\,,\, \xi(s_0)\rangle=\langle K(s_0,s_0,)h_k\,,\, h_k\rangle
\]
as $s\to s_0$ since $K$ is continuous as a kernel. Thus $s\mapsto \langle K(s,s,)h_k\,,\, h_k\rangle$ is continuous. Hence $\tr(K(s,s))=\Sigma_{k=1}^\infty \langle K(s,s,)h_k\,,\, h_k\rangle$ is a limit of continuous functions and is measurable.

Next, note that
\begin{align}
\int_G&\tr(K_\sigma(s,s))\, d\nu(s)
=\frac{1}{4}\int_G\tr\Big(\int_G
F(s\cdot\sigma)^2 F(v\cdot\sigma)^2\notag\\
&\hskip2cm \cdot\big( b(vs^{-1},s\cdot\sigma) s\cdot\sigma(a)+b(sv^{-1},v\cdot\sigma)v\cdot\sigma(a)\big)^2\, d\nu(v) \Big)\ d\nu(s).\label{eq-trace1}
\end{align}
We now show that we can move the trace through the inner integral.
 For  fixed $s\in G$, define $A:G\to B(\H_\sigma)^+$ by
\[
A(v)=F(s\cdot\sigma) F(v\cdot\sigma)\big( b(vs^{-1},s\cdot\sigma) s\cdot\sigma(a)+b(sv^{-1},v\cdot\sigma)v\cdot\sigma(a)\big).
\]
Then $A$ is norm continuous, so for each $k\geq 1$, the function $v\mapsto \langle A(v)^2h_k\,,\, h_k\rangle$ is continuous and non-negative on $G$.
Using  the Monotone Convergence Theorem we obtain that
\begin{align*}
\tr(K_\sigma(s,s))
&=\sum_{k=1}^\infty\langle K(s,s)h_k\,,\,h_k\rangle
=\sum_{k=1}^\infty\int_G \langle A(v)^2h_k\,,\, h_k\rangle\, d\nu(v)\\
&=\int_G \sum_{k=1}^\infty\langle A(v)^2h_k\,,\, h_k\rangle\, d\nu(v)=\int_G\tr(A(v)^2\, d\nu(v).
\end{align*}
It now follows from \eqref{eq-trace1} that
\begin{align}
\int_G&\tr(K_\sigma(s,s))\, d\nu(s)=\frac{1}{4}\int_G\int_G
F(s\cdot\sigma)^2 F(v\cdot\sigma)^2\notag\\
&\hskip1.5cm \cdot\tr\big(\big( b(vs^{-1},s\cdot\sigma) s\cdot\sigma(a)+b(sv^{-1},v\cdot\sigma)v\cdot\sigma(a)\big)^2\big)\, d\nu(v) \, d\nu(s).\label{trace-formula}
\end{align}
If $s,v\in\phi_\sigma^{-1}(V_2)$, then $s\cdot\sigma(a)$ and $v\cdot\sigma(a)$ are projections of rank at most $u$. So
by \cite[Lemma~3.4.10]{ped},
\begin{align}
\tr\big(\big( b(vs^{-1},s\cdot\sigma) s\cdot\sigma(a)&+b(sv^{-1},v\cdot\sigma)v\cdot\sigma(a)\big)^2\big)\notag\\
&\leq \|b(vs^{-1},s\cdot\sigma) s\cdot\sigma(a)+b(sv^{-1},v\cdot\sigma)v\cdot\sigma(a)  \|\notag\\
&\hskip2cm\cdot\tr\big(b(vs^{-1},s\cdot\sigma) s\cdot\sigma(a)+b(sv^{-1},v\cdot\sigma)v\cdot\sigma(a)\big)\notag\\
&\leq 2(u+u)=4u\label{trace-estimate1}
\end{align}
whenever $s,v\in\phi_\sigma^{-1}(V_2)$, since $0\leq b\leq 1$.
Combining \eqref{trace-formula} and  \eqref{trace-estimate1} we obtain
\[
\int_G\tr(K_\sigma(s,s))\, d\nu(s)\leq u\left(\int_G
F(s\cdot\sigma)^2\, d\nu(s)\right)^2\leq u\|F\|_\infty^4\nu(\phi_\sigma^{-1}(V_2))^2<\infty.
\]
So by Theorem~\ref{thm-duflo},
\[
\tr(\Ind\sigma(E*E)) =\int_G\tr(K_\sigma(s,s))\, d\nu(s).
\]
This holds in particular when we replace $\sigma$ by $\pi$ or  $\pi_{n_i}$,
and we first use it to obtain an upper bound for the trace of $\Ind\pi_{n_i}(E*E)$:
\begin{align}
\tr(\Ind\pi_{n_i}(E*E)) &=\int_G\tr(K_{\pi_{n_i}}(s,s))\, d\nu(s)\label{eq-for3.5'}\\
&\leq u\left(\int_G
F(s\cdot\pi_{n_i})^2\, d\nu(s)\right)^2\notag\\
&\leq  \frac{uM^2(\phi_\pi^{-1}(V_2))^2}{{\|g\|^4_2}}\quad\text{(using \eqref{eq-moved})}\notag\\
&<uM^2(1+\epsilon)^2\quad\text{(using \eqref{eq-estimate})}.\label{bound1}
\end{align}

Next, we look for a lower bound of the trace of $\Ind\pi(E*E)$.  Let $R_\pi=\{s\in G\colon F(s\cdot\pi\neq0\}$. If $s,v\in R_\pi$ then $s,v\in\phi_\pi^{-1}(V_2)$ so that $\|v\cdot\pi(a)-s\cdot\pi(a)\|<2\epsilon$  by \eqref{action-estimate}; thus
$s\cdot\pi(a)-2\epsilon 1\leq v\cdot\pi(a)$ and hence
\begin{align*}
\tr(s\cdot\pi(a) v\cdot\pi(a))
&=\tr( (s\cdot\pi(a)^{1/2} v\cdot\pi(a)  (s\cdot\pi(a)^{1/2})\notag\\
&\geq\tr( (s\cdot\pi(a)^{1/2} (s\cdot\pi(a) -2\epsilon 1) (s\cdot\pi(a)^{1/2}   )\\
&=\tr( (s\cdot\pi(a))^2 -2\epsilon s\cdot\pi(a)   )\\
&\geq 1-2\epsilon
\end{align*}
because $s\cdot\pi(a)$ is a non-zero projection.  Since $R_\pi\subset \supp F_\pi$ and $b$ is identically
one on the set $(\supp F_\pi)(\supp F_\pi)^{-1}\times\supp F$,
\begin{align*}
\frac{1}{4}\tr\big( &\big( b(vs^{-1},s\cdot\pi) s\cdot\pi(a)+b(sv^{-1},v\cdot\pi)v\cdot\pi(a)\big)^2 \big)\\
&=\frac{1}{4}\tr\big(  s\cdot\pi(a)+s\cdot\pi(a)v\cdot\pi(a)+v\cdot\pi(a)s\cdot\pi(a)+v\cdot\pi(a) \big)\geq 1-\epsilon
\end{align*}
whenever $s,v\in R_\pi$. Together with \eqref{trace-formula} we obtain
\begin{align}
\tr(\Ind\pi(E*E))&\geq (1-\epsilon)\left(\int_{R_\pi}
F(s\cdot\pi)^2\, d\nu(s)\right)^2
\notag\\
&=(1-\epsilon)\|F_\pi\|_2^4=1-\epsilon.\label{bound2}
\end{align}
Finally, using \eqref{bound1} and \eqref{bound2} we have
\begin{align*}
uM^2(1+\epsilon)^2&
\geq\liminf_n \tr(\Ind\pi_n(E*E))\\
&\geq
M_L(\Ind\pi,(\Ind\pi_n))\tr(\Ind\pi(E^**E)))
\\
&\geq M_L(\Ind\pi,(\Ind\pi_n))(1-\epsilon),
\end{align*}
and hence
\[
M_L(\Ind\pi,(\Ind\pi_n))\leq\frac{uM^2(1+\epsilon)^2}{1-\epsilon}<\lfloor M^2u\rfloor +1
\]
by our choice of $\epsilon$.  Thus $M_L(\Ind\pi,(\Ind\pi_n))\leq \lfloor M^2u\rfloor$.
\end{proof}

%%%%%%%%%%%%%%%%%%%%%%%%%%%%%%%%%%%%%%%%%%%%%%%%%%%%%%%%%%%%%%%%%%%%%%%%%%%%%%%%%%%%%%%%%%%%%%%%%%%%%%%%%%%%%%%%%%%%%%%%%%%%%%

Theorem~\ref{3.1'} and \cite[Theorem~3.5]{AaH} now provide all the ideas needed to sharpen the upper bound of Theorem~\ref{3.1'} from $\lfloor M^2u\rfloor$ to $\lfloor Mu\rfloor$. Theorem~\ref{3.5'} has the same hypotheses as
Theorem~\ref{3.1'} but a stronger conclusion; the
proof requires Theorem~\ref{3.1'} and  builds heavily on the proof of \cite[Theorem~3.5]{AaH}.
We refer to \cite[Theorem~3.5]{AaH} for quite a few of the details.
%%%%%%%%%%%%%%%%%%%%%%%%%%%%%%%%%%%%%%%%%%%%%%%%%%%%%%%%%%%%%%%%%%%%%%%%%%%%%%%%%%%%%%%%%%%%%%%%%%%%%%%%%%%%%%%%%%%%%%%%%%%%%%%%%%

\begin{thm}\label{3.5'}
Let $(A,G,\alpha)$ be a separable $C^*$-dynamical system such that   $\hat A$ is Hausdorff and the induced action of $G$ on $\hat A\simeq\Prim A$  is free.    Let $\pi\in\hat A$ such that $G\cdot\pi$ is locally closed in  $\hat A$ and  $M_U(\pi)=u<\infty$, and  $(\pi_n)_n$ a sequence in $\hat A$. Let $M\in\R$ with $M\geq 1$.  Suppose that for every open neighbourhood $V$ of $\pi$ in $\hat A$ there exists an open neighbourhood $V_1$ of $\pi$ such that $V_1\subset V$ and
\[
\nu(\{s\in G:s\cdot\pi_n\in V_1\})\leq M\nu(\{s\in G:s\cdot\pi\in V_1\})
\]
frequently.  Then $M_L(\Ind\pi,(\Ind\pi_n))\leq \lfloor Mu\rfloor$.
\end{thm}

\begin{proof}
If $\Ind\pi_n\not\to\Ind\pi$ then $M_L(\Ind\pi,(\Ind\pi_n))=0<\lfloor Mu\rfloor$. So we assume from now on that $\Ind\pi_n\to\Ind\pi$, or equivalently, since $\Ind$ induces a homeomorphism of $\hat A/G$ onto $(A\rtimes_\alpha G)^\wedge$, that $G\cdot \pi_n\to G\cdot \pi$.  By Theorem~\ref{3.1'}, $M_L(\Ind\pi,(\Ind\pi_n))\leq \lfloor Mu\rfloor<\infty$, so by \cite[Proposition~3.4]{AaH} there exists an open neighbourhood $U$ of $\Ind\pi$ such that $\Ind\pi$ is the unique limit of $(\Ind\pi_n)_n$ in $U$.  Let $q:\hat A\to \hat A/G$ be the quotient map and let $J$ be the closed two-sided $G$-invariant ideal of $A$ such that $\hat J\simeq q^{-1}(U)$.  Then  $(J\rtimes_\alpha G)^\wedge \simeq U$. Note that $\pi_n\in q^{-1}(U)$ eventually.  The multiplicity $M_L(\Ind\pi,(\Ind\pi_n))$ is the same whether we compute it in $J\rtimes_\alpha G$ or in $A\rtimes_\alpha G$ (see \cite[Proposition~5.3]{ASS}).  If $V$ is an open neighbourhood of $z$ in $q^{-1}(U)$ then $V$ is an open neighbourhood of $z$ in $X$ and there exists a neighbourhood $V_1\subset V$ such that $\nu(\phi_{\pi_n}^{-1}(V_1))\leq M\nu(\phi_\pi^{-1}(V_1))$;  note that $V_1$ is also a neighbourhood of $q^{-1}(U)$.
Thus we may replace  $A\rtimes_\alpha G$ by $J\rtimes_\alpha G$ where $\Ind\pi$ is the unique limit of $(\Ind\pi_n)_n$ in $(J\rtimes_\alpha G)^\wedge$.

We now revert back to $A\rtimes_\alpha G$ and assume that $G\cdot \pi$ is the unique limit of $(G\cdot\pi_n)_n$ in $\hat A/G$.
Fix $\epsilon>0$ such that \[\frac{Mu(1+\epsilon)^2}{1-\epsilon}<\lfloor Mu\rfloor+1.\]
Let $a$ and $V_2$ be as at the beginning of the proof of \ref{3.1'}.
Thus  $a\in A^+$ and $V_2$ is an open neighbourhood of $\pi$ in $\hat A$ such that $\|a\|=1$, $\pi(a)$ is a rank-one projection, and $\sigma(a)$ is a finite-rank projection with rank at most $u$ for all $\sigma\in V_2$.  Further, $\phi_\pi^{-1}(V_2)$ is relatively compact and $\|s\cdot\pi(a)-\pi(a)\|=\|\pi(\alpha_{s^{-1}}(a))-\pi(a)\|<\epsilon$ for all $s\in\phi_\pi^{-1}(V_2)$.
Moreover,
$
\nu(\phi_{\pi_n}^{-1}(V_2))\leq M\nu(\phi_\pi^{-1}(V_2))<\infty
$
frequently, so there exists a subsequence $\{\pi_{n_i}\}$ such that
\[
\nu(\phi_{\pi_{n_i}}^{-1}(V_2))\leq M\nu(\phi_\pi^{-1}(V_2))<\infty
\]
for all $i\geq 1$.

Now we use estimates from \cite[Theorem~3.5]{AaH}. Choose $\gamma>0$ such that
\begin{equation*}\label{eq-Mgamma}
\gamma<\frac{\epsilon\nu(\phi_\pi^{-1}(V_2))}{1+\epsilon}.
%<\nu(\phi_\pi^{-1}(V_2)).
\end{equation*}
By \cite[Lemma~3.3]{AaH} there exists an open  relatively compact
neighbourhood $V_3$ of $\pi$ such that $\overline{V_3}\subset V_2$ and
\begin{equation*}\label{eq-MV_1}
0<\nu(\phi_\pi^{-1}(V_2))-\gamma
<\nu(\phi_\pi^{-1}(V_3))
\leq\nu(\phi_\pi^{-1}(\overline{V_3}))
\leq\nu(\phi_\pi^{-1}(V_2)).
%<\nu(\phi_\pi^{-1}(V_1))+\gamma.
\end{equation*}
Then
\begin{equation}\label{eq-smallV}
\nu(\phi_{\pi_{n_i}}^{-1}(V_3))
<M(1+\epsilon)\nu(\phi_\pi^{-1}(V_3))
\end{equation}
for all $i\geq 1$ (see the calculation culminating at \cite[Equation~3.9]{AaH}), and there exists
$\delta>0$ such that $\delta<\nu(\phi_\pi^{-1}(V_3))$ and
\begin{equation}\label{eq-delta}
\frac{\nu(\phi_\pi^{-1}(V_3))\big(\nu(\phi_\pi^{-1}(\overline{V_3}))+\delta)
\big)}{\big( \nu(\phi_\pi^{-1}(V_3))-\delta \big)^2}
%<\frac{\nu(\phi_\pi^{-1}(V_3))\big(\nu(\phi_\pi^{-1}(V_3))+\gamma+\delta)
%\big)}{\big( \nu(\phi_\pi^{-1}(V_3))-\delta \big)^2}
<1+\epsilon.
\end{equation}

Let $g\in C_c(G)$ and  $f, F\in C_c(\hat A)$ be the functions constructed in Theorem~\ref{3.1'} (but replacing $V_2$ with $V_3$).
Thus
\begin{gather}
\nu(\phi_\pi^{-1}(V_3))-\delta< \int_G g(t)^2\, dt=\|g\|_2^2,\label{eq-g}\\
\|f_\pi\|^2_2=\int_G f(t\cdot \pi)^2\, dt
 \geq\|g\|^2_2\label{eq-reciprocal2}
\end{gather}
and $\|F_\pi\|_2=1$ and $F_\sigma(s)=F(s\cdot \sigma)\neq 0$ implies $s\in\phi_\sigma^{-1}(V_3)$.

Let $K$ be an open relatively compact symmetric neighbourhood of
$(\supp F_\pi)(\supp F_\pi)^{-1}$ in $G$ and $L$ an open relatively
compact neighbourhood of $\supp F$ in $\hat A$. Choose $b\in C_c(G\times
\hat A)$ such that $0\leq b\leq 1$ and $b$ is identically one on the set
$(\supp F_\pi)(\supp F_\pi)^{-1}\times\supp F$ and $b$ is identically
zero off $K\times L$. (Thus $b$ is as in Theorem~\ref{3.1'},
but we have rounded it off.) For $t\in G$ and $\sigma\in\hat A$ set
\begin{equation*}
B(t,\sigma)=F(\sigma)F(t^{-1}\cdot \sigma)b(t^{-1},\sigma)\Delta(t)^{-1/2},\ \
C(t)=\Psi(B(t,\cdot))a\text{\ \ and\ }
E=\frac{1}{2}(C+C^*).
\end{equation*}
>From \eqref{eq-for3.5'} and \eqref{trace-formula},
\begin{align*}
\tr(\pi_{n_i}(E*E)
&=\frac{1}{4}\int_G
F(s\cdot\pi_{n_i})^2 \Bigg(\int_G F(v\cdot\pi_{n_i})^2\notag\\
&\hskip1.1cm \cdot\tr\big(\big( b(vs^{-1},s\cdot\pi_{n_i}) s\cdot\pi_{n_i}(a)+b(sv^{-1},v\cdot\pi_{n_i})v\cdot\pi_{n_i}(a)\big)^2\big)\, d\nu(v)\Bigg) \, d\nu(s).
\end{align*}
The inner integrand is zero unless $v\in \phi_{\pi_{n_i}}^{-1}(V_3)\cap Ks$ because $F_{\pi_{n_i}}(v)=F(v\cdot \pi_{n_i})\neq 0$ implies $v\in\phi_{\pi_{n_i}}^{-1}(V_3)$ and $b$ is identically zero off $K\times L$. Thus (see the trace estimate at \eqref{trace-estimate1} in Theorem~\ref{3.1'}),
\begin{align*}
\tr(\Ind\epsilon_{\pi_{n_i}}(E**E))
&\leq u\int_{s\in\phi_{\pi_{n_i}}^{-1}(V_3)}F(s\cdot \pi_{n_i})^2\Bigg( \int_{v\in\phi_{\pi_{n_i}}^{-1}(V_3)\cap Ks} F(v\cdot \pi_{n_i})^2\, d\nu(v)\Bigg)\, d\nu(s)\\
&\leq \frac{u}{\|f_\pi\|_2^4}\int_{s\in\phi_{\pi_{n_i}}^{-1}(V_3)}1\Bigg( \int_{v\in\phi_{\pi_{n_i}}^{-1}(V_3)\cap Ks}1\, d\nu(v)\Bigg)\, d\nu(s).
\end{align*}

Choose an open neighbourhood $U$ of $\phi_\pi^{-1}(\overline{V_3})$ such that
$\nu(U)<\nu(\phi_\pi^{-1}(\overline{V_3}))+\delta$. By
\cite[Lemma~3.2]{AaH}, applied with $\overline{V_3}$, $\overline{K}$ and $U$,
there exists $i_0$ such that, for every $i\geq
i_0$ and every $s\in \phi_{\pi_{n_i}}^{-1}(\overline{V_3})$ there exists
$r\in\phi_\pi^{-1}(\overline{V_3})$ with
$\overline{K}s\cap\phi_{\pi_{n_i}}^{-1}(\overline{V_3})\subset Ur^{-1}s$.
It follows that
\[
\nu\big( \overline{K}s\cap\phi_{\pi_{n_i}}^{-1}( \overline{V_3})
\big)\leq\nu(U)<\nu(\phi_\pi^{-1}(\overline{V_3}))+\delta
\]
by the right-invariance of $\nu$.
So, provided $i\geq i_0$,
\begin{align*}
\tr(\Ind\epsilon_{\pi_{n_i}}(E*E))&\leq \frac{u\nu(\phi_{\pi_{n_i}}^{-1}(V_3))\big(\nu(\phi_\pi^{-1}(\overline{V_3}))+\delta \big)}{\|f_\pi\|_2^4}\\
&<\frac{uM(1+\epsilon)\nu(\phi_\pi^{-1}(V_3))\big(\nu(\phi_\pi^{-1}(\overline{V_3}))+\delta
\big)}{\|g\|_2^4}
\text{\quad\quad using \eqref{eq-smallV} and \eqref{eq-reciprocal2}}\\
&\leq\frac{uM(1+\epsilon)\nu(\phi_\pi^{-1}(V_3))\big(\nu(\phi_\pi^{-1}(\overline{V_3}))+\delta
\big)}{(\nu(\phi_\pi^{-1}(V_3))-\delta)^2}
\text{\quad\quad using \eqref{eq-g}}\\
&<uM(1+\epsilon)^2\text{\qquad\qquad\qquad\qquad\qquad\qquad\qquad using \eqref{eq-delta}}.
\end{align*}
On the other hand, as in Theorem~\ref{3.1'},
\begin{equation*}
\tr(\Ind\pi(E*E))
\geq (1-\epsilon)\|F_\pi\|_2^4=1-\epsilon.
\end{equation*}
Thus
\begin{align*}
uM(1+\epsilon)^2&
\geq\liminf_n \tr(\Ind\pi_n(E*E))\\
&\geq
M_L(\Ind\pi,(\Ind\pi_n))\tr(\Ind\pi(E^**E)))
\\
&\geq M_L(\Ind\pi,(\Ind\pi_n))(1-\epsilon),
\end{align*}
and hence
\[
M_L(\Ind\pi,(\Ind\pi_n))\leq\frac{uM(1+\epsilon)^2}{1-\epsilon}<\lfloor Mu\rfloor +1
\]
by our choice of $\epsilon$.  Thus $M_L(\Ind\pi,(\Ind\pi_n))\leq \lfloor Mu\rfloor$.
\end{proof}
%%%%%%%%%%%%%%%%%%%%%%%%%%%%%%%%%%%%%%%%%%%%%%%%%%%%%%%%%%%%%%%%%%%%%%%%%%%%%%%%%%%%%%%%%%%%%%%%%%%%%%%%%%%%%%%

%%%%%%%%%%%%%%%%%%%%%%%%%%%%%%%%%%%%%%%%%%%%%%%%%%%%%%%%%%%%%%%%%%%%%%%%%%%%%%%%%%%%%%%%%%%%%%%%%%%%%%%%%%%%%%%%
\begin{cor}\label{cor-Mupper}
Let $(A,G,\alpha)$ be a separable $C^*$-dynamical system such that   $\hat A$ is Hausdorff and the induced action of $G$ on $\hat A\simeq\Prim A$  is free.    Let $\pi\in\hat A$ such that $G\cdot\pi$ is locally closed in $\hat A$ and   $M_U(\pi)=u<\infty$, and $(\pi_n)_n$ a sequence in $\hat A$. Let $M\in\R$ with $M\geq 1$.  Suppose that for every open neighbourhood $V$ of $\pi$ in $\hat A$ there exists an open neighbourhood $V_1$ of $\pi$ such that $V_1\subset V$ and
\[
\nu(\{s\in G:s\cdot\pi_n\in V_1\})\leq M\nu(\{s\in G:s\cdot\pi\in V_1\})
\]
eventually.  Then $M_U(\Ind\pi,(\Ind\pi_n))\leq \lfloor Mu\rfloor$.
\end{cor}

\begin{proof}
Since $A\rtimes_\alpha G$ is separable, by \cite[Lemma~A1]{AaH}   there exists a
subsequence $(\Ind\pi_{n_i})_i$ such that
\[M_U(\Ind\pi,(\Ind\pi_n))=M_U(\Ind\pi,(\Ind\pi_{n_i}))=M_L(\Ind\pi,(\Ind\pi_{n_i})).\]
 By Theorem~\ref{3.5'}, $M_L(\Ind\pi,(\Ind\pi_{n_i}))\leq \lfloor
uM\rfloor$, and hence $M_U(\Ind\pi,(\Ind\pi_n))\leq \lfloor Mu\rfloor$.
\end{proof}
%%%%%%%%%%%%%%%%%%%%%%%%%%%%%%%%%%%%%%%%%%%%%%%%%%%%%%%%%%%%%%%%%%%%%%%%%%%%%%%%%%%%%%%%%%%%%%%%%%%%%%%%%%%%%%%%%%%
\begin{cor}\label{cor-Mabsolute}
Let $(A,G,\alpha)$ be a separable $C^*$-dynamical system such that   $\hat A$ is Hausdorff, the induced action of $G$ on $\hat A\simeq\Prim A$  is free  and all orbits are locally closed in $\hat A$.      Let $\pi\in\hat A$ with $M_U(\pi)=u<\infty$.  Suppose that for every sequence $(\pi_n)_n$ in $\hat A$ converging to $\pi$ and every open neighbourhood $V$ of $\pi$ in $\hat A$ there exists an open neighbourhood $V_1$ of $\pi$ such that $V_1\subset V$ and
\[
\nu(\{s\in G:s\cdot\pi_n\in V_1\})\leq M\nu(\{s\in G:s\cdot\pi\in V_1\})
\]
frequently.  Then $M_U(\Ind\pi)\leq \lfloor Mu\rfloor$.
\end{cor}

\begin{proof}
Since  $A\rtimes_\alpha G$ is separable  it follows from
\cite[Lemma~1.2]{AK} that there exists a sequence $(\psi_n)_{n\geq
1}$ in $(A\rtimes_\alpha G)^\wedge$ converging to $\Ind\pi$, such that
\begin{equation*}
M_L(\Ind\pi,(\psi_n))=M_U(\Ind\pi,(\psi_n))=M_U(\Ind\pi).
\end{equation*}
For each $\sigma\in\hat A$, the orbit $G\cdot\sigma$ is locally closed in $\hat A$, so it follows from
\cite[Theorem~1]{Gli}, applied to the locally compact Hausdorff
transformation group $(G, G\cdot \sigma)$, that $s\mapsto s\cdot\sigma$ is a
homeomorphism of $G$ onto $G\cdot \sigma$.
Thus $\sigma\mapsto
\Ind\sigma$ induces a homeomorphism of $\hat A/G$ onto $(A\rtimes_\alpha G)^\wedge$ by Theorem~\ref{thm-green}.
So there exists
a sequence $(\pi_i)_{i\geq 1}$ in $\hat A$ converging to $\pi$ such that
$(\Ind\pi_i)_{i\geq 1}$ is a subsequence of $(\psi_n)_{n\geq
1}$. By Theorem~\ref{3.5'},
$M_L(\Ind\pi,(\Ind\pi_i))\leq \lfloor Mu\rfloor$.
Since
\begin{align*}
M_U(\Ind\pi)&=M_L(\Ind\pi,(\psi_n))\leq M_L(\Ind\pi,(\Ind\pi_i))\leq M_U(\Ind\pi,(\Ind\pi_i))\\
&\leq M_U(\Ind\pi,(\psi_n))=M_U(\Ind\pi),
\end{align*}
we obtain $M_U(\Ind\pi)\leq\lfloor Mu\rfloor$.
\end{proof}
%%%%%%%%%%%%%%%%%%%%%%%%%%%%%%%%%%%%%%%%%%%%%%%%%%%%%%%%%%%%%%%%%%%%%%%%%%%%%%%%%%%%%%%%%%%%%%%%%%%%%%%%%%%%%%%%%%%

%%%%%%%%%%%%%%%%%%%%%%%%%%%%%%%%%%%%%%%%%%%%%%%%%%%%%%%%%%%%%%%%%%%%%%%%%%%%%%%%%%%%%%%%%%%%%%%%%%%%%%%%%%%%%%%%%%%
\begin{thm}\label{thm-circle2} Let $(A,G,\alpha)$ be a separable $C^*$-dynamical system
such that   $\hat A$ is Hausdorff and the induced action of $G$ on $\hat A\simeq\Prim A$  is free.   
Let $k$ be a positive integer, let $\pi\in \hat A$ such that $G\cdot \pi$ is locally
closed in $\hat A$  and $M_U(\pi)=m<\infty$. Let  $(\pi_n)$ be a sequence in $\hat A$.
Consider the following eight conditions:
\begin{enumerate}
\item the sequence $(\pi_n)_n$ converges $k$-times in $\hat A/G$ to $\pi$;
\item in $(C_0(\hat A)\rtimes_{\lt} G)^\wedge$, $M_L(\Ind\epsilon_\pi,(\Ind\epsilon_{\pi_n}))\geq k$;
\item for every open  neighbourhood $V$ of $\pi$ in $\hat A$  such that
$\{s\in G:s\cdot \pi\in V\}$  is relatively compact we have
\[
\liminf_n\nu(\{s\in G:s\cdot \pi_n\in V\})\geq k\nu(\{s\in G:s\cdot \pi\in V\});
\]
\item there exists a real number $R>k-1$ such that for every open
neighbourhood $V$ of $\pi$ in $\hat A$ with $\{s\in G:s\cdot \pi\in V\}$ relatively
compact we have
\[
\liminf_n\nu(\{s\in G:s\cdot \pi_n\in V\})\geq R\nu(\{s\in G:s\cdot
\pi\in V\});
\]
\item there exists a decreasing sequence of basic compact
neighbourhoods $(W_m)_{m\geq 1}$ of $\pi$ in $\hat A$ such that, for each
$m\geq1$,
\[
\liminf_n\nu(\{s\in G:s\cdot \pi_n\in W_m\})> (k-1)\nu(\{s\in G:s\cdot \pi\in
W_m\});
\]
\item in $(A\rtimes_\alpha G)^\wedge$, $M_L(\Ind\pi,(\Ind\pi_n)_n)\geq km$;
\item given $\epsilon>0$, there exists an open neighbourhood $U_\epsilon$ of $\pi$ in $\hat A$ such that, for all open neighbourhoods $U$ of $\pi$ with $U\subset U_\epsilon$, we have
\[
\nu(\{s\in G:s\cdot \pi_n\in U\})>(k-\epsilon)\nu(\{s\in G:s\cdot
\pi\in U\})
\]
eventually;
\item there exist a real number $R>k-1$ and an open neighbourhood $V$ of $\pi$ in $\hat A$ such that, for all open neighbourhoods $U$ of $\pi$ with $U\subset V$, we have
\[
\liminf_n\nu(\{s\in G:s\cdot \pi_n\in U\})\geq R\nu(\{s\in G:s\cdot
\pi\in U\}).
\]
\end{enumerate}
Then (1)--(5) are equivalent, and (6) $\Longrightarrow$ (7) $\Longrightarrow$ (8) $\Longrightarrow $ (5).  If $m=M_U(\pi)=1$ then  (1) $\Longrightarrow $(6), and hence (1)--(8) are equivalent.
\end{thm}

%%%%%%%%%%%%%%%%%%%%%%%%%%%%%%%%%%%%%%%%%%%%%%%%%%%%%%%%%%%%%%%%%%%%%%%%%%%%%%%%%%%%%%%%%%%%%%%%%%%%%%%%%%%%%%%%%%%%%%%%%
\begin{proof}
The equivalence of (1)--(5) is \cite[Theorem~1.1]{AaH}.  We start by showing that (6) $\Longrightarrow $ (7) $\Longrightarrow$ (8) $\Longrightarrow $ (5).

(6) $\Longrightarrow$ (7).
Assume that (7) fails.  So
for some $\epsilon>0$ there is no such $U_\epsilon$.  Thus, for every open neighbourhood $V$ of $\pi$ in $\hat A$ there exists an open neighbourhood $V_1$ of $\pi$ with $V_1\subset V$ such that
\[
\nu(\phi_{\pi_n}^{-1}(V_1))\leq(k-\epsilon)\nu(\phi_{\pi}^{-1}(V_1))
\]
frequently. Thus by Theorem~\ref{3.5'}, $M_L(\Ind\pi,(\Ind\pi_n))\leq\lfloor (k-\epsilon)m\rfloor<km$, that is (6) fails.

(7) $\Longrightarrow$ (8).
Assume (7).  Let $\epsilon=1/2$ and set $V=U_{1/2}$.  Then for all open neighbourhoods $U$ of $\pi$ with $U\subset V$ we have
\[
\nu(\phi_{\pi_n}^{-1}(U))>(k-1/2)\nu(\phi_{\pi}^{-1}(U))
\]
eventually.  Thus
\[
\liminf_n\nu(\phi_{\pi_n}^{-1}(U))\geq (k-1/2)\nu(\phi_{\pi}^{-1}(U))
\]
and so we may take $R=k-1/2$.

(8) $\Longrightarrow$ (5).  Assume (8).  Let $(V_j)_j$ be a decreasing sequence of basic open neighbourhoods of $\pi\in\hat A$ such that $V_1\subset V$ and $\phi_\pi^{-1}(V_1)$ is relatively compact (such a $V_1$  exists by \cite[Lemma~2.1]{AaH} because $G\cdot\pi$ is locally closed). Arguing as in the \emph{proof} of \cite[Lemma~5.1]{AaH} there exists a compact neighbourhood $W_1$ of $\pi$ such that $W_1\subset V_1$ and
\[
\liminf_n\nu(\phi_{\pi_n}^{-1}(W_1))>(k-1)\nu(\phi_\pi^{-1}(W_1)).
\]
Now assume there are compact neighbourhoods $W_1,\dots,W_m$ of $z$ with $W_1\supset\dots\supset W_m$ such that $W_i\subset V_i$ and $\nu(\phi_{\pi_n}^{-1}(W_i))>(k-1)\nu(\phi_\pi^{-1}(W_i))$ for $1\leq i\leq m$.  Note that $U:=(\interior W_m)\cap V_{m+1}\subset V_1\subset V$, so $\phi_\pi^{-1}(U)$ is relatively compact  and \[\liminf_n\nu(\phi_{\pi_n}^{-1}(U))\geq R\nu(\phi_\pi^{-1}(U)).\]
Arguing as in the proof of \cite[Lemma~5.1]{AaH} again, there exists a compact neighbourhood $W_{m+1}\subset U$ such that
\[\liminf_n\nu(\phi_{\pi_n}^{-1}(W_{m+1}))>(k-1)\nu(\phi_\pi^{-1}(W_{m+1})).\]

Finally, assume that $m=M_U(\pi)=1$ and that (1) holds.
We apply Theorem~\ref{thm-a} to $(\pi_n)_n$ and $\pi$ with $m_1=m_2=\cdots =m_k=1$ to conclude that $M_L(\Ind\pi,(\Ind\pi_n))\geq k$.
Thus (1) $\Longrightarrow$ (6)  if $m=M_U(\pi)=1$.  Hence (1)--(8) are equivalent if $m=M_U(\pi)=1$.
\end{proof}

\begin{cor}
Let $(A,G,\alpha)$ be a separable $C^*$-dynamical system such that   $\hat A$ is Hausdorff and the induced action of $G$ on $\hat A\simeq\Prim A$  is free.     Let $s,k\in\P$ with $k\geq 2$ and suppose that
\begin{enumerate}
\item $M_L(\Ind\pi,(\Ind\pi_n))\geq s$;
\item $(\pi_n)_n$ does not converge $k$-times to $\pi$ in $\hat A/G$.
\end{enumerate}
Then
\[
s\leq M_L(\Ind\pi,(\Ind\pi_n))\leq (k-1)M_U(\pi).
\]
In particular, $M_U(\pi)\geq \lceil\frac{s}{k-1}\rceil$.
\end{cor}

\begin{proof}
Immediate from the negation of the (6) $\Longrightarrow$ (2) direction of Theorem~\ref{thm-circle2}.
\end{proof}
%%%%%%%%%%%%%%%%%%%%%%%%%%%%%%%%%%%%%%%%%%%%%%%%%%%%%%%%%%%%%%%%%%%%%%%%%%%%%%%%%%%%%%%%%%%%%%%%%%%%%%%%%%%%%%%%%%%
\section{Characterising Type I properties of $A\rtimes_\alpha G$}\label{sec-3}
%%%%%%%%%%%%%%%%%%%%%%%%%%%%%%%%%%%%%%%%%%%%%%%%%%%%%%%%%%%%%%%%%%%%%%%%%%%%%%%%%%%%%%%%%%%%%%%%%%%%%%%%%%%%%%%%%%%

Recall that it is possible for a crossed product $A\rtimes_\alpha G$ to be Type I even if $A$ is not.

%%%%%%%%%%%%%%%%%%%%%%%%%%%%%%%%%%%%%%%%%%%%%%%%%%%%%%%%%%%%%%%%%%%%%%%%%%%%%%%%%%%%%%%%%%%%%%%%%%%%%%%%%%%%%%%%
\begin{thm}\label{thm-tak3}\textup{(Glimm and Takesaki)}
Let $(A,G,\alpha)$ be a separable $C^*$-dynamical system such that    $A$ is Type I and the induced action of $G$ on $\hat A\simeq \Prim A$ is free.      The following are equivalent:
\begin{enumerate}
\item $A\rtimes_\alpha G$ is Type I;
\item the action of $G$ on  $A$ is smooth (in  the sense that the quotient space $\hat A/G$ is countably separated);
\item the orbits in $\hat A$ are locally closed;
\item for each $\pi\in\hat A$, $s\mapsto s\cdot\pi$ induces a homeomorphism of $G$ onto $G\cdot\pi$.
\end{enumerate}
\end{thm}
%%%%%%%%%%%%%%%%%%%%%%%%%%%%%%%%%%%%%%%%%%%%%%%%%%%%%%%%%%%%%%%%%%%%%%%%%%%%%%%%%%%%%%%%%%%%%%%%%%%%%%%%%%%%%%%%
\begin{proof}
Since $A$ is Type I, $\hat A$ is almost Hausdorff. So by \cite[Theorem~1]{Gli}, items (2)--(4) are equivalent.
That the action of $G$ on $A$ is smooth if and only if $A\rtimes_\alpha G$ is Type I follows by combining \cite[Theorem~6.3]{tak} and \cite[Theorem~8.1]{tak}.
\end{proof}
%%%%%%%%%%%%%%%%%%%%%%%%%%%%%%%%%%%%%%%%%%%%%%%%%%%%%%%%%%%%%%%%%%%%%%%%%%%%%%%%%%%%%%%%%%%%%%%%%%%%%%%%%%%%%%%%%%%

\begin{cor}\label{cor-ccr}
Let $(A,G,\alpha)$ be a separable $C^*$-dynamical system such that   $A$ is liminal and the induced action of $G$ on $\hat A\simeq\Prim A$ is free.      Then $A\rtimes_\alpha G$ is liminal if and only if orbits in $\hat A$ are  closed.
\end{cor}
%%%%%%%%%%%%%%%%%%%%%%%%%%%%%%%%%%%%%%%%%%%%%%%%%%%%%%%%%%%%%%%%%%%%%%%%%%%%%%%%%%%%%%%%%%%%%%%%%%%%%%%%%%%%%%%%%%%

\begin{proof}
Suppose that $A\rtimes_\alpha G$ is liminal.  By assumption, $A$ is  also liminal and the action of $G$ on $\hat A$ is free, so  by Theorem~\ref{thm-tak3} the orbits in $\hat A$ are locally closed and  canonically homeomorphic to $G$. Thus $(A\rtimes_\alpha G)^\wedge\simeq\hat A/G$ by Theorem~\ref{thm-green}.  Since $A\rtimes_\alpha G$ is liminal,  $(A\rtimes_\alpha G)^\wedge$ is $\textrm{T}_1$. So  $\hat A/G$ is $\textrm{T}_1$ as well, and hence the orbits in $\hat A$ are closed.

Conversely, suppose that the orbits in $\hat A$ are closed. By Theorem~\ref{thm-tak3}, $A\rtimes_\alpha G$ is Type I and the orbits are canonically homeomorphic to $G$. Thus $(A\rtimes_\alpha G)^\wedge\simeq\hat A/G$ by Theorem~\ref{thm-green}. Since the orbits are closed, $\hat A/G$ and hence $(A\rtimes_\alpha G)^\wedge$ is $\textrm{T}_1$.  Now $A\rtimes_\alpha G$ is Type I with $\textrm{T}_1$ spectrum and hence  is liminal.
\end{proof}
%%%%%%%%%%%%%%%%%%%%%%%%%%%%%%%%%%%%%%%%%%%%%%%%%%%%%%%%%%%%%%%%%%%%%%%%%%%%%%%%%%%%%%%%%%%%%%%%%%%%%%%%%%%%%%%%
See also \cite[Theorem~6.1]{tak}, \cite[Theorem~3.3]{Goot}, \cite[Theorem~3.1]{W2} and \cite[Propositions~7.32--7.33]{tf^2b} for various results characterising when $A\rtimes_\alpha G$ is liminal or Type I when the induced action of $G$ on    $\Prim A$    is not free.

If $G$ is amenable  we can improve Theorem~\ref{thm-tak3} and Corollary~\ref{cor-ccr} by incorporating the Type I and liminal assumptions on $A$ into the if-and-only-if statements (see Theorem~\ref{thm-better} below); for its proof we need the following lemma. Although we shall apply it only when
G acts freely on   $\Prim A$,    Lemma~\ref{lem-antiliminal} holds for an almost
free action.
%%%%%%%%%%%%%%%%%%%%%%%%%%%%%%%%%%%%%%%%%%%%%%%%%%%%%%%%%%%%%%%%%%%%%%%%%%%%%%%%%%%%%%%%%%%%%%%%%%%%%%%%%%%%%%%%
\begin{lemma}\label{lem-antiliminal}
Let $(A,G,\alpha)$ be a separable $C^*$-dynamical system such that the induced action of $G$ on   $\Prim A$    is almost free. Also assume that  $G$ is amenable. If $A$ is antiliminal then $A\rtimes_\alpha G$ is antiliminal.
\end{lemma}
\begin{proof}
Suppose that $A\rtimes_\alpha G$ is not antiliminal. Let $I$ be the largest Fell ideal of $A\rtimes_\alpha G$ (note $I\neq\{0\}$) by the supposition).  Since   the action of the amenable group $G$ is almost free on $\Prim A$,    there exists a non-zero closed $G$-invariant ideal $J$ of $A$ such that $J\rtimes_\alpha G\subset I$ \cite[Theorem~9]{LN}.  Moreover, there exists a closed $G$-invariant essential ideal $K$ of $A$ such that the induced action of $G$  on   $\Prim K$    is free \cite[Proposition~5(v)]{LN}.  Since $K$ is essential,  $L=K\cap J$ is  a non-zero $G$-invariant ideal of $A$.  Now $G$ acts freely on   $\Prim  L$   
and $L\rtimes_\alpha G$ is a Fell algebra since it is an ideal of $I$.  By Corollary~\ref{cor-A}, $L$ is a Fell algebra.  Since $L\neq \{0\}$, $A$ is not antiliminal.
\end{proof}
%%%%%%%%%%%%%%%%%%%%%%%%%%%%%%%%%%%%%%%%%%%%%%%%%%%%%%%%%%%%%%%%%%%%%%%%%%%%%%%%%%%%%%%%%%%%%%%%%%%%%%%%%%%%%%%%
\begin{thm}\label{thm-better} Let $(A,G,\alpha)$ be a separable $C^*$-dynamical system such that  the induced action of $G$ on   $\Prim A$    is free. Also assume that $G$ is amenable.
\begin{enumerate}
\item  $A\rtimes_\alpha G$ is Type I if and only if  $A$ is Type I and the orbits in $\hat A$ are locally closed.
\item  $A\rtimes_\alpha G$ is liminal  if and only if $A$ is liminal and the orbits in $\hat A$ are closed.
\end{enumerate}
\end{thm}
%%%%%%%%%%%%%%%%%%%%%%%%%%%%%%%%%%%%%%%%%%%%%%%%%%%%%%%%%%%%%%%%%%%%%%%%%%%%%%%%%%%%%%%%%%%%%%%%%%%%%%%%%%%%%%%%
\begin{proof}
(1)  In view of Theorem~\ref{thm-tak3} it remains to show that if $A\rtimes_\alpha G$ is Type I then so is $A$. Suppose that $A$ is not Type I.  Let $J$ be the largest Type I ideal of $A$.  Note that $J$ is $G$-invariant and $(A\rtimes_\alpha G)/(J\rtimes_\alpha G)\cong A/J\rtimes_\alpha G$.  Since $A/J$ is antiliminal and $G$ acts freely on   $\Prim (A/J)$,    $(A/J)\rtimes_\alpha G$ is antiliminal by Lemma~\ref{lem-antiliminal}.  Since  $(A/J)\rtimes_\alpha G$ is isomorphic to a quotient of  $A\rtimes_\alpha G$, and since quotients of Type I algebras are Type I, it follows that $A\rtimes_\alpha G$ is not Type I either.

(2) In view of Corollary~\ref{cor-ccr} it remains to show that if
$A\rtimes_\alpha G$ is liminal then so is $A$.  Suppose that
$A\rtimes_\alpha G$ is liminal but that $A$ is not liminal.  Since
$G$ is amenable and the action of $G$ on   $\Prim A$    is free, $A$ must
be Type I by (1).  So there exists $\pi\in\hat A$ such that
$\pi(A)\supsetneq \K(\H_\pi)$ and hence there exists $\sigma\in \hat
A$ such that $\ker\sigma\supsetneq\ker\pi$. Since $\Ind$ is
continuous and $\sigma\in\overline{\{\pi\}}$,
$\Ind\sigma\in\overline{\{\Ind\pi\}}$.  The spectrum of
$A\rtimes_\alpha G$  is $\textrm{T}_1$, so $\Ind\sigma\simeq\Ind\pi$.
Since $A$ and $A\rtimes_\alpha G$ are Type I and the action of $G$
on $\hat A$ is free, the orbits are locally closed and canonically homeomorphic to $G$ by Theorem~\ref{thm-tak3}. Hence $(A\rtimes_\alpha G)^\wedge \simeq\hat A/G$ by
Theorem~\ref{thm-green}.  Thus $\sigma\in G\cdot\pi$, that is,
$\sigma$ is equivalent to $\pi\circ\alpha_s$ for some $s\in G$.

Note that $\alpha_s(\ker\sigma)\subset\ker\pi\subsetneq\ker\sigma$.  Thus $s\neq e$ and
\[
\ker\sigma\supsetneq\alpha_s(\ker\sigma)\supsetneq\alpha_{s^2}(\ker\sigma)\supsetneq\cdots.
\]
Let
\[
J=\bigcap_{n=1}^\infty\alpha_{s^n}(\ker\sigma).
\]
Since each $\alpha_{s^n}(\ker\sigma)$ is primitive, $J$ is prime and hence primitive by separability of $A$.  Routine calculations show that $\alpha_s(J)=J$.  This contradicts that $G$ acts freely on $\Prim A$. So $A$ must be liminal.
\end{proof}

%%%%%%%%%%%%%%%%%%%%%%%%%%%%%%%%%%%%%%%%%%%%%%%%%%%%%%%%%%%%%%%%%%%%%%%%%%%%%%%%%%%%%%%%%%%%%%%%%%%%%%%%%%%%%%%%%%%
Let $(A, G,\alpha)$ be a separable $C^*$-dynamical system such that  the induced action of $G$ on $\hat A$ is free. Raeburn and Rosenberg showed that if $A$ has continuous trace and $G$ acts properly on $\hat A$, then $A\rtimes_\alpha G$  has continuous trace \cite[Theorem~1.1(3)]{RR}; Olesen and Raeburn proved the converse in \cite[Theorem~3.1]{OR} for  abelian $G$ assuming $A$ has continuous trace. Deicke used non-abelian duality to extend the Olesen-Raeburn result to non-abelian groups in \cite[Theorem~4.7]{D}.  Thus, provided $A$ has continuous trace, $A\rtimes_\alpha G$ has continuous trace if and only if the induced action of $G$ on $\hat A$ is proper.

An interesting test question for our techniques is whether we can
use them to recover the  Deicke-Olesen-Raeburn-Rosenberg result.
Indeed, not only can we recover their theorem, we can improve it by
removing the standing assumption that $A$ has continuous trace and
incorporating it into the if-and-only-if statement (see
Theorem~\ref{thm-ct} below).  The next theorem is a key ingredient
in the proof of the new part of Theorem~\ref{thm-ct}. If
$A\rtimes_\alpha G$ has continuous trace then $A$ must be a Fell
algebra by Corollary~\ref{cor-A};  Theorem~\ref{thm-1a} will help to
establish that   $\hat A\simeq\Prim A$    must actually be  Hausdorff, and thus $A$ is
in fact a continuous-trace $C^*$-algebra.

%%%%%%%%%%%%%%%%%%%%%%%%%%%%%%%%%%%%%%%%%%%%%%%%%%%%%%%%%%%%%%%%%%%%%%%%%%%%%%%%%%%%%%%%%%%%%%%%%%%%%%%%%%%%%%%%%%%

\begin{thm}\label{thm-1a}
Let $(A,G,\alpha)$ be a separable $C^*$-dynamical system such that  the induced action
of $G$ on   $\Prim A$    is free. Suppose $\sigma_0, \sigma_1,\dots, \sigma_k$  are distinct
points of $\hat A$ such that $\sigma_0\in\overline{G\cdot\sigma_i}$ for $1\leq i\leq k$ and  that $(\pi_j)_j$ is a sequence in $\hat A$ such that $\pi_j\to\sigma_i$ for $0\leq i\leq k$.  Then $M_L(\Ind\sigma_0,(\Ind\pi_j))\geq k+1$.
\end{thm}
%%%%%%%%%%%%%%%%%%%%%%%%%%%%%%%%%%%%%%%%%%%%%%%%%%%%%%%%%%%%%%%%%%%%%%%%%%%%%%%%%%%%%%%%%%%%%%%%%%%%%%%%%%%%%%%%%%%

We immediately obtain the following corollary (which will be illustrated by Example~\ref{ex-fell} below):

%%%%%%%%%%%%%%%%%%%%%%%%%%%%%%%%%%%%%%%%%%%%%%%%%%%%%%%%%%%%%%%%%%%%%%%%%%%%%%%%%%%%%%%%%%%%%%%%%%%%%%%%%%%%%%%%%%%
\begin{cor}\label{cor-1}
Let $(A,G,\alpha)$ be a separable $C^*$-dynamical system such that  the induced action of $G$ on   $\Prim A$    is free. Suppose $\sigma_0$ and $\sigma_1$ are distinct points of $\hat A$ which cannot be separated by disjoint open sets, and $\sigma_0\in\overline{G\cdot\sigma_1}$.  Then $M_U(\Ind\sigma_0)\geq 2$.
\end{cor}

%%%%%%%%%%%%%%%%%%%%%%%%%%%%%%%%%%%%%%%%%%%%%%%%%%%%%%%%%%%%%%%%%%%%%%%%%%%%%%%%%%%%%%%%%%%%%%%%%%%%%%%%%%%%%%%%%%%
To prove Theorem~\ref{thm-1a}, we need two lemmas. The first  reworks ideas from the proof of \cite[Theorem~1]{A-pems} into a form that will be convenient. The second is based on the Gram-Schmidt process.

%%%%%%%%%%%%%%%%%%%%%%%%%%%%%%%%%%%%%%%%%%%%%%%%%%%%%%%%%%%%%%%%%%%%%%%%%%%%%%%%%%%%%%%%%%%%%%%%%%%%%%%%%%%%%%%%%%%%%%%%%%%%%%

\begin{lemma}\label{lem-pems2}
Let $A$ be a $C^*$-algebra, $(\pi_\lambda)_{\lambda\in\Lambda}$ be a net in $\hat A$ and suppose that $\sigma_1,\dots, \sigma_k$ are distinct limits of $(\pi_\lambda)_\lambda$ in $\hat A$. For $1\leq i\leq k$, let $\phi_i$ be a pure state of $A$ associated with $\sigma_i$ and let $N_i$ be a $w^*$-open neighbourhood of $0$ in $A^*$. Let $\lambda_0\in\Lambda$ and $\epsilon>0$.  Then there exist $\lambda\in\Lambda$ and unit vectors $\xi_1,\dots \xi_k$ in the Hilbert space for $\pi_\lambda$  such that $\lambda\geq \lambda_0$ and
\begin{gather*}|\langle\xi_i\,,\,\xi_j\rangle|<\epsilon\quad\quad(1\leq i<j\leq k);\\
\langle\pi_\lambda(\cdot)\xi_i\,,\,\xi_i\rangle\in\phi_i+N_i\quad\quad(1\leq i\leq k).
\end{gather*}
\end{lemma}
%%%%%%%%%%%%%%%%%%%%%%%%%%%%%%%%%%%%%%%%%%%%%%%%%%%%%%%%%%%%%%%%%%%%%%%%%%%%%%%%%%%%%%%%%%%%%%%%%%%%%%%%%%%%%%%%%%%%%%%
\begin{proof}
Let $N_0=\cap_{i=1}^k N_i$, a $w^*$-open neighbourhood of $0$ in $A^*$.  Since $\phi_i$ and $\phi_j$ ($i\neq j$) are inequivalent, the transition probability $\langle\phi_i\,,\,\phi_j\rangle=0$.  Since there are only a finite number of such pairs $\{i,j\}$, it follows from a restricted continuity property for transition probabilities (see \cite[Remarks following Corollary~2.4, item 2]{ASh}) that there exists a $w^*$-open neighbourhood $N$ of $0$ in $A^*$ such that $N\subset N_0$ and, for all pure states $\psi_i\in\phi_i+N$,
\begin{equation}\label{eq-lem-pems}
\langle\psi_i\,,\,\psi_j\rangle<\epsilon^2 \quad(i\neq j).
\end{equation}
Let $V_i$ be the canonical image of $(\phi_i+N)\cap P(A)$ in $\hat A$; in particular $V_i$ is an open neighbourhood of $\sigma_i\ (1\leq i\leq k)$.  There exists $\lambda\geq\lambda_0$ such that $\pi_\lambda\in\cap_{i=1}^k V_i$.  Hence there exist unit vectors $\xi_i$ such that
\[
\langle\pi_\lambda(\cdot)\xi_i\,,\,\xi_i\rangle\in\phi_i+N\subset \phi_i+N_i
\]
for $1\leq i\leq k$.  If $i\neq j$ then $|\langle\xi_i\,,\,\xi_j\rangle|^2<\epsilon^2$ by \eqref{eq-lem-pems}, and hence $|\langle\xi_i\,,\,\xi_j\rangle|<\epsilon$.
\end{proof}
%%%%%%%%%%%%%%%%%%%%%%%%%%%%%%%%%%%%%%%%%%%%%%%%%%%%%%%%%%%%%%%%%%%%%%%%%%%%%%%%%%%%%%%%%%%%%%%%%%%%%%%%%%%%%%%%%%%
\begin{lemma}\label{lem-gs}
Let $k\geq 1$, $A$ a $C^*$-algebra and $\phi\in S(A)$. Suppose that there is a net $(\pi_\lambda)_{\lambda\in\Lambda}$ of representations of $A$ on Hilbert spaces $\H_\lambda$ and unit vectors $\xi_\lambda^{(i)}$\ $(1\leq i\leq k)$ in each $\H_\lambda$ such that
\begin{gather}
\langle\pi_\lambda(\cdot)\xi_\lambda^{(i)}\,,\,\xi_\lambda^{(i)}\rangle\to\phi\quad(1\leq i\leq k)\quad\text{and}\label{eq-lem-gs-1}\\
\langle \xi_\lambda^{(i)}\,,\, \xi_\lambda^{(j)}\rangle \to 0\quad (1\leq i<j\leq k).\label{eq-lem-gs-2}
\end{gather}
Then there exists $\lambda_1\in\Lambda$ such that for each $\lambda\geq \lambda_1$ there is an orthonormal set $\{\eta_\lambda^{(1)},\dots,\eta_\lambda^{(k)}\}\subset \H_\lambda$ satisfying
\begin{equation}\label{eq-lem-gs-3}
\langle\pi_\lambda(\cdot)\eta_\lambda^{(i)}\,,\,\eta_\lambda^{(i)}\rangle\to\phi\quad\text{and}\quad \|\xi_\lambda^{(i)}-\eta_\lambda^{(i)}\|\to 0\quad(1\leq i\leq k).
\end{equation}
\end{lemma}
%%%%%%%%%%%%%%%%%%%%%%%%%%%%%%%%%%%%%%%%%%%%%%%%%%%%%%%%%%%%%%%%%%%%%%%%%%%%%%%%%%%%%%%%%%%%%%%%%%%%%%%%%%%%%%%%%%%
\begin{proof}
We use induction on $k$.  When $k=1$ we may simply define $\eta_\lambda^{(1)}=\xi_\lambda^{(1)}$ for all $\lambda$.  Now suppose that $n\geq 1$, that the result is true when $k=n$, and that we have the data in \eqref{eq-lem-gs-1} and \eqref{eq-lem-gs-2} for $1\leq i\leq n+1$ and $1\leq i<j\leq n+1$.  By the induction hypothesis there exists $\lambda_0\in\Lambda$ such that for each $\lambda\geq \lambda_0$ there is an orthonormal set $\{\eta_\lambda^{(1)},\dots, \eta_\lambda^{(n)}\}$ such that \eqref{eq-lem-gs-3} holds for $1\leq i\leq n$.  For each $\lambda\geq \lambda_0$, let
\[
v_\lambda=\xi^{n+1}_\lambda-\sum_{i=1}^n\langle\xi_\lambda^{(n+1)}\,,\,\eta_\lambda^{(i)}\rangle\eta_\lambda^{(i)}
\]
in $\H_\lambda$.  Note that, since $\langle \xi_\lambda^{(n+1)}\,,\, \xi_\lambda^{(i)}\rangle \to 0$  and $\|\xi_\lambda^{(i)}-\eta_{\lambda}^{(i)}\|\to 0\|$ for $1\leq i\leq n$, we have $\langle\xi_\lambda^{(n+1)}\,,\,\eta_\lambda^{(i)}\rangle\to 0$ for $1\leq i\leq n$.  Thus $\|v_\lambda-\xi^{(n+1)}_\lambda\|\to 0$ and so $\|v_\lambda\|^2\to 1$ and $\langle\pi_\lambda(\cdot)v_\lambda\,,\,v_\lambda\rangle\to\phi$.  There exists $\lambda_1\geq \lambda_0$ such that $\|v_\lambda\|>0$ for all $\lambda>\lambda_1$, and so we complete the inductive step by defining $\eta_\lambda^{(n+1)}=v_\lambda/\|v_\lambda\|$ for $\lambda\geq \lambda_1$.
\end{proof}

%%%%%%%%%%%%%%%%%%%%%%%%%%%%%%%%%%%%%%%%%%%%%%%%%%%%%%%%%%%%%%%%%%%%%%%%%%%%%%%%%%%%%%%%%%%%%%%%%%%%%%%%%%%%%%%%%%%
%%%%%%%%%%%%%%%%%%%%%%%%%%%%%%%%%%%%%%%%%%%%%%%%%%%%%%%%%%%%%%%%%%%%%%%%%%%%%%%%%%%%%%%%%%%%%%%%%%%%%%%%%%%%%%%%%%%
%%%%%%%%%%%%%%%%%%%%%%%%%%%%%%%%%%%%%%%%%%%%%%%%%%%%%%%%%%%%%%%%%%%%%%%%%%%%%%%%%%%%%%%%%%%%%%%%%%%%%%%%%%%%%%%%%%%

\begin{proof}[Proof of Theorem~\ref{thm-1a}]   Since the induced action of $G$ on $\Prim A$ is free, $\Ind\sigma$ is irreducible for all $\sigma\in\hat A$ by by \cite[Proposition~1.7]{EW}.    We will use Lemma~\ref{lem-j}(2) to show that $M_L(\Ind\sigma_0,(\Ind\pi_j))\geq k+1$.  Note that our data ($\pi_j\to\sigma_i$ for $0\leq i\leq k$ and $\sigma_0\in\overline{G\cdot\sigma_i}$ for $1\leq i\leq k$) are unchanged by passing to a subsequence; so we pass to some subsequence and relabel to continue working with $(\pi_j)_j$.

Let $\phi$ be a pure state of $A$ associated to $\sigma_0$ and let $\xi_\phi$ be the corresponding GNS vector.  Since $A$ is separable, there exists a decreasing sequence $(N_n)_{n\geq 1}$ of basic $w^*$-open neighbourhoods of $0$ in $A^*$.

Temporarily fix $n\geq 1$ and consider the canonical image in $\hat A$ of the open neighbourhood $(\phi+\frac{1}{2}N_n)\cap P(A)$ of $\phi$ in $P(A)$.  This image is an open neighbourhood of $\sigma_0$ in $\hat A$ and so contains $g_n^{(i)}\cdot\sigma_i$ for some $g_n^{(i)}\in G$ for $1\leq i\leq k$. It follows that for each $i\in\{1,\dots,k\}$ there is a unit vector $u_n^{(i)}$ in the Hilbert space of $\sigma_i$ such that
\[
\psi_n^{(i)}:=\langle(g_n^{(i)}\cdot\sigma_i)(\cdot )u_n^{(i)}\,,\, u_n^{(i)}\rangle\in\phi+\frac{1}{2}N_n.
\]
For $1\leq i\leq k$, let
\begin{equation}\label{eq-thm1a-0}
M_n^{(i)}=\big\{\rho\circ\alpha_{g_n^{(i)}}:\rho\in \frac{1}{2}N_n\big\},
\end{equation}
a $w^*$-open neighbourhood of $0$ in $A^*$, and note that
\[
\psi_n^{(i)}\circ\alpha_{g_n^{(i)}}=\langle\sigma_i(\cdot )u_n^{(i)}\,,\, u_n^{(i)}\rangle
\]
is a pure state associated with $\sigma_i$.

We will inductively construct a sequence $j_1<j_2<\cdots<j_n<\cdots$ such that, for each $n$, there exist $k+1$ unit vectors $\xi_n^{(i)}$ $(0\leq i\leq k)$ in the Hilbert space of $\pi_{j_n}$  satisfying
\begin{gather}
\phi_n:=\langle\pi_{j_n}(\cdot)\xi_n^{(0)}\,,\,\xi_n^{(0)}\rangle\in\phi+\frac{1}{2}N_n,\label{eq-thm1a-1}\\
\langle\pi_{j_n}(\cdot)\xi_n^{(i)}\,,\,\xi_n^{(i)}\rangle\in\psi_n^{(i)}\circ\alpha_{g_n^{(i)}}+M_n^{(i)} \quad \quad(1\leq i\leq j),\label{eq-thm1a-2}\\
|\langle\xi_n^{(i)}\,,\,\xi_n^{(j)}\rangle|<\frac{1}{n} \quad\quad(0\leq i<j\leq k).\label{eq-thm1a-3}
\end{gather}
When $n=1$ we apply Lemma~\ref{lem-pems2} to the sequence $(\pi_j)_j$ and the $k+1$ inequivalent states $\phi$ and $\psi_1^{(i)}\circ\alpha_{g_1^{(i)}}\ (1\leq i\leq k)$ with $\epsilon=1$ and neighbourhoods $\frac{1}{2}N_1$ and $M_1^{(i)}\ (1\leq i\leq k)$, to get $j_1>0$ and unit vectors $\xi_1^{(i)}\ (0\leq i\leq k)$ on the Hilbert space of $\pi_{j_1}$ such that \eqref{eq-thm1a-1}, \eqref{eq-thm1a-2} and \eqref{eq-thm1a-3} hold with $n=1$.

Now assume that $j_1<j_2<\cdots<j_r$ and corresponding vectors $\xi_n^{(i)}\ (0\leq i\leq k)$ in the Hilbert space of $\pi_{j_n}$ have been chosen such that \eqref{eq-thm1a-1}, \eqref{eq-thm1a-2} and \eqref{eq-thm1a-3} hold for each $n\in\{1,\dots r\}$.  We apply Lemma~\ref{lem-pems2} to the sequence $(\pi_j)_j$ and the $k+1$ inequivalent states $\phi$ and $\psi_{r+1}^{(i)}\circ\alpha_{g_{r+1}^{(i)}}\ (1\leq i\leq k)$ with $\epsilon=1/(r+1)$ and neighbourhoods $\frac{1}{2}N_{r+1}$ and $M_{r+1}^{(i)}\ (1\leq i\leq k)$, to get $j_{r+1}>j_r$ and unit vectors $\xi_{r+1}^{(i)}\ (0\leq i\leq k)$ in the Hilbert space of $\pi_{j_{r+1}}$ such that \eqref{eq-thm1a-1}, \eqref{eq-thm1a-2} and \eqref{eq-thm1a-3} hold with $n=r+1$.

By our choice of $M_n^{(i)}$ at \eqref{eq-thm1a-0}, for each $n\geq 1$,
\[
\rho_n^{(i)}:=\langle(g_n^{(i)}\cdot\pi_{j_n})(\cdot )\xi_n^{(i)}\,,\, \xi_n^{(i)}\rangle\in\psi_n^{(i)}+\frac{1}{2}N_n\subset\phi+N_n\quad\quad (1\leq i\leq k).
\]

Now let $n$ vary and recall that the sequence $(N_n)_n$ is decreasing, so that
\begin{equation*}
%\label{eq-0}
\phi_n\to\phi,\quad\rho_n^{(i)}\to\phi\quad(1\leq i\leq k)\quad\text{and}
\quad\langle\xi_n^{(i)}\,,\,\xi_n^{(j)}\rangle\to 0\quad (0\leq i<j\leq k)
\end{equation*}
as $n\to\infty$.

Let $W$ be a compact symmetric neighbourhood of $e$ in $G$ and let
\[
\tilde\xi_\phi=\nu(W)^{-1/2}\chi_W\otimes\xi_\phi.
\]
in $L^2(G,\nu)\otimes\H_\pi$.  Then $\tilde\xi_\phi$ is a unit vector and, for $n\geq 1$, we similarly define unit vectors
\[
\tilde\xi_n^{(0)}=\nu(W)^{-1/2}\chi_W\otimes\xi_n^{(0)}\quad\text{and}\quad\tilde\xi_n^{(i)}=\nu(W)^{-1/2}\chi_{Wg_n^{(i)}}\otimes\xi_n^{(i)}\quad(1\leq i\leq k)
\]
in $L^2(G,\nu)\otimes\H_{\pi_{j_n}}$.  Note that $\langle\tilde\xi_n^{(i)}\,,\,\tilde\xi_n^{(j)}\rangle\to_n 0$ as $n\to\infty$ for $0\leq i<j\leq k$.
We aim to show that
\begin{equation}
\langle\Ind\pi_{j_n}(\cdot )\tilde\xi_n^{(i)}\,,\, \tilde\xi_n^{(i)}\rangle\to\langle\Ind\pi(\cdot )\tilde\xi_\phi\,,\, \tilde\xi_\phi\rangle\quad\quad(0\leq i\leq k).\label{eq-thm1a-4}\\
\end{equation}
Since all the functionals involved have unit norm, it suffices to check \eqref{eq-thm1a-4}  on $C_c(G,A)$, a dense subalgebra of $A\rtimes_\alpha G$. For $f\in C_c(G,A)$, routine calculations as in the proof of Theorem~\ref{thm-a} yield that
\begin{align*}
\langle\Ind\pi_{j_n}(f)\tilde\xi_n^{(0)}\,,\,\tilde\xi_n^{(0)}\rangle &
=\int_{t\in W}\int_{s\in W}\Delta(st^{-1})^{-1/2}\langle\pi_{j_n}(\alpha_t^{-1}(f(ts^{-1})))\xi_n^{(0)}\,,\,\xi_n^{(0)}\rangle\, d\nu(s)\, d\nu(t)\\
&=\int_{t\in W}\int_{s\in W}\Delta(st^{-1})^{-1/2} \phi_n(\alpha_t^{-1}(f(ts^{-1}))) \, d\nu(s)\, d\nu(t)
\end{align*}
and
\begin{align*}
\langle\Ind\pi(f)\tilde\xi_\phi\,,\,\tilde\xi_\phi\rangle
&=\int_{t\in W}\int_{s\in W}\Delta(st^{-1})^{-1/2}\langle\pi(\alpha_t^{-1}(t(ts^{-1})))\xi_\phi\,,\,\xi_\phi\rangle\, d\nu(s)\, d\nu(t)\\
&=\int_{t\in W}\int_{s\in W}\Delta(st^{-1})^{-1/2} \phi(\alpha_t^{-1}(f(ts^{-1}))) \, d\nu(s)\, d\nu(t).
\end{align*}
Since $\phi_n\to\phi$,  Equation~\eqref{eq-thm1a-4} for $i=0$  now follows using the Bounded Convergence Theorem.

On the other hand, for $1\leq i\leq k$,
\begin{align*}
\langle&\Ind\pi_{j_n}(f)\tilde\xi_n^{(i)}\,,\, \tilde\xi_n^{(i)}\rangle\\
&=
\nu(W)^{-1}\int_{t\in Wg_n^{(i)}}\int_{s\in Wg_n^{(i)}}\Delta(st^{-1})^{-1/2}\langle\pi_{j_n}(\alpha_t^{-1}(f(ts^{-1})))\xi_n^{(i)}\,,\,\xi_n^{(i)}\rangle\, d\nu(s)\, d\nu(t)\\
&=\nu(W)^{-1}\int_{t\in Wg_n^{(i)}}\int_{p\in W}\Delta(pg_n^{(i)}t^{-1})^{-1/2}\langle\pi_{j_n}(\alpha_t^{-1}(f(t(g_n^{(i)})^{-1}p^{-1})))\xi_n^{(i)}\,,\,\xi_n^{(i)}\rangle\, d\nu(p)\, d\nu(t)\\
&=\nu(W)^{-1}\int_{q\in W}\int_{p\in W}\Delta(pq^{-1})^{-1/2}\langle\pi_{j_n}\circ\alpha_{g_n^{(i)}}^{-1}(\alpha_q^{-1}(f(qp^{-1})))\xi_n^{(i)}\,,\,\xi_n^{(i)}\rangle\, d\nu(p)\, d\nu(q)\\
&=\nu(W)^{-1}\int_{q\in W}\int_{p\in W}\Delta(pq^{-1})^{-1/2}\rho_n^{(i)}(\alpha_q^{-1}(f(qp^{-1})))\, d\nu(p)\, d\nu(q).
\end{align*}
Since $\rho_n^{(i)}\to\phi$, Equation~\eqref{eq-thm1a-4} for $1\leq i\leq k$ now follows using the Bounded Convergence Theorem. Since also $\langle\tilde\xi_n^{(i)}\,,\,\tilde\xi_n^{(j)}\rangle\to 0$ if $i\neq j$, it now follows from Lemma~\ref{lem-gs} that there exists $n_0$ such that for $n\geq n_0$ there is an orthonormal set $\{\eta_n^{(i)}:0\leq i\leq k\}\subset L^2(G,\nu)\otimes\H_{\pi_{j_n}}$ such that
\[
\langle\Ind\pi_{j_n}(\cdot)\eta_n^{(i)}\,,\,\eta_n^{(i)}\rangle \to \langle\Ind\pi(\cdot )\tilde\xi_\phi\,,\, \tilde\xi_\phi\rangle\quad(0\leq i\leq k).
\]
Thus $M_L(\Ind\sigma_0,(\Ind\pi_j))\geq k+1$  by Lemma~\ref{lem-j}(2).
\end{proof}

%%%%%%%%%%%%%%%%%%%%%%%%%%%%%%%%%%%%%%%%%%%%%%%%%%%%%%%%%%%%%%%%%%%%%%%%%%%%%%%%%%%%%%%%%%%%%%%%%%%%%%%%%%%%%%%%%%%
It is tempting to conjecture that Theorem~\ref{thm-a} could be used to derive the result in Theorem~\ref{thm-1a}, but Example~\ref{ex-fell} below illustrates that this is not so.
%%%%%%%%%%%%%%%%%%%%%%%%%%%%%%%%%%%%%%%%%%%%%%%%%%%%%%%%%%%%%%%%%%%%%%%%%%%%%%%%%%%%%%%%%%%%%%%%%%%%%%%%%%%%%%%%%%%

\medskip

The following theorem is mainly due to the combined efforts of
Deicke, Olesen, Raeburn and Rosenberg (see the discussion earlier in
this section). Our contribution to the theorem is that if the action
of $G$ on   $\Prim A$    is free and $A\rtimes_\alpha G$ has continuous
trace then $A$ too must have continuous trace. We also show how the
other parts of the theorem can be obtained by using some of our
earlier results.
%%%%%%%%%%%%%%%%%%%%%%%%%%%%%%%%%%%%%%%%%%%%%%%%%%%%%%%%%%%%%%%%%%%%%%%%%%%%%%%%%%%%%%%%%%%%%%%%%%%%%%%%%%%%%%%%%%%
\begin{thm}\label{thm-ct}
Let $(A,G,\alpha)$ be a separable $C^*$-dynamical system such that the induced action of $G$ on   $\Prim A$    is free. Then $A\rtimes_\alpha G$ has continuous trace if and only if $A$ has continuous trace and the action of $G$ on $\hat A$ is proper.
\end{thm}

%%%%%%%%%%%%%%%%%%%%%%%%%%%%%%%%%%%%%%%%%%%%%%%%%%%%%%%%%%%%%%%%%%%%%%%%%%%%%%%%%%%%%%%%%%%%%%%%%%%%%%%%%%%%%%%%%%%
\begin{proof}
Suppose that $A\rtimes_\alpha G$ has continuous trace.
  Since the action of $G$ on $\Prim A$ is free,  $A$ is a Fell algebra and $\hat A\simeq\Prim A$ is a Cartan $G$-space by Corollary~\ref{cor-A}.   
We start by showing  that $\hat A$ is Hausdorff.  Suppose not, that is suppose there exist distinct points $\pi$ and $\sigma$ in $\hat A$ which cannot be separated by disjoint open sets.  Then there exists a sequence $(\pi_n)$ in $\hat A$ such that $\pi_n\to\sigma$ and $\pi_n\to\pi$. Then $\Ind\pi_n\to\Ind\sigma$ and $\Ind\pi_n\to\Ind\pi$ in $(A\rtimes_\alpha G)^\wedge$.  But $(A\rtimes_\alpha G)^\wedge$ is Hausdorff, so $\Ind\sigma\simeq\Ind\pi$.  Since both $A$ and $A\rtimes_\alpha G$ are Type I, the orbits in $\hat A$ are locally closed and canonically isomorphic to $G$ by Theorem~\ref{thm-tak3}. Thus $\Ind:\hat A\to (A\rtimes_\alpha G)^\wedge$ induces a homeomorphism of $\hat A/G$ onto $(A\rtimes_\alpha G)^\wedge$ by Theorem~\ref{thm-green}. It follows that $\sigma\in G\cdot\pi$.  But now $M_U(\Ind\sigma)\geq 2$ by Corollary~\ref{cor-1}, which is impossible since $A\rtimes_\alpha G$ has continuous trace.  Thus $\hat A$ is Hausdorff and hence $A$ has continuous trace.

Since $\hat A$ is a locally compact Hausdorff space it is completely regular. Moreover, $(G,\hat A)$ is Cartan and $\hat A/G$ is Hausdorff (since it is homeomorphic to $(A\rtimes_\alpha G)^\wedge$), so $G$ acts properly on $\hat A$ by \cite[Theorem~1.2.9]{palais}.
%\footnote{Alternatively, we could have appealed to the combination of \cite[Theorem~3.1]{OR} and \cite[Theorem~4.7]{D} here.}

Conversely, assume that $A$ has continuous trace and that $G$ acts properly on $\hat A$.  Since   $\hat A\simeq\Prim A$    is Hausdorff and the action is proper, the orbits are closed.  So  the orbits are canonically homeomorphic to $G$ by Theorem~\ref{thm-tak3} and then $(A\rtimes_\alpha G)^\wedge\simeq \hat A/G$ by Theorem~\ref{thm-green}.  Since the action of $G$ on $\hat A$  is proper and $\hat A$ is locally compact Hausdorff (hence completely regular), $\hat A/G$ is Hausdorff by \cite[Theorem~1.2.9]{palais}.  Thus $(A\rtimes_\alpha G)^\wedge$ is Hausdorff as well.

It remains to show that $A\rtimes_\alpha G$ is a Fell algebra.  Suppose it is not, that is suppose that there exists $\sigma\in (A\rtimes_\alpha G)^\wedge$ such that $M_U(\sigma)\geq 2$. By \cite[Lemma~1.2]{AK} there exists a sequence $(\sigma_n)_n$ in $(A\rtimes_\alpha G)^\wedge$ such that
$
2\leq M_U(\sigma)=M_L(\sigma, (\sigma_n)).
$
Since $(A\rtimes_\alpha G)^\wedge\simeq \hat A/G$, there exist $\pi\in\hat A$ and $(\pi_n)_n\subset\hat A$ such that $\sigma=\Ind\pi$ and $\sigma_n=\Ind\pi_n$.
By the (6) $\Longrightarrow$ (2) direction of Theorem~\ref{thm-circle2} (note that $A$ has continuous trace so $m=1$) it follows that
\[
M_L(\Ind\epsilon_\pi,(\Ind(\epsilon_{\pi_n}))\geq 2
\]
in $(C_0(\hat A)\rtimes_\lt G)^\wedge$.  But now $M_U(\Ind\epsilon_\pi,(\Ind(\epsilon_{\pi_n}))\geq 2$ and hence $C_0(\hat A)\rtimes_\lt G$ does not have continuous trace, contradicting that $G$ acts freely and properly (see \cite[Theorem~17]{green1}).  Thus $A\rtimes_\alpha G$ is a Fell algebra with Hausdorff spectrum and hence has continuous trace.
\end{proof}

%%%%%%%%%%%%%%%%%%%%%%%%%%%%%%%%%%%%%%%%%%%%%%%%%%%%%%%%%%%%%%%%%%%%%%%%%%%%%%%%%%%%%%%%%%%%%%%%%%%%%%%%%%%%%%%%%%%
If $X$ is completely regular and $(G,X)$ is Cartan then all orbits are closed in $X$ by \cite[Proposition~1.14]{palais} and $s\mapsto s\cdot x$ is an open map of $G$ onto the orbit $G\cdot x$ by \cite[Lemma on p. 298]{palais}.  If $X=\hat A$ is the spectrum of a $C^*$-algebra $A$ then it can easily fail to be completely regular, even if $A$ is a Fell algebra; in this situation extra data is needed to characterise when $A\rtimes_\alpha G$ is a Fell algebra.
If $A=C_0(X)$ then $\hat A$ is Hausdorff and Theorem~\ref{thm-fell} reduces to  \cite[Theorem~1]{aH-Fell}.
%%%%%%%%%%%%%%%%%%%%%%%%%%%%%%%%%%%%%%%%%%%%%%%%%%%%%%%%%%%%%%%%%%%%%%%%%%%%%%%%%%%%%%%%%%%%%%%%%%%%%%%%%%%%%%%%%%%

\begin{thm}\label{thm-fell}
Let $(A,G,\alpha)$ be a separable $C^*$-dynamical system such that  the induced action of $G$ on    $\Prim A$    is free. Then $A\rtimes_\alpha G$ is a Fell algebra if and only if
\begin{enumerate}
\item $A$ is a Fell algebra; and
\item $(G,\hat A)$ is Cartan; and
\item every $\pi\in\hat A$ has a $G$-invariant open Hausdorff neighbourhood.
\end{enumerate}
\end{thm}

%%%%%%%%%%%%%%%%%%%%%%%%%%%%%%%%%%%%%%%%%%%%%%%%%%%%%%%%%%%%%%%%%%%%%%%%%%%%%%%%%%%%%%%%%%%%%%%%%%%%%%%%%%%%%%%%%%%

\begin{proof}
Suppose that $A\rtimes_\alpha G$ is a Fell algebra.  Then $A$ is a Fell algebra and  $(G, \hat A)$ is  Cartan  by Corollary~\ref{cor-A}. Since both $A$ and $A\rtimes_\alpha G$ are Type I, the orbits in $\hat A$ are locally closed  and canonically homeomorphic to $G$ by Theorem~\ref{thm-tak3}. Thus  $\Ind:\hat A\to (A\rtimes_\alpha G)^\wedge$ induces a homeomorphism of $\hat A/G$ onto $(A\rtimes_\alpha G)^\wedge$ by Theorem~\ref{thm-green}.  Let $\pi\in\hat A$.  Since $A\rtimes_\alpha G$ is a Fell algebra, $\Ind\pi$ has an open Hausdorff neighbourhood $W$ \cite[Corollary~3.4]{ASo}.  Then $\Ind^{-1}(W)$ is an open  $G$-invariant neighbourhood of $\pi$ in $\hat A$, and hence equals $\hat J$ for some closed two-sided   $G$-invariant    ideal $J$ of $A$.  We claim that $(J\rtimes_\alpha G)^\wedge\subset W$.  To see this, suppose that $\rho\in (A\rtimes_\alpha G)^\wedge\setminus W$. Note that $\rho=\Ind\sigma$ for some $\sigma\in\hat A$.  Since $\rho\notin W$, $\sigma\notin \hat J$. Hence $\sigma(J)=\{0\}$ and so $\rho(J\rtimes_\alpha G)=\{0\}$, that is $\rho\notin (J\rtimes_\alpha G)^\wedge$. Thus $(J\rtimes_\alpha G)^\wedge\subset W$ as claimed.
Now $J\rtimes_\alpha G$ is a Fell algebra with Hausdorff spectrum, hence has continuous trace.    Since $G$ acts freely on $\Prim J$,   $J$ has continuous trace  by Theorem~\ref{thm-ct}.     Thus $\hat J$ is an open  $G$-invariant Hausdorff neighbourhood of $\pi$.

Conversely, assume (1), (2) and (3).  We start by showing that (2)
and (3) together imply that the orbits in $\hat A$ are closed.  Let
$\pi\in\hat A$ and let $\sigma$ be in the closure of $G\cdot\pi$.
So there exists a sequence $(s_n)_n$ in $G$ such that
$s_n\cdot\pi\to\sigma$.  Let $U$ be a wandering neighbourhood of
$\sigma$.  There exists $n_0$ such that $s_n\cdot\pi\in U$ for
$n\geq n_0$.  Since $(s_{n_0}s_n^{-1})\cdot
(s_n\cdot\pi)=s_{n_0}\cdot\pi\in U$ we have
$s_{n_0}s_n^{-1}\in\{t\in G:t\cdot U\cap U\neq\emptyset\}$ whenever
$n\geq n_0$.  Since $U$ is wandering there exists a subsequence
$(s_{n_i})_i$ of $(s_n)_n$ which converges to some $s\in G$.  Now
$s_{n_i}\cdot\pi\to\sigma$ and $s_{n_i}\cdot\pi\to s\cdot\pi$.  Let
$V$ be a $G$-invariant Hausdorff neighbourhood of $\sigma$.  Since
$s_{n_i}\cdot\pi\in V$ eventually, $s\cdot\pi\in V$ by the
$G$-invariance.  Now $\sigma=s\cdot\pi\in G\cdot \pi$, so
$G\cdot\pi$ is closed.
%\footnote{Similarly, (2) and (3) imply that $s\mapsto s\cdot\pi$ is an open map of $G$ onto $G\cdot\pi$.}

Let $\sigma\in(A\rtimes_\alpha G)^\wedge$; we will show that $M_U(\sigma)=1$.  Since the orbits in $\hat A$ are closed and $A$ is Type I, the orbits   are canonically homeomorphic to $G$ by Theorem~\ref{thm-tak3}. So $(A\rtimes_\alpha G)^\wedge\simeq\hat A/G$ by Theorem~\ref{thm-green}, and hence $\sigma=\Ind\pi$ for some $\pi\in\hat A$.  By (3) there exists  a $G$-invariant open Hausdorff neighbourhood $V$ of $\pi$ in $\hat A$.  By (2) and shrinking $V$ if necessary, we may assume that $(G,V)$ is proper (the proof of \cite[Poposition~1.2.4]{palais} does not require complete regularity of $\hat A$). Let $J$ be the $G$-invariant closed two-sided ideal of $A$ such that $V=\hat J$. Since $A$ is a Fell algebra by (1), so is $J$, and since $J$ has Hausdorff spectrum it has continuous trace.  By Theorem~\ref{thm-ct} $J\rtimes_\alpha G$ has continuous trace.  Since $\sigma\in  (J\rtimes_\alpha G)^\wedge$, $M_U(\sigma)=1$ in $(J\rtimes_\alpha G)^\wedge$ and hence in $(A\rtimes_\alpha G)^\wedge$ by \cite[Lemma~2.7]{ASS}.
\end{proof}

%%%%%%%%%%%%%%%%%%%%%%%%%%%%%%%%%%%%%%%%%%%%%%%%%%%%%%%%%%%%%%%%%%%%%%%%%%%%%%%%%%%%%%%%%%%%%%%%%%%%%%%%%%%%%%%%%%%%%%%%%%%%%%

\begin{example}\label{ex-fell}
This example illustrates Corollary~\ref{cor-1} and also shows that condition (3) in Theorem~\ref{thm-fell} is not implied by conditions (1) and (2).  Let
\[A=\left\{x=(x_n)_{n\geq 0}: x_n\in M_2(\C), x_n\to\begin{pmatrix} \lambda(x)&0\\0&\mu(x) \end{pmatrix}\right\}
\]
and $G=\{1,-1\}$.  We define $\alpha : G \to \Aut(A)$ by specifying $\alpha_{-1}(x)$
for $x\in A$:
\[
(\alpha_{-1}(x))_{2m}=u^*x_{2m+1}u\quad\text{and}\quad (\alpha_{-1}(x))_{2m+1}=u^*x_{2m}u\quad\text{where}\quad u=\begin{pmatrix} 0&1\\1&0 \end{pmatrix}.
\]
So $\alpha_{-1}(x)$ is the sequence $u^*x_1u$, $u^*x_0u$, $u^*x_3u$, $u^*x_2u,\dots$ and converges to $\begin{pmatrix} \mu(x)&0\\0&\lambda(x)
\end{pmatrix}$.
Here  $\hat A\simeq\{\pi_m:m\in\NN\}\cup\{\lambda,\mu\}$    where $\{\pi_m:m\in\NN\}$ is discrete in the relative topology and $\pi_m\to\lambda,\mu$ as $m\to\infty$. Note
\[
-1\cdot\lambda=\mu,\quad- 1\cdot\mu=\lambda\quad\text{and}\quad -1\cdot\pi_{2m}=\pi_{2m+1},
\]
so the induced action of $G$ on   $\hat A\simeq \Prim  A$    is free (and $(G,\hat A)$ is trivially Cartan since $G$ is compact).  Note that $\lambda$ and $\mu$ do not have Hausdorff $G$-invariant neighbourhoods.

Each orbit is closed and homeomorphic to $G$, so \[ \hat A/G\simeq(A\rtimes_\alpha G)^\wedge=\{\Ind\pi_{2m+1}:m\in\NN\}\cup\{\Ind\lambda\}.\]
Let $\{e_1,e_2\}$ be the standard basis in $\C^2$ and define $\eta_1,\eta_2:G\to \C^2$ by
\[\eta_1(1)=e_1,\ \eta_1(-1)=0\quad\text{and}\quad\eta_2(1)=0,\ \eta_2(-1)=e_2.\] Then $\langle\eta_1\,,\,\eta_2\rangle =0$. Also let $\xi:G\to \C$ be the function $\xi(1)=1$ and $\xi(-1)=0$.
It is straightforward to verify that
\[\lim_{m\to\infty}\langle\Ind\pi_{2m+1}(\cdot)\eta_i\,,\,\eta_i\rangle=\langle\Ind\lambda(\cdot)\xi\,,\,\xi\rangle\quad(i=1,2)\]
and hence $M_U(\Ind\lambda)\geq 2$ (as also predicted by Corollary~\ref{cor-1}).  Thus $A\rtimes_\alpha G$ is not a Fell algebra. Note that $A\rtimes_\alpha G$  has bounded trace by Theorem~\ref{thm-bt1} below.
% \footnote{We believe that $A\rtimes_\alpha G$ is isomorphic to
% \[
% \{
% (x_n)_{n\geq 0}:x_n\in M_4(\C),
% x_n\to\textrm{diag}
% (\sigma(x),\sigma(x))
% \text{ for some $\sigma(x)\in M_2(\C)$}
% \}
% \]
% or a variant where one of the summands is $u\sigma(x)u^*$ or something similar.}
\end{example}
%%%%%%%%%%%%%%%%%%%%%%%%%%%%%%%%%%%%%%%%%%%%%%%%%%%%%%%%%%%%%%%%%%%%%%%%%%%%%%%%%%%%%%%%%%%%%%%%%%%%%%%%%%%%%%%%%%%%%%%%%%%%%%
Our next theorem uses and  improves \cite[Corollary~3.9]{AD} by
incorporating into the if-and-only-if statement the standing
assumptions that $A$ has bounded trace and orbits in $\hat A$ are
locally closed.

\begin{thm}\label{thm-bt1}
Suppose that $(A,G,\alpha)$ is a separable $C^*$-dynamical system such that  the induced action of $G$ on   $\Prim A$    is free. Then $A\rtimes_\alpha G$ has bounded trace if and only if
\begin{enumerate}
\item $A$ has bounded trace; and
\item the action of $G$ on $\hat A$ is integrable; and
\item the orbits in $\hat A$ are locally closed.
\end{enumerate}
\end{thm}
%%%%%%%%%%%%%%%%%%%%%%%%%%%%%%%%%%%%%%%%%%%%%%%%%%%%%%%%%%%%%%%%%%%%%%%%%%%%%%%%%%%%%%%%%%%%%%%%%%%%%%%%%%%%%%%%%%%%%%%%%%%%%%
\begin{proof}
Suppose that $A\rtimes_\alpha G$ has bounded trace.  Then $A$ has bounded trace and $G$ acts integrably on $\hat A$ by Corollary~\ref{cor-A}. Since the action of $G$ on $\hat A$ is free and both $A$ and $A\rtimes_\alpha G$ are Type 1, the orbits are locally closed by Theorem~\ref{thm-tak3}.

The converse is \cite[Theorem~3.6]{AD}.
\end{proof}

%%%%%%%%%%%%%%%%%%%%%%%%%%%%%%%%%%%%%%%%%%%%%%%%%%%%%%%%%%%%%%%%%%%%%%%%%%%%%%%%%%%%%%%%%%%%%%%%%%%%%%%%%%%%%%%%%%%
In view of Theorem~\ref{thm-fell} and Example~\ref{ex-fell} it is natural to ask whether there is an example  for which conditions (1) and (2) of Theorem~\ref{thm-bt1} hold but (3) fails or whether (1) and (2) together imply (3).  We have been unable to answer this question, but the following theorem gives related information in terms of Hausdorff orbits.

%%%%%%%%%%%%%%%%%%%%%%%%%%%%%%%%%%%%%%%%%%%%%%%%%%%%%%%%%%%%%%%%%%%%%%%%%%%%%%%%%%%%%%%%%%%%%%%%%%%%%%%%%%%%%%%%%%%%%%%%%%%%%
\begin{thm}\label{thm-bt2}
Suppose that $(A,G,\alpha)$ is a separable $C^*$-dynamical system such that   $A$ has bounded trace, $G$ acts integrably on $\hat A$, and  the induced action of $G$ on  $\hat A\simeq \Prim A$  is free.     Then the following conditions are equivalent:
\begin{enumerate}
\item each orbit in $\hat A$ is Hausdorff (in the relative topology);
\item each orbit in $\hat A$ is locally closed;
\item $A\rtimes_\alpha G$ has bounded trace;
\item $A\rtimes_\alpha G$ is Type I;
\item for each $\pi\in\hat A$, the map $s\mapsto s\cdot \pi$ of $G$ onto $G\cdot\pi$ is a homeomorphism.
\end{enumerate}
\end{thm}
%%%%%%%%%%%%%%%%%%%%%%%%%%%%%%%%%%%%%%%%%%%%%%%%%%%%%%%%%%%%%%%%%%%%%%%%%%%%%%%%%%%%%%%%%%%%%%%%%%%%%%%%%%%%%%%%%%%%%%%%%%%%%
\begin{proof}
That (2) implies (3) is \cite[Theorem~3.6]{AD} (or see
Theorem~\ref{thm-bt1}). Every bounded-trace $C^*$-algebra is
liminal, hence Type I, so (3) implies (4). Items (4), (5) and (2) are equivalent by
Theorem~\ref{thm-tak3}.
%That (4) implies (5) follows from  Theorem~\ref{thm-tak2} since the orbits are
%locally closed in $\hat A$ if and only if the orbits  are canonically homeomorphic
%to $G$ by \cite[Theorem~1]{Gli}.
That (5) implies (1) is immediate since $G$ is Hausdorff.
It remains to prove that (1) implies (2).

Suppose that each orbit is Hausdorff but that for some $\pi\in\hat A$, $G\cdot\pi$ is not locally closed.   We consider two cases, and show that each leads to a contradiction.

First, suppose that for every open neighbourhood $V$ of $\pi$,
$\phi_\pi^{-1}(V)$ is not relatively compact. Then, even though
$\hat A$ might not be Hausdorff,  it follows as in the proof of the
(1) $\Longrightarrow$ (2) direction of \cite[Lemma~2.1]{AaH} that
for each $k\in\P$ the sequence $\pi$,$\pi$,$\pi,\dots$ converges
$k$-times to $\pi$ in $\hat A/G$. Since the action of $G$ on $\hat
A$ is jointly continuous \cite[Lemma 7.1]{raewill}, we may proceed
as in the (2) $\Longrightarrow$ (3) direction of
\cite[Lemma~2.1]{AaH} (even though $\hat A$ might not be Hausdorff)
to show that for every open neighbourhood $V$ of $\pi$,
$\nu(\phi_\pi^{-1}(V))=\infty$.
%\footnote{The transformation group $(G,X)$ of \cite[Lemma~2.1]{AaH} is assumed to be
%locally compact Hausdorff. In the (1) $\Longrightarrow$ (2) the assumption that
%$\phi_\pi^{-1}(V)$ is not relatively compact for any open neighbourhood $V$ of $\pi$
%replaces the  Hausdorff assumption on $X$ there, and since $\hat A$ is quasi
%locally-compact we do not need  Hausdorffness of $X$ in the (2) $\Longrightarrow$ (3)
%direction.}
By \cite[Lemma~3.5]{AD} this contradicts the integrability of the action of $G$ on $\hat A$.

Second, suppose that there exists an open neighbourhood $V$ of $\pi$ such that $\phi_\pi^{-1}(V)$ is relatively compact.   Since $G\cdot\pi$ is not locally closed there is a sequence $(\sigma_n)_n$ in $V\cap(\overline{G\cdot\pi}\setminus G\cdot\pi)$ which is convergent to $\pi$. Since $K=\overline{\phi_\pi^{-1}(V)}$ is compact,
for each $n\geq 1$  there is a convergent sequence $(t_j^{(n)})_j$ in $K$ (with limit $s_n$ in $K$ say) such that $t_j^{(n)}\cdot\pi\to_j\sigma_n$.  Note also that $t_j^{(n)}\cdot\pi\to_j s_n\cdot\pi$.

Replacing $(\sigma_n)$ by a subsequence, we may assume without loss of generality that there exists $s\in K$ such that $s_n\to_n s$.  We claim that $s\cdot\pi= \pi$. To see the claim, suppose that $s\cdot\pi\neq\pi$. Since $G\cdot\pi$ is Hausdorff there exist open neighbourhoods $W_1$ and $W_2$ in $\hat A$ of $s\cdot\pi$ and $\pi$ respectively, such that
\begin{equation}\label{eq-thm4} W_1\cap W_2\cap G\cdot\pi=\emptyset.\end{equation}
Since $s_n\cdot\pi\to s\cdot\pi$ and $\sigma_n\to\pi$ there exists $n_0$ such that $s_{n_0}\cdot\pi\in W_1$ and $\sigma_{n_0}\in W_2$.  Since $t_j^{(n_0)}\cdot\pi\to_j\sigma_{n_0}$ and $t_j^{n_0}\cdot\pi\to_js_{n_0}\cdot\pi$, there exists $j_0$ such that $t_{j_0}^{(n_0)}\cdot\pi\in W_1\cap W_2\cap G\cdot\pi$, contradicting \eqref{eq-thm4}, and hence $s\cdot\pi=\pi$ as claimed. It follows that $s_n\cdot \pi\to\pi$.

Let $(N_k)_k$ be a decreasing base of neighbourhoods of $0$ in $A^*$ and let $\phi$ be a pure state of $A$ associated with $\pi$.
Temporarily fix $k\geq 1$.  The canonical image of $(\phi+\frac{1}{2}N_k)\cap P(A)$ in $\hat A$ is an open neighbourhood of $\pi$, and since  $\sigma_n\to\pi$ and $s_n\cdot\pi\to\pi$ this neighbourhood contains $\sigma_{n_k}$ and $s_{n_k}\cdot\pi$ for some $n_k$.
So there are unit vectors $u_k$ and $v_k$, in the Hilbert spaces of $\sigma_{n_k}$ and $\pi$, respectively, such that
\[
\psi_k^{(1)}:=\langle\sigma_{n_k}(\cdot)u_k\,,\, u_k\rangle\in\phi +\frac{1}{2}N_k\quad\text{and}\quad
\psi_k^{(2)}:=\langle (s_{n_k}\cdot\pi)(\cdot)v_k\,,\, v_k\rangle\in\phi +\frac{1}{2}N_k.
\]

The sequence $(t^{(n_k)}_j\cdot\pi)_j$  converges to both of the inequivalent representations $\sigma_{n_k}$ and $s_{n_k}\cdot\pi$.
We apply Lemma~\ref{lem-pems2} to the sequence $(t^{(n_k)}_j\cdot\pi)_j$, the states $\psi_k^{(1)}$ and $\psi_k^{(2)}$, the neighbourhoods $\frac{1}{2}N_k$ and  $\frac{1}{2}N_k$, and $\epsilon=1/k$, to obtain $r_k:=t_{j_k}^{(n_k)}\in G$ and unit vectors $\xi_k^{(1)}$ and $\xi_k^{(2)}$ in $\H_\pi$ such that
\begin{gather*}
\langle (r_k\cdot\pi)(\cdot)\xi^{(i)}_k\,,\,\xi^{(i)}_k\rangle\in\psi_k^{(i)}+\frac{1}{2}N_k\subset\phi+N_k\quad\quad(i=1,2),\label{eq-thm4-1}\\
\text{and\ }|\langle \xi^{(1)}_k\,,\, \xi^{(2)}_k\rangle|<\frac{1}{k}\label{eq-thm4-3}.
\end{gather*}

Now let $k$ vary: since $(N_k)_k$ is a decreasing sequence of neighbourhood of $0$,
\[
\langle \xi_k^{(1)}\,,\, \xi_k^{(2)}\rangle\to 0\quad\text{and}\quad\langle (r_k\cdot\pi)(\cdot)\xi_k^{(i)}\,,\,\xi_k^{(i)}\rangle\to\phi\quad(i=1,2)
\]
as $k\to\infty$.
By Lemma~\ref{lem-gs} and \cite[Lemma~5.2(i)]{ASS}, $M_U(\pi,(r_k\cdot\pi))\geq 2$.

Now recall that $A$ is assumed to have bounded trace and that $M_U(r_k\cdot\pi)=M_U(\pi)$.  We now have
\[\infty>M_U(\pi)\geq M_U(\pi,(r_k\cdot\pi))M_U(r_k\cdot\pi)\geq 2M_U(\pi)\]
by \cite[Theorem~1.5]{AKLSS}, a contradiction.

Both cases have led to contradictions, so every orbit in $\hat A$  must be locally closed.
\end{proof}
%%%%%%%%%%%%%%%%%%%%%%%%%%%%%%%%%%%%%%%%%%%%%%%%%%%%%%%%%%%%%%%%%%%%%%%%%%%%%%%%%%%%%%%%%%%%%%%%%%%%%%%%%%%%%%%%%%%%%%%%%%%%%

\section{An example}\label{sec-4}
In this section we build an example that further illustrates the power of Theorems~\ref{thm-a} and \ref{3.5'}. In more detail, we consider a $C^*$-algebra $A=C_0(X)\otimes B$ with a diagonal action $\alpha=\gamma\otimes\beta$ of $G=\R$.  While $B\rtimes_\beta G$ does not have bounded trace, the action of $G$ on $X$ is sufficiently well-behaved (it is integrable!) to ensure that $A\rtimes_\alpha G$ does have bounded trace. We have chosen $X$ and $B$ in such a way that there is a sequence $\epsilon_{x_n}\otimes\epsilon_{y_n}$ in $\hat A=X\otimes\hat B$ which converges $2$-times to some $\epsilon_{x_0}\otimes\lambda\in\hat A$, and  such that \[M(\epsilon_{x_0}\otimes\lambda,(\epsilon_{x_n}\otimes\epsilon_{y_n}))=3\quad\text{and}\quad M(\epsilon_{x_0}\otimes\lambda,(\epsilon_{s_n\cdot x_n }\otimes \epsilon_{s_n\cdot y_n}))=5,\]
where $(s_n)\subset G$ implements the $2$-times convergence of $\epsilon_{x_n}\otimes\epsilon_{y_n}$ to $\epsilon_{x_0}\otimes\lambda$.  (It will be clear from the construction of B that 3 and 5 could be
replaced by any two positive integers.)
We show below that Theorems~\ref{thm-a} and \ref{3.5'} imply that $M_L(\Ind(\epsilon_{x_0}\otimes\lambda),(\Ind(\epsilon_{x_n}\otimes\epsilon_{y_n}))=8,9$ or $10$, and  further investigation shows  $M_L(\Ind(\epsilon_{x_0}\otimes\lambda),(\Ind(\epsilon_{x_n}\otimes\epsilon_{y_n}))=8$.

 Fix an orthonormal basis $\{e_i\}_{i\in\mathbb{N}}$ for $l^2(\Z)$, and let $B$ be the $C^*$-subalgebra of $C_b(\R,\K(l^2(\Z))$ consisting of $f$ such that, for some $\lambda(f)\in\R$,
\begin{gather*}
\lim_{t\to-\infty} f(t)=\lambda(f)\sum_{i=1}^3 e_i\otimes\overline{e_i}=\lambda(f)\diag(1,1,1,0,\dots)\\
\lim_{t\to\infty}f(t)=\lambda(f)\sum_{i=1}^5 e_i\otimes\overline{e_i}=\lambda(f)\diag(1,1,1,1,1,0,0,\dots).
\end{gather*}
 Let $G=\R$ and, for $s,t\in G$ and $f\in B$ set
\[(\beta_s(f))(t)=f(t-s).\]

\begin{lemma}
Let $(B,G,\beta)$ be as above. Then $\beta$ is a strongly continuous action of $G=\R$ on $B$.
\end{lemma}
\begin{proof}
It is easy to see that  $\beta_s(f)\in B$ and  $\lambda(\beta_s(f))=\lambda(f)$ for all $s\in G$ and $f\in B$, and  that $\beta:G\to\Aut B$ is a homomorphism.

To see that $\beta$ is strongly continuous, fix $f\in B$ and $\epsilon>0$.  There exists $T_0>0$ such that
\begin{gather*}
\|f(r)-\lambda(f)\diag(1,1,1,0,\dots)\|<\epsilon/2\text{\ for all $r<-T_0$\  and }\\
\|f(r)-\lambda(f)\diag(1,1,1,1,1,0,0,\dots)\|<\epsilon/2\text{\ for all $r>T_0$.}
\end{gather*}
So if $r,r'>T_0$ or $r,r'<-T_0$ then $\|f(r)-f(r')\|<\epsilon$.
Moreover, since $f$ is uniformly continuous on $[-T_0-1,T_0+1]$ there exists $\delta\in (0,1/2)$ such that $\|f(r)-f(r')\|<\epsilon$ whenever $r,r'\in [-t_0-1,t_0+1]$ with $|r-r'|<\delta$.
Now let  $s,s'\in G$ such that $\|s-s'\|<\delta$ and $t\in \R$. If $t-s, t-s'\in [-T_0-1,T_0+1]$ then $\|f(t-s)-f(t-s')\|<\epsilon$.  If one of $t-s, t-s'$ is not in $[-T_0-1,T_0+1]$ then either $t-s, t-s'>T_0$ or $t-s, t-s'<-T_0$, and again $\|f(t-s)-f(t-s')\|<\epsilon$. Thus
\[
\|\beta_s(f)-\beta_{s'}(f)\|=\sup_{t\in\R}\|f(t-s)-f(t-s')\|<\epsilon
\]
whenever $|s-s'|<\delta$, and hence $\beta$ is a strongly continuous action of $G$ on $B$.
\end{proof}

We have  $\hat B\simeq \R\cup\{\lambda\}$, the one-point compactification of $\R$. There are two orbits under the  action of $G$ on $\R\cup\{\lambda\}$ induced from the homeomorphism, $\R$ and $\{\lambda\}$; the first orbit is not closed but is locally closed and the second orbit is closed.
To see that $\lambda$ has no integrable neighbourhood let $U$ be any neighbourhood of $\lambda$ in $\R\cup\{\lambda\}$.  Then there exists $t_0>0$ such that $(-\infty, -t_0)\cup(t_0,\infty)\subset U$.  Let $t>t_0$.  Then $s\cdot t=t+s\in U$ for all $s\geq 0$.   So $\nu(\{s\in G:s\cdot t\in U\})=\infty$.

Fix $k\in\P$. If $(t_n)_n$ is a sequence in $\R$  such that $t_n\to\lambda$ (that is, $|t_n|\to\infty$) then $(t_n)$ is $k$-times convergent in $(\R\cup\{\lambda\})/G$ to $\lambda$.  To see this it suffices to note that, for example,
\[
j|t_n|\cdot t_n=t_n+j|t_n|\to\lambda\quad\quad(j=2,\dots,k)
\]
and that $(j-i)|t_n|\to\infty$ for $1\leq i<j\leq k$.

Note that
\[
M_L(\lambda,(\epsilon_{j|t_n|\cdot t_n})_n)=5=M_U(\lambda,(\epsilon_{j|t_n|\cdot t_n})_n)\quad\quad (2\leq j\leq k).
\]
The numbers $M_L(\lambda, (\epsilon_{t_n}))$ and $M_U(\lambda, (\epsilon_{t_n}))$ depend on how $t_n\to\lambda$: for example, if $t_n\to-\infty$ then $M_L(\lambda, (\epsilon_{t_n}))=3=M_U(\lambda, (\epsilon_{t_n}))$; if $t_n=(-1)^nn$ then $M_L(\lambda, (\epsilon_{t_n}))=3$ and $M_U(\lambda, (\epsilon_{t_n}))=5$.
In particular, $B$ has bounded trace because $M_U(\sigma)\leq 5$ for all $\sigma\in\hat B$. The action $\beta$ is not free on $\hat B$ since $\lambda$ is a fixed point, so \cite[Corollary~3.9]{AD}  does not apply;  it turns out that $B\rtimes_\beta G$ does not have bounded trace.

%%%%%%%%%%%%%%%%%%%%%%%%%%%%%%%%%%%%%%%%%%%%%%%%

\begin{lemma}
\begin{enumerate}
\item $B\rtimes_\beta G$ is an extension of $C_0(\R)$ by the compact operators $\K(L^2(G))$.
\item $B\rtimes_\beta G$  does not have bounded trace.
\end{enumerate}
\end{lemma}
\begin{proof}
(1) Consider the ideal  \[J=C_0(\R,\K(l^2(\Z)))=\{f\in B:f(t)\to 0\text{\ as \ }|t|\to\infty\}.\]  Then $J$ is $G$-invariant and $G$ acts freely on $\hat J\simeq\R$. By the $G$-invariance we obtain the short-exact sequence
\[
0\to J\rtimes_\beta G\to B\rtimes_\beta G\to (B/J)\rtimes_\beta G\to 0.
\]
Since the action of $G$ on $\hat J$ is free with the unique orbit  closed and $J$ is Type I, we have  $(J\rtimes_\beta G)^\wedge\simeq\hat J/G$, a singleton. Since $J\rtimes_\beta G$ is separable, it is isomorphic to $\K(L^2(\R))$.
%\footnote{Alternatively,
%\begin{align*}
%J\rtimes_\beta G
%&\cong (C_0(\R)\otimes \K(l^2(\Z)))\times_{\lt\otimes\id}G
%%\\&
%\cong (C_0(\R)\times_{\lt} \R)\otimes \K(l^2(\Z))
%%\\&
%\cong \K(L^2(\R))\otimes \K(l^2(\Z))\cong \K(L^2(\R)).
%\end{align*}}
Moreover, $(B/J)\rtimes_\beta G\cong \C\otimes C^*(\R)\cong C_0(\R)$.

(2) Let $\pi\in (B\rtimes_\beta G)^\wedge$ such that $\{\pi\}=(J\rtimes_\beta G)^\wedge$. We will show that $\pi$ is faithful.  This suffices because then  $\pi(B\rtimes_\beta G)\supsetneq \K(\H_\pi)$ and $M_U(\sigma)=\infty$ for all $\sigma\in\overline{\{\pi\}}\setminus\{\pi\}$ by \cite[Proposition~4.8]{A}, and hence $B\rtimes_\alpha G$ does not have bounded trace by \cite[Theorem~2.6]{ASS}.

Suppose that $\pi$ is not faithful.   Since the action of $G$ on $\hat B$ is almost free and $G=\R$ is amenable, by \cite[Theorem~9]{LN} there exists a non-zero $G$-invariant ideal $I$ of $B$ such that $\ker\pi\supset I\rtimes_\beta G$. Since $\ker\pi\cap(J\rtimes_\beta G)=\{0\}$ we have $I\cap J=\{0\}$, which contradicts that $J$ is an essential ideal of $B$.  Hence $\pi$ is faithful as required.
\end{proof}
%%%%%%%%%%%%%%%%%%%%%%%%%%%%%%%%%%%%%%%%%%%%
Now let $(G, C_0(X), \gamma)$ be the $C^*$-dynamical system arising from Green's example in \cite{green1} and set
\[
A=C_0(X)\otimes B\quad\text{and}\quad\alpha=\gamma\otimes\beta.
\]
Note that   $\alpha$ induces a free and integrable action of $G=\R$ on    \[\Prim A\simeq\hat A\simeq X\times\hat B\simeq X\times(\R\cup\{\lambda\})\]    because the action of $G$ on $X$ is free and integrable.   Moreover,    the orbits in $\hat A$ are closed by \cite[Proposition~1.17]{rieffel}  because $\hat A$ is Hausdorff and the action of $G$ on $\hat A$ is integrable.  Thus $A\rtimes_\alpha G$ has bounded trace by \cite[Theorem~3.6]{AD}.

In Green's example  $x_n=(2^{-2n},0,0)\to x_0=(0,0,0)$ and also $(2n+\pi)\cdot x_n\to x_0$.
Let $y_n=-n\in\hat B$, and consider $(x_n,y_n)\in X\times(\R\cup\{\lambda\})$. Then
\[
(x_n,y_n)\to(x_0,\lambda)\ \text{\ and \ }
(2n+\pi)\cdot(x_n,y_n)=((2n+\pi)\cdot x_n,n+\pi)\to(x_0,\lambda),
\]
so $(x_n,y_n)$ converges $2$-times in $X\times(\R\cup\{\lambda\})/G$  to $(x_0,\lambda)$.
Note that
\begin{gather*}
M_L(\epsilon_{x_0},(\epsilon_{x_n}))=1=M_U(\epsilon_{x_0},(\epsilon_{x_n}))\\
M_L(\epsilon_{x_0},(\epsilon_{(2n+\pi)\cdot x_n }))=1=M_U(\epsilon_{x_0},(\epsilon_{(2n+\pi)\cdot x_n }))\\
M_L(\lambda,(\epsilon_{y_n}))=3=M_U(\lambda,(\epsilon_{y_n}))\\
M_L(\lambda,(\epsilon_{(2n+\pi)\cdot y_n }))=5=M_U(\lambda,(\epsilon_{(2n+\pi)\cdot y_n })),
\end{gather*}
so by Corollary~\ref{cor-mult-tensor} we have
\[
M(\epsilon_{x_0}\otimes\lambda,(\epsilon_{x_n}\otimes\epsilon_{y_n}))=3\quad\text{and}\quad M(\epsilon_{x_0}\otimes\lambda,(\epsilon_{(2n+\pi)\cdot x_n }\otimes \epsilon_{(2n+\pi)\cdot y_n}))=5.
\]
So by Theorem~\ref{thm-a},
\begin{equation}\label{eq-ex-lowerbound}
M_L(\Ind(\epsilon_{x_0}\otimes\lambda),(\Ind(\epsilon_{x_n}\otimes\epsilon_{y_n})))\geq 3+5=8.
\end{equation}

We will use Theorem~\ref{3.5'} to obtain an  upper bound for $M_L(\Ind(\epsilon_{x_0}\otimes\lambda),(\Ind(\epsilon_{x_n}\otimes\epsilon_{y_n})))$. Note that $M_U(\epsilon_{x_0}\otimes\lambda)=M_U(\epsilon_{x_0})M_U(\lambda)=5$.
Let $V$ be any open neighbourhood of $(x_0,\lambda)$ in $X\times(\R\cup\{\lambda\})$. There exist open neighbourhoods $U$ of $x_0$ in $X$ and $W$ of $\lambda$ in $\R\cup\{\lambda\}$ such that $U\times W\subset V$.  Moreover, we can assume that \[U=X\cap((-\delta,\delta)\times(-\delta,\delta)\times(-\delta,\delta))\] for some $0<\delta<1/4$ and that $W\supset (-\infty, -T)\cup(T,\infty)$ for some $T>0$.

Since $\delta<1/4$,  $U$ never meets the arc joining the two line segments of an orbit $G\cdot x_n$ in $X$.  Choose $n_0$ such that $n> n_0$ implies $2^{-2n}<\delta$. Then
\[
U\cap G\cdot x_n=\{s\cdot x_n:s\in (-\delta,\delta)\cup(-\delta+2n+\pi,\delta+2n+\pi)\}
\]
whenever $n\geq n_0$.
If $s\in(-\delta,\delta)$ then $s\cdot y_n=s-n\in(-\delta-n,\delta-n)$, and we note that
\begin{gather*}
(-\delta-n,\delta-n)\cap(-\infty, -T)=(-\delta-n,\delta-n)\quad\text{provided $n>T+\delta$};\\
(-\delta-n,\delta-n)\cap(T,\infty)=\emptyset\quad\text{for all $n> 0$}.
\end{gather*}
On the other hand, if $s\in(-\delta+2n+\pi,\delta+2n+\pi)$ then $s\cdot y_n=s-n\in(-\delta+n+\pi,\delta+n+\pi)$, and we note that
\begin{gather*}
(-\delta+n+\pi,\delta+n+\pi)\cap(T,\infty)=(-\delta+n+\pi,\delta+n+\pi)\quad\text{provided $n>T+\delta-\pi$};\\
(-\delta+n+\pi,\delta+n+\pi)\cap(-\infty,-T)=\emptyset\quad\text{for all $n> 0$}.
\end{gather*}
Let $n_1>T+\delta$.  Thus  if  $n\geq\max\{n_0,n_1\}$ we have $\phi_{x_n}^{-1}(U)\subset \phi_{y_n}^{-1}(W)$. Since $\lambda$ is a fixed point,
\begin{align*}
\nu(\{s\in G:s\cdot(x_n,y_n)\in U\times W\}
&=
\nu(\{s\in G:s\cdot x_n\in U\text{\ and\ }s\cdot y_n\in  W\}\\
%&=\nu(\{s\in((-\delta,\delta)\cup(-\delta+2n+\pi,\delta+2n+\pi)):s\cdot y_n\in W\}\\
&=\nu((-\delta,\delta)\cup(-\delta+2n+\pi,\delta+2n+\pi))\\
&=2\nu\{s\in G:s\cdot x_0\in U\}\\
&=2\nu\{s\in G:s\cdot x_0\in U\text{\ and\ } s\cdot\lambda=\lambda\}\\
&=2\nu\{s\in G:s\cdot(x_0,\lambda)\in U\times W\},
\end{align*}
whenever $n\geq\max\{n_0,n_1\}$.
So by Corollary~\ref{cor-Mupper},
\[
M_L(\Ind(\epsilon_{x_0}\otimes\lambda),(\Ind(\epsilon_{x_n}\otimes\epsilon_{y_n})))\leq M_U(\Ind(\epsilon_{x_0}\otimes\lambda),(\Ind(\epsilon_{x_n}\otimes\epsilon_{y_n})))\leq \lfloor 2M_U((x_0,\lambda))\rfloor=10.
\]
Combining with \eqref{eq-ex-lowerbound} we obtain $M_L(\Ind(\epsilon_{x_0}\otimes\lambda),(\Ind(\epsilon_{x_n}\otimes\epsilon_{y_n})))=8,9$ or $10$.

\begin{remark}
If the statements of Theorem~\ref{3.1'} and Theorem~\ref{3.5'} were to hold with $u=M_U(\pi, (\pi_n))$ instead of $u=M_U(\pi)$, then in the example above we would get $8\leq M_L(\Ind(\epsilon_{x_0}\otimes\lambda),(\Ind(\epsilon_{x_n}\otimes\epsilon_{y_n})))\leq 6=2\cdot 3$ because $M(\epsilon_{x_0}\otimes\lambda,(\epsilon_{x_n}\otimes\epsilon_{y_n}))=3$.  The point is that $\tr(\pi_n(a))$ and $\tr(s\cdot\pi(a))$ can vary widely as $s$ does.

This tells us that to sharpen the estimates of Theorems~\ref{3.1'} and~\ref{3.5'} we need to know specifics about the action of $G$ on $\hat A$, as we do in this example.
\end{remark}

%%%%%%%%%%%%%%%%%%%%%%%%%%%%%%%%%%%%%%%%%%%%%%%%%
\begin{lemma}\label{lem-ex-upper-bound}
With notation as above, \[M_L(\Ind(\epsilon_{x_0}\otimes\lambda),(\Ind(\epsilon_{x_n}\otimes\epsilon_{y_n})))\leq M_U(\Ind(\epsilon_{x_0}\otimes\lambda),(\Ind(\epsilon_{x_n}\otimes\epsilon_{y_n})))\leq 8.\]
\end{lemma}
\begin{proof}
To save some typing we set
\[\pi:=\epsilon_{x_0}\otimes\lambda\quad\text{and}\quad\pi_n:=\epsilon_{x_n}\otimes\epsilon_{y_n}.\]
We start by showing that $M_L(\Ind\pi,(\Ind\pi_n))\leq 8$ using  ideas from the proof of Theorems~\ref{3.1'}; the point is that we can be very specific about the functions $F$ and $b$ and the element $a$ (see below) and that we will know bounds for $\tr(s\cdot\pi_n(a))$ for all $s$ of interest.

Fix $\delta$ such that $0<\delta<1/16$. Choose $F\in C_c(\hat A)^+$ such that
\begin{enumerate}
\item[i)] $\supp F$ is contained in a relatively compact neighbourhood $V=U\times W$ of $\pi$, where $U=X\cap\big((-1/2,1)\times(-1/2-\delta,1/2+\delta)\times(-1/4,1/4)\big)$;
\item[ii)] $0\leq F\leq 1$;
\item[iii)] $F(s\cdot\pi)=F(\pi\circ\lt_{s^{-1}})=F(\epsilon_{s\cdot x_0}\otimes\lambda)=1$ if $s\in(-1/2,1/2)$.
\end{enumerate}
Then choose $b\in C_c(G\times\hat A)$ such that
\begin{enumerate}
\item[i)] $0\leq b\leq 1$;
\item[ii)] $b$ is  one on $(\supp F_\pi)(\supp F_\pi)^{-1}\times \supp F$; and
\item[iii)] $b$ is zero off $(-3,3)\times M$, where $M$ is some open neighbourhood of $\supp F$.
\end{enumerate}
Finally, let $a=a_1\otimes a_2\in A^+$ where
\begin{enumerate}
\item[i)] $a_1\in C_c(X)^+$ with $\supp a_1\subset X\cap\big([-1/2,1]\times[-1,1]\times\{0\}\big)$;
\item[ii)] $0\leq a_1\leq 1$ and $a_1(s\cdot x_0)=1$ for $s\in (-1/2,1/2)$;
\item[iii)] $a_2\in B$ is defined by
\[a_2(t)=\begin{cases}
\diag(1,1,1,0,0,\dots), &\text{if $t\leq 0$};\\
\diag(1,1,1,t^2,t^2,0,0,\dots), &\text{if $0<t<1$};\\
\diag(1,1,1,1,1,0,0\dots), &\text{if $t\geq 1$}.
\end{cases}
\]
Note that $a_2$ is positive, $\|a_2\|=1$, and  $\lambda(a_2)=1$.
\end{enumerate}

There are three consequences of our choices above that are crucial for the argument that follows.  First, for $n\geq 1$,
\begin{equation*}\label{eq-cons1}
s\cdot\pi_n(a)=\epsilon_{x_n}\otimes\epsilon_{y_n}(\lt_{s^{-1}}(a_1)\otimes\beta_{s^{-1}}(a_2))=a_1(s\cdot x_n)a_2(-n+s)
\end{equation*}
is zero if $s\notin [-1,1]\cup [-1+2n+\pi,1+2n+\pi]$ because $\supp a_1\subset X\cap\big([-1/2,1]\times[-1,1]\times\{0\}\big)$.
Second, for $n\geq 1$,
\begin{equation*}
\tr(s\cdot\pi_n(a))=a_1(s\cdot x_n)\tr(a_2(-n+s))=\begin{cases}
3 &\text{if $s\in[-1,1]$};\\
5 &\text{if $s\in[-1+2n+\pi,1+2n+\pi]$}.
\end{cases}
\end{equation*}
Third, for any $n\geq 1$, if  $s\in [-1+2n+\pi,1+2n+\pi]$ and $v\in [-1,1]$ then $s-v\in [-2+2n+\pi, 2+2n+\pi]$ for $n\geq 1$; in particular $s-v, v-s\notin (-3,3)$ and hence $b(s-v,\sigma)=b(v-s,\sigma)=0$ for any $\sigma\in\hat A$.  Similarly, if  $s\in [-1,1]$ and $v\in [-1+2n+\pi,1+2n+\pi]$ implies $b(s-v,\sigma)=b(v-s,\sigma)=0$ for any $\sigma\in\hat A$.

Using the three consequences just discussed above we have, for $n\geq 1$,
\begin{align*}
\tr\big(\big( b(v-s,s\cdot\pi_n) &s\cdot\pi_n(a)+b(s-v,v\cdot\pi_n)v\cdot\pi_n(a)\big)^2\big)\notag\\
&\leq \|b(v-s,s\cdot\pi_n) s\cdot\pi_n(a)+b(sv^{-1},v\cdot\pi_n)v\cdot\pi_n(a)  \|\notag\\
&\hskip2cm\cdot\tr\big(b(v-s,s\cdot\pi_n) s\cdot\pi_n(a)+b(s-s,v\cdot\pi_n)v\cdot\pi_n(a)\big)\notag\\
&\leq\begin{cases}
2(3+3), &\text{if $s,v\in[-1,1]$};\\
2(5+5), &\text{if $s,v\in[-1+2n+\pi,1+2n+\pi]$};\\
0, &\text{else}.
\end{cases}
\end{align*}

For $t\in G$ and $\sigma\in\hat A$ set
$B(t,\sigma)=F(\sigma)F(t^{-1}\cdot \sigma)b(t^{-1},\sigma)\Delta(t)^{-1/2}$,
$C(t)=\Psi(B(t,\cdot))a$, and
$E=\frac{1}{2}(C+C^*)$. Then a formula for
$\tr(\Ind\pi_n(E*E))$ follows from Theorem~\ref{thm-duflo} as in Theorem~\ref{3.1'}.  Using \eqref{trace-formula} we have, for $n\geq 1$,
\begin{align*}
\tr(\Ind\pi_n(E*E))
&=\frac{1}{4}\int_G\int_G
F(s\cdot\pi_n)^2 F(v\cdot\pi_n)^2\notag\\
&\hskip0.8cm \cdot\tr\big(\big( b(v-s,s\cdot\pi_n) s\cdot\pi_n(a)+b(s-v,v\cdot\pi_n)v\cdot\pi_n(a)\big)^2\big)\, d\nu(v) \, d\nu(s)\\
&\leq\frac{12}{4}\left(\int_{s\in [-1,1]\cap\phi_{\pi_n}^{-1}(V)} F(s\cdot \pi_n)^2  \, d\nu(s)\right)^2
\\
&\hskip3cm
+ \frac{20}{4}\left(\int_{s\in [-1+2n+\pi,1+2n+\pi]\cap\phi_{\pi_n}^{-1}(V)} F(s\cdot\pi_n)^2  \, d\nu(s)\right)^2\\
&\leq (3+5)(1+2\delta)
\end{align*}
because $[-1,1]\cap\phi_{\pi_n}^{-1}(V)\subset (-1/2-\delta, 1/2+\delta)$ and $[-1+2n+\pi,1+2n+\pi]\cap\phi_{\pi_n}^{-1}(V)\subset (-1/2+2n+\pi-\delta, 1/2+2n+\pi+\delta)$ by choice of $V$.

On the other hand, since $F(s\cdot \pi)=1=a_1(s\cdot x_0)$ for $s\in (-1/2,1/2)$ and $b$ is  identically one on $(\supp F_\pi)(\supp F_\pi)^{-1}\times \supp F$,
\begin{align*}
\tr(\Ind\pi(E*E))
&=\frac{1}{4}\int_G\int_G
F(s\cdot\pi)^2 F(v\cdot\pi)^2\notag\\
&\hskip1.5cm \cdot\tr\big(\big( b(v-s,s\cdot\pi) s\cdot\pi(a)+b(s-v,v\cdot\pi)v\cdot\pi(a)\big)^2\big)\, d\nu(v) \, d\nu(s)\\
&\geq \frac{1}{4}\int_{s\in(-1/2,1/2)}\int_{v\in(-1/2,1/2)}\tr\big( (s\cdot\pi(a)+v\cdot\pi(a))^2 \big)\, d\nu(v) \, d\nu(s)\\
&=\frac{1}{4}\int_{s\in(-1/2,1/2)}\int_{v\in(-1/2,1/2)}(a_1(s\cdot x_0)+a_1(v\cdot x_0))^2\, d\nu(v) \, d\nu(s)\\
&=1.
\end{align*}
Now
\begin{align*}
 9>8(1+2\delta)&\geq \liminf_n \tr(\Ind\pi_n(E*E))\\
&\geq
M_L(\Ind\pi,(\Ind\pi_n))\tr(\Ind\pi(E^**E)))\\
&\geq M_L(\Ind\pi,(\Ind\pi_n))
\end{align*}
and hence  $M_L(\Ind\pi,(\Ind\pi_n))\leq 8$.

By \cite[Lemma~A.1]{AaH} there exists a subsequence $(\pi_{n_i})_i$ such that
\[
M_L(\pi,(\pi_{n_i}))
=M_U(\pi,(\pi_{n_i}))=M_U(\pi,(\pi_n)).
\]
But $M_L(\pi,(\pi_{n_i}))\leq 8$ by the argument above, so $M_U(\pi,(\pi_n))\leq 8$.
\end{proof}
Combining Lemma~\ref{lem-ex-upper-bound} with \eqref{eq-ex-lowerbound} we obtain \[M_L(\Ind(\epsilon_{x_0}\otimes\lambda),(\Ind(\epsilon_{x_n}\otimes\epsilon_{y_n})))
=M_U(\Ind(\epsilon_{x_0}\otimes\lambda),(\Ind(\epsilon_{x_n}\otimes\epsilon_{y_n})))=8.\]
%\begin{remark}
%It's fun to briefly consider  what other relative multiplicities can occur for $\Ind(\epsilon_{x_0}\otimes\lambda)$.   For example,
%\begin{align*}
%&M_L(\Ind(\epsilon_{x_0}\otimes\lambda),(\Ind(\epsilon_{x_n}\otimes\lambda)))=2;\\
%&M_L(\Ind(\epsilon_{x_0}\otimes\lambda),(\Ind(\epsilon_{x_0}\otimes\epsilon_{-n})))=3;\\
%&M_L(\Ind(\epsilon_{x_0}\otimes\lambda),(\Ind(\epsilon_{x_0}\otimes\epsilon_{n})))=5;\\
%&M_L(\Ind(\epsilon_{x_0}\otimes\lambda),(\Ind(\epsilon_{x_n}\otimes\epsilon_{n})))=10.
%\end{align*}
%Multiplicity one occurs trivially and multiplicities $4,7$ and $9$ do not occur -- this is intuitively obvious, but how would we prove it if we had to?
%\end{remark}

Clearly, instead of using the numbers $3$ and $5$ in the definition of $B$ we could have chosen  any two positive integers.  By replacing Green's example with one with $k$-times convergence instead of $2$-times convergence, we could arrange to have $k$-times convergence in $\hat A/G$, but the $m_1,m_2,\dots,m_k$ would all be either $3$ or $5$, depending on how  $y_n$ converges to $\lambda$.
Can we build an example for Theorem~\ref{thm-a} where not only $k$ but also $m_1,m_2,\dots,m_k$ take prescribed values?

%%%%%%%%%%%%%%%%%%%%%%%%%%%%%%%%%%%%%%%%%%%%%%%%%%%%%%%%%%%%%%%%%%%%%%%%%%%%%%%%%%%%%%%%%%%%%%%%%%%%%%%%%%%%%%%%%%%
%%%%%%%%%%%%%%%%%%%%%%%%%%%%%%%%%%%%%%%%%%%%%%%%%%%%%%%%%%%%%%%%%%%%%%%%%%%%%%%%%%%%%%%%%%%%%%%%%%%%%%%%%%%%%%%%%%%
%%%%%%%%%%%%%%%%%%%%%%%%%%%%%%%%%%%%%%%%%%%%%%%%%%%%%%%%%%%%%%%%%%%%%%%%%%%%%%%%%%%%%%%%%%%%%%%%%%%%%%%%%%%%%%%%%%%

\appendix\section{Relative multiplicity and tensor products }\label{appendix}

Corollary~\ref{cor-mult-tensor} below is needed in Section~\ref{sec-4}.

\begin{lemma}\label{lem-mult-tensor}
Suppose that $A$ and $B$ are  $C^*$-algebras. Let $\pi\in\hat A$, $(\pi_\alpha)$ a net in $\hat A$, $\sigma\in\hat B$ and $(\sigma_\alpha)$ a net in $\hat B$.
Then, for the spatial tensor product,
\[
M_U(\pi,(\pi_\alpha))M_U(\sigma,(\sigma_\alpha))\geq M_U(\pi\otimes\sigma, (\pi_\alpha\otimes\sigma_\alpha))\geq M_U(\pi,(\pi_\alpha))M_L(\sigma,(\sigma_\alpha)).
\]
\end{lemma}

\begin{proof}
We start by showing that $M_U(\pi\otimes\sigma, (\pi_\alpha\otimes\sigma_\alpha))\geq M_U(\pi,(\pi_\alpha))M_L(\sigma,(\sigma_\alpha))$.  Choose $k,l\in\P$ such that $M_U(\pi,(\pi_\alpha))\geq k$ and $M_L(\sigma,(\sigma_\alpha)\geq l$.  Let $\phi$ and $\psi$ be pure states associated with $\pi$ and $\sigma$, respectively. Since $M_U(\pi,(\pi_\alpha))\geq k$, by \cite[Lemma~5.2(i)]{ASS} there is a subnet $(\pi_{\alpha_\beta})$ of $(\pi_\alpha)$ and, for each $\beta$, an orthonormal subset $\{ \eta_1^\beta,\dots,\eta_k^\beta \}$ of the Hilbert space of $\pi_{\alpha_\beta}$ such that
\[
\phi=\lim_\beta\langle \pi_{\alpha_\beta}(\cdot)\eta_i^\beta\,,\, \eta_i^\beta\rangle\quad (1\leq i\leq k).
\]
Since $M_L(\pi,(\pi_\alpha))\geq l$, by \cite[Lemma~5.2(ii)]{ASS} we may  assume that there is also an orthonormal subset $\{ \xi_1^\beta,\dots,\xi_l^\beta \}$ of the Hilbert space of $\sigma_{\alpha_\beta}$ such that
\[
\psi=\lim_\beta\langle \sigma_{\alpha_\beta}(\cdot)\xi_j^\beta\,,\, \xi_j^\beta\rangle\quad (1\leq j\leq l).
\]
Now consider the orthonormal subset $\{\eta_i^\beta\otimes\xi_j^\beta:1\leq i\leq k, 1\leq j\leq l  \}$ and the subnet $(\pi_{\alpha_\beta}\otimes \sigma_{\alpha_\beta})$ of $(\pi_\alpha\otimes\sigma_\alpha)$.  For $a\otimes b\in A\otimes B$ we have
\begin{align*}
\lim_\beta\langle \pi_{\alpha_\beta}\otimes_{\alpha_\beta}(a\otimes b)\eta_i^\beta\otimes\xi_j^\beta\,,\, \eta_i^\beta\otimes\xi_j^\beta\rangle
&=\lim_\beta\left(\langle \pi_{\alpha_\beta}(a)\eta_i^\beta\,,\, \eta_i^\beta\rangle \langle \sigma_{\alpha_\beta}(bt)\xi_j^\beta\,,\, \xi_j^\beta\rangle \right)\\
&=\phi(a)\psi(b)=\phi\otimes\psi(a\otimes b).
\end{align*}
So by continuity and linearity, $\lim_\beta\langle \pi_{\alpha_\beta}\otimes_{\alpha_\beta}(\cdot)\eta_i^\beta\otimes\xi_j^\beta\,,\, \eta_i^\beta\otimes\xi_j^\beta\rangle=\phi\otimes\psi$, and hence  $M_U(\pi\otimes\sigma, (\pi_\alpha\otimes\sigma_\alpha))\geq kl$ \cite[Lemma~5.2(i)]{ASS}.

It now suffices to assume that $M_U(\pi,(\pi_\alpha))=u_1\in\P$ and $M_U(\sigma,(\sigma_\alpha))=u_2\in\P$ and to show that $M_U(\pi\otimes\sigma, (\pi_\alpha\otimes\sigma_\alpha))\leq u_1u_2$. By \cite[Theorem~2.4]{ASS} there exist $a\in A^+$ and $b\in B^+$ such that $\pi(a)$ and $\sigma(b)$ are non-zero projections and eventually $\tr(\pi_\alpha(a))\leq u_1$ and $\tr(\sigma_\alpha(b))\leq u_2$.  Let $c=a\otimes b$; note $c$ is positive in $A\otimes B$ and that $\pi\otimes\sigma(c)=\pi(a)\otimes\sigma(b)$ is a non-zero projection. Moreover,
$\tr(\pi_\alpha\otimes\sigma_\alpha(c))=\tr(\pi_\alpha(a))\tr(\sigma_\alpha(b))\leq u_1u_2$ eventually.  So $M_U(\pi\otimes\sigma, (\pi_\alpha\otimes\sigma_\alpha))\leq u_1u_2$ by \cite[Theorem~2.4]{ASS}.
\end{proof}

%%%%%%%%%%%%%%%%%%%%%%%%%%%%%%%%%%%%%%%%%%%%%%%%%%%%%%%%%%%%%%%%%%%%%%%%%%%%%%%%%%%%%%%%%%%%%%%%%%%%%%%%%%%%%%%%%%%%%%

\begin{cor}\label{cor-mult-tensor}
Suppose that $A$ and $B$ are  $C^*$-algebras. Let $\pi\in\hat A$, $(\pi_\alpha)$ a net in $\hat A$, $\sigma\in\hat B$ and $(\sigma_\alpha)$ a net in $\hat B$.
Also suppose that $M_L(\pi,(\pi_\alpha))=M_U(\pi,(\pi_\alpha))$ and $M_L(\sigma,(\sigma_\alpha))=M_U(\sigma,(\sigma_\alpha))$.
Then, for the spatial tensor product,
\[
M_L(\pi\otimes\sigma, (\pi_\alpha\otimes\sigma_\alpha))=M_U(\pi\otimes\sigma, (\pi_\alpha\otimes\sigma_\alpha))=M(\pi,(\pi_\alpha))M(\sigma,(\sigma_\alpha).
\]
\end{cor}
\begin{proof}By \cite[Proposition~2.3]{AS} there exists a subnet $\pi_{\alpha_\beta}\otimes \sigma_{\alpha_\beta}$ of $\pi_\alpha\otimes\sigma_\alpha$ such that
\[M_L(\pi\otimes\sigma, (\pi_\alpha\otimes\sigma_\alpha))=M(\pi\otimes\sigma, (\pi_{\alpha_\beta}\otimes\sigma_{\alpha_\beta})).\]
Thus
\begin{align*}
M_L(\pi\otimes\sigma, (\pi_\alpha\otimes\sigma_\alpha))
&=M_U(\pi\otimes\sigma, (\pi_{\alpha_\beta}\otimes\sigma_{\alpha_\beta})) \quad\text{(by Lemma~\ref{lem-mult-tensor})}\\
&=M_U(\pi,(\pi_{\alpha_\beta}))M_U(\sigma, (\sigma_{\alpha_\beta}))\\
&=M_U(\pi,(\pi_{\alpha}))M_U(\sigma, (\sigma_{\alpha}))\\
&=M_U(\pi\otimes\sigma, (\pi_\alpha\otimes\sigma_\alpha)) \quad\text{(by Lemma~\ref{lem-mult-tensor})},
\end{align*}
so the result follows.
\end{proof}

%%%%%%%%%%%%%%%%%%%%%%%%%%%%%%%%%%%%%%%%%%%%%%%%%%%%%%%%%%%%%%%%%%%%%%%%%%%%%%%%%%%%%%%%%%%%%%%%%%%%%%%%%%%%%%%%%%%%%%
%%%%%%%%%%%%%%%%%%%%%%%%%%%%%%%%%%%%%%%%%%%%%%%%%%%%%%%%%%%%%%%%%%%%%%%%%%%%%%%%%%%%%%%%%%%%%%%%%%%%%%%%%%%%%%%%%%%%

\end{document}